\documentstyle[amscd]{amsart}
\newsymbol\varkappa 207B
\def\dj{d\kern-.30em\raise1.25ex\vbox{\hrule width .3em height .03em}}
\def\Dj{D\kern-.75em\raise0.75ex\vbox{\hrule width .3em height .03em}
\kern.03em}
\newsymbol\restriction 1316
\newsymbol\between 1347
\makeatletter
\renewcommand{\subsection}{\@startsection{subsection}{2}{\z@}%
{\baselineskip}{0.5\baselineskip}{\defaultfont\bf}}
\makeatother
\newcommand{\rig}{\wp}
\def\restr{{\restriction}}
\def\Der{\mbox{Der}}
\def\ellG{\ell}
\def\bla#1{$(${\it #1\/{}}$)$}
\def\vrkpa{\varkappa}
\newcommand{\grten}{\mathbin{\widehat{\otimes}}}
\newcommand{\inv}{i\!\hspace{0.8pt}n\!\hspace{0.6pt}v}
\newcommand{\Uni}{\between}
\newcommand{\lie}{\mbox{\family{euf}\shape{n}\selectfont lie}}
\newcommand{\Sum}{\displaystyle{\sum}}
\newcommand{\hor}{\mbox{\family{euf}\shape{n}\selectfont hor}}
\newcommand{\ver}{\mbox{\family{euf}\shape{n}\selectfont ver}}
\newcommand{\con}{\mbox{\family{euf}\shape{n}\selectfont con}}

\newcommand{\vh}{\mbox{\family{euf}\shape{n}\selectfont vh}}

\newcommand{\ad}{\mbox{\shape{n}\selectfont ad}}
\newcommand{\id}{\mbox{\shape{n}\selectfont id}}
\newcommand{\im}{\mbox{\shape{n}\selectfont im}}
\newcommand{\adj}{\varpi}
\newcommand{\k}{\kappa}
\newcommand{\e}{\epsilon}

\newcommand{\lin}{\mbox{lin}}
\def\sstar{{\raise0.2ex\hbox{$\scriptstyle\star$}}}
\renewcommand{\L}{\cal{L}}
\newtheorem{thm}{Theorem}[section]
\newtheorem{lem}[thm]{Lemma}
\newtheorem{pro}[thm]{Proposition}
\theoremstyle{definition}
\newtheorem{defn}{Definition}[section]
\numberwithin{equation}{section}
\begin{document}
\title[QUANTUM PRINCIPAL BUNDLES]{GEOMETRY OF QUANTUM PRINCIPAL BUNDLES I}
\author{Mi\'co \Dj ur\Dj evi\'c}
\address{Instituto de Matematicas, UNAM,
Area de la Investigacion Cientifica,
Circuito Exterior,
Ciudad Universitaria,
M\'exico DF, CP 04510, MEXICO\newline
\mbox{}\newline
\indent {\it Written In}\newline
\indent Faculty of Physics, University of Belgrade, Pbox 550,
Studentski Trg 12, 11001 Beograd, YUGOSLAVIA\newline
\indent {\it
Summer 1991--Abstract Theory\newline
\indent Summer 1992--Examples}
}
\maketitle
\begin{abstract}
A theory of principal bundles possessing quantum
structure  groups  and classical base manifolds is
presented. Structural analysis of such quantum  principal  bundles
is performed. A differential calculus  is  constructed,  combining
differential forms  on  the  base  manifold  with  an  appropriate
differential calculus on the structure  quantum  group.  Relations
between the calculus on the group and the calculus on  the  bundle
are investigated. A concept of (pseudo)tensoriality is formulated.
The  formalism  of  connections  is  developed.   In   particular,
operators  of  horizontal  projection,  covariant  derivative  and
curvature  are  constructed  and  analyzed.   Generalizations   of
the first structure equation and of the Bianchi  identity  are  found.
Illustrative examples are presented.
\end{abstract}
\renewcommand{\thepage}{}
\tableofcontents
\filbreak
\section{Introduction}
In diversity of mathematical concepts and theories a fundamental
role is played by those giving a unified  treatment  of  different
and at a first sight mutually independent circles of problems.

As far as classical differential geometry  is  concerned,  such  a
fundamental role is given to the theory of principal bundles \cite{KN}.
Various basic concepts of theoretical physics are  also  naturally
expressible  in  the  language  of  principal  bundles.  Classical
gauge theory is a paradigmic example.

In this work a quantum generalization of the theory  of  principal
bundles will be presented. All  constructions  and  considerations
will be performed within a conceptual framework of  noncommutative
differential geometry \cite{C},\cite{K}.

The generalization will be twofold. First of all,  quantum  groups
will  play  the  role  of  structure  groups.   Secondly,
appropriate quantum spaces will play the role of base manifolds.

This paper is devoted to the study of  quantum  principal  bundles
over classical smooth manifolds.

The paper is organized as follows.

\smallskip
\renewcommand{\thepage}{\arabic{page}}
Section 2 begins with a definition of quantum  principal  bundles.
For technical reasons, it will be assumed that a base manifold $M$
is compact. Concerning a structure quantum group $G$, it will be a
compact matrix quantum group (pseudogroup), in the sense of \cite{W2}.

We shall prove that, as a consequence of an  inherent  geometrical
inhomogeneity  of  quantum   groups,   there   exists   a   natural
correspondence between quantum principal  bundles,  and  classical
principal bundles over the same manifold $M$, with  the  structure
group $G_{\!cl}$ consisting  of  ``classical  points''  of  $G$.  Informally
speaking, if we start from a quantum principal bundle $P$ then the
corresponding  classical  principal   bundle   $P_{\!cl}$   consists
precisely of ``classical points'' of $P$. Conversely, starting  from
a  $G_{\!cl}$-bundle  $P_{\!cl}$,  the  bundle  $P$  can  be  recovered
applying  a  variant  of  the  classical  procedure  of  extending
structure groups.

Section 3 is devoted to the  study  of  differential  calculus  on
quantum  principal  bundles.  As  first,  general  properties  for
differential  calculus  on  $P$  will  be  formulated,   including
relations with differential structures over $M$ and $G$. The main
idea  is  that  local  trivializations  of  the   bundle   locally
trivialize the calculus, too.

A differential  calculus  over  $M$  will  be  the  standard  one,
specified by differential forms. A differential  calculus  on  the
structure quantum  group  $G$  will  be  based  on  the
{\it universal
envelope} of  an  appropriate  first-order  differential  calculus
$\Gamma$.  This  universal  envelope  can  be  constructed  by
applying an extended bimodule technique \cite{W1,W3}. As we shall see,
the mentioned local triviality property of  the  calculus  on  the
bundle implies certain restrictions on the calculus  $\Gamma$.
Informally speaking, $\Gamma$ should be  compatible  with  all
possible  ``transition  functions''  for  $P$.  Motivated  by   this
observation, we shall introduce a notion  of  {\it admissibility\/}  to
distinguish first-order differential structures on $G$  for  which
the mentioned compatibility holds.

The next theme of Section 3 is a construction of the  calculus  on
$P$,  starting  from  differential  forms  on  $M$  and  a   given
admissible first-order calculus $\Gamma$ over $G$. As a result
we obtain a graded  differential  algebra  $\Omega(P,\Gamma)$,
representing the calculus on $P$. We shall prove the uniqueness
of this algebra.

After this, various properties of $\Omega(P,\Gamma)$  will  be
studied (the existence of $*$-structures, the right covariance and
the  existence  of  the  graded-differential  extension   of   the
dualized right action of  $G$  on  $P$).  These  properties  are
closely related to similar  properties  of  $\Gamma$.  On  the
other hand, independently of  the  choice  of  $\Gamma$  there
exists a natural left coaction of $G$  on  $\Omega(P,\Gamma)$,
becoming trivial in the classical case.

In Section 3 the structure of  admissible  calculi is studied,
too. In particular,  left-covariant  admissible  calculi  are
characterized in terms of the corresponding right  ideals  in  the
algebra $\cal{A}$ of ``polynomial functions'' on $G$. It  turns  out
that  there  exists  the  ``simplest''   left-covariant   admissible
calculus (which is automatically bicovariant and *-covariant).

Finally, at the  end  of  Section  3  we  introduce and
briefly analyze analogs  of horizontal and verticalized
differential forms on the bundle.

The study of connections on quantum principal bundles is the  main
topic of Sections 4 and 5.  Through these sections we shall assume
that  $\Gamma$  is  the  simplest  left-covariant  admissible
calculus.

In  Section  4  we  shall   first   generalize
the classical  concept  of  (pseudo)tensorialilty.  Together  with
certain considerations performed in Section 3 this will anable  us
to introduce connection forms, in analogy with classical geometry.
We then pass to the study of local representations of  connections,
in terms of gauge potentials.

Further, we shall prove that  each  connection  on  $P$  admits  a
decomposition into a ``classical connection'', interpretable  as  an
ordinary connection  on  $P_{\!cl}$,  and  an  appropriate  ``purely
quantum'' tensorial 1-form.

Each connection decomposes the algebra $\Omega(P,\Gamma)$  into
a tensor product of spaces of horizontal forms and  left-invariant
forms on $G$.  With  the  help  of  this  decomposition  we  shall
introduce the horizontal  projection operator. This will
enable us to  define  the  analogs  of covariant  derivative  and
curvature operators, which will be studied in Section 5.  In  particular,
we shall analyze local representations of covariant derivative and
curvature, and find counterparts of the first  Structure  Equation
and the Bianchi identity.

In Section 5 some concrete examples are worked out. Considerations
are mainly confined to specific  properties  of  the  calculus  on
structure quantum group $G$, and to the presentation  of  ``quantum
phenomena'' appearing at the level  of  connections.  A  particular
care is  devoted  to  the  example  with  quantum  $SU(2)$  group.
Finally, we shall briefly discuss  a  possible  formulation  of  a
``gauge theory'' in the framework of quantum principal bundles.

The paper ends with three technical appendices.
In Appendix A relevant properties of the set $G_{\!cl}$ of classical
points of $G$ are collected. Some concrete examples are computed.

In the  second  appendix  properties  of  universal  envelopes  of
first-order differential structures are analyzed in details. It is
important to mention that, in  the  general  case,  the  universal
envelope of a bicovariant first-order calculus does  not  coincide
with the exterior algebra constructed in \cite{W4},  although  in  the
case of ordinary Lie groups (and ordinary  1-forms  on  them)  two
structures coincide. We shall see that,  quite  generally, the
universal envelope coincides with the graded-differential  algebra
constructed  by  applying  the   mentioned   extended   bimodule
technique. A reason for our choice of higher-order calculus on $G$
lies in the conceptual simplicity of the universal calculus, which
is independent of the group structure on $G$ (in contrast  to  the
exterior  algebra  construction).   Because   of   this,   similar
considerations can be  applied  to  more  general  fiberings,  for
example of  the  type  of  associated  bundles  where  fibers  are
diffeomorphic to an arbitrary quantum space. On the other hand, we
are able  to  consider  examples  in  which  $\Gamma$  is  not
bicovariant.

We shall also prove that $\Omega(P,\Gamma)$ can be  understood
as the universal envelope over its first-order part.

In  Appendix  C  some  properties  of  already  mentioned  minimal
admissible first-order calculi are collected.

Concerning the notation of quantum group entities, we shall follow
\cite{W2}. A quantum group $G$ will be represented as a  pair  $G=(A,u)$,
where $A$ is the $C^*$-algebra of ``continuous  functions''  on  the
space $G$ and $u\in M_n(A)$ is the  matrix  determining  the
group  structure.   The   $*$-algebra   representing   ``polynomial
functions'' on $G$ will be denoted by $\cal  A$.  This
$*$-algebra  is
generated by entries of $u$. The comultiplication, the counit  and
the antipode will be denoted by $\phi$, $\e$ and $\k$ respectively.

We shall write symbolically
$$\phi(a)=a^{(1)}\otimes  a^{(2)}$$
for each $a\in \cal A$. Similarly, the symbol $a^{(1)}\otimes \dots
a^{(n)}$ denotes the result of a $(n-1)$-fold comultiplication
of $a\in\cal A$ (due to the  coassociativity  property  of  $\phi$
this is independent of the way  in  which  comultiplications  are
performed).

We shall denote by $\ad\colon\cal A\rightarrow\cal A\otimes\cal  A$  the
adjoint action of $G$ on itself. Explicitly, this map is given by
$$\ad(a)=a^{(2)}\otimes \k(a^{(1)})a^{(3)}.$$

If $M$ is  a  smooth  manifold  we  shall  denote  by  $S(M)$  the
$*$-algebra of complex smooth functions on $M$. Similarly, $S_{\!c}(M)$
will be  the  $*$-algebra  consisting  of  smooth  functions having a
compact support.
\section{Structure of Quantum Principal Bundles}

Let us consider a compact matrix quantum group $G$. Let $M$ be
a compact smooth manifold.
\begin{defn}\label{def:bun}
A {\it  (quantum)  principal  $G$-bundle}  over  $M$  is  a  triplet
of the form $P=(\cal B,i,F)$ where  $\cal B$ is a (unital)
$*$-algebra, $i:S(M)\rightarrow\cal B$ is a unital linear  map
and $F:\cal B\rightarrow\cal B\otimes\cal A$ is a linear map  such
that for each $x\in M$ there exists an  open  set  $U\subseteq  M$
containing $x$ and a homomorphism $\pi_U\colon\cal B\rightarrow
S(U)\otimes\cal A$ such that the following properties hold:
\par\smallskip
\bla{qpb1} We have
$$\pi_Ui(f)=(f\restr_U)\otimes 1$$ for each $f\in S(M)$. \par\smallskip
\bla{qpb2} If $q=i(\varphi)b$ where
$\varphi\in S_{\!c}(U)$ then  $\pi_U(q)=0$  implies
$q=0$.\par\smallskip
\bla{qpb3} We have $$(\id\otimes\phi)\pi_U=(\pi_U\otimes \id)F\qquad
\pi_U(\cal B)\supseteq S_{\!c}(U)\otimes \cal A.$$
\end{defn}

A motivation for this definition comes from classical differential
geometry. The  map  $i\colon S(M)\rightarrow\cal  B$  is
interpretable   as   the
``dualized projection'' of the bundle $P$ on its base $M$.
The map $F$  plays the role of a dualized right action of $G$
on $P$.  Finally,  maps
$\pi_U$ are dualized local trivializations of the bundle.

Let $P=(\cal B,i,F)$ be a principal $G$-bundle over $M$.
\begin{defn} A {\it  local  trivialization}  for  $P$  is  a  pair
$(U,\pi_U)$ consisting of a non\-empty open  set  $U\subseteq  M$
and a $*$-homomorphism $\pi_U\colon \cal B\rightarrow S(U)\otimes\cal A$
such that properties listed in the  previous  definition  hold.  A
{\it   trivialization   system}    for    $P$    is    a    family
$\tau=(\pi_U)_{U\in\cal U}$, where $\cal U$ is a finite open cover
of $M$ and for each $U\in \cal U$ the pair $(U,\pi_U)$ is a  local
trivialization for $P$.
\end{defn}

Let $\tau=(\pi_U)_{U\in\cal U}$ be  a  trivialization  system  for
$P$.

\begin{lem}\label{lem:21}
The family $\tau$ distinguishes elements of $\cal B$.
\end{lem}
\begin{pf} Let us consider a partition of unity
$\varpi=(\varphi_U)_{U\in\cal U}$ for $\cal{U}$. In other words
$\varphi_U\in S_{\!c}(U)$ and
$$ \sum_{U\in\cal U}\varphi_U=1_M. $$
According to Definition~\ref{def:bun} if
$b$ belongs to the intesection of kernels of maps $\pi_U$ then
$i(\varphi_Ub)=0$, and hence $\varphi_Ub=0$,
for each $U\in\cal U$. Summing over $\cal U$  we
conclude that $b=0$.
\end{pf}
\begin{lem}\label{lem:22} \bla{i} The map $i\colon S(M)\rightarrow\cal B$ is
a $*$-monomorphism.

\smallskip
\bla{ii} The image $i\bigl(S(M)\bigr)$ is contained
in the centre of $\cal B$.
\end{lem}
\begin{pf} The following equalities hold
\begin{gather*}
\pi_U\bigl(i(fg)-i(f)i(g)\bigr)=(fg\restr_U)\otimes1-(f\restr_U)
(g\restr_U)\otimes 1=0\\
\pi_U\bigl(i(f^*)\bigr)-\pi_U\bigl(i(f)^*\bigr)=(f^*\restr_U)\otimes 1
-(f\restr_U)^*\otimes 1=0\\
\pi_U\bigl(i(f)b-bi(f)\bigr)
=\bigl((f\restr_U)\otimes 1\bigr)\pi_U(b)-\pi_U(b)
\bigl((f\restr_U)\otimes 1\bigr)=0.
\end{gather*}
Using   Lemma~\ref{lem:21}   we   conclude   that   $i$    is    a
$*$-homomorphism and  that  \bla{ii}  holds.  If  $f\in  \ker(i)$  then
$f\restr_U=0$ for each $U\in \cal U$ and hence $f=0$.
\end{pf}
\begin{lem}\label{lem:23} \bla{i} The map $F$ is a unital *-homomorphism.

\smallskip
\bla{ii} The following identities hold
\begin{gather}
(F\otimes \id)F=(\id\otimes\phi)F\label{coas} \\
(\id\otimes \e)F =\id.\label{ide}
\end{gather}

\bla{iii} An element $b\in \cal B$ belongs to $i\bigl(S(M)\bigr)$ iff
\begin{equation}
F(b)=b\otimes 1.\label{orb}
\end{equation}

In other words $F$ defines a right  action  of  $G$  on  $P$.  The
corresponding ``orbit space'' coincides with the base manifold $M$.
\end{lem}
\begin{pf} According to Definition~\ref{def:bun},
\begin{multline*}
(\pi_U\otimes \id)F(b^*)=(\id\otimes\phi)\pi_U(b^*)=
\bigl((\id\otimes\phi)\pi_U(b)\bigr)^*\\
=\bigl((\pi_U\otimes \id)F(b)\bigr)^*=
(\pi_U\otimes \id)\left(F(b)\right)^*
\end{multline*}
as well as
\begin{equation*}
\begin{split}
(\pi_U\otimes \id)F(bq)&=(\id\otimes\phi)\pi_U(bq)=
(\id\otimes\phi)\bigl(\pi_U(b)\pi_U(q)\bigr)\\
&=\bigl((\id\otimes\phi)\pi_U(b)\bigr)\bigl((\id\otimes\phi)\pi_U(q)
\bigr)\\
&=\bigl((\pi_U\otimes \id)F(b)\bigr)\bigl((\pi_U\otimes
\id)F(q)\bigr)\\
&=(\pi_U\otimes \id)\bigl(F(b)F(q)\bigr)
\end{split}
\end{equation*}
for each $U\in \cal U$. Hence, $F$ is a *-homomorphism.  Equations
\eqref{coas}--\eqref{ide} as well as the identity
$$ Fi(f)=i(f)\otimes 1$$
can be checked in a similar way.

Let us assume that $F(b)=b\otimes 1$. We have then
$$  (\pi_U\otimes  \id)F\bigl(i(\varphi_U)b\bigr)
=\pi_U\bigl(i(\varphi_U)b\bigr)\otimes
1=(\id\otimes\phi)\pi_U(i(\varphi_U)b),$$
where $\left(\varphi_U\right)_{U\in \cal  U}$  is  a  partition  of
unity for $\cal U$.

Acting by $\id\otimes \e\otimes \id$ on  the  second
equality we obtain
$$ \pi_U\bigl(i(\varphi_U)b\bigr)=\bigl[(\id\otimes \e)\pi_U\bigl(
i(\varphi_U)b\bigr)\bigr]\otimes1.$$
It follows that
$$i(\varphi_U)b=i(\eta_U),$$
where $\eta_U=(\id\otimes\e)\pi_U\bigl(i(\varphi_U)b\bigr)$.
Summing over $U$'s we obtain
$$ b=i\biggl(\sum_{U\in\cal U}\eta_U\biggr).$$
Finally, the unitality of $F$ directly follows from
\bla{iii}  and  from the unitality of $i$.
\end{pf}

We pass to the study of internal structure of quantum principal
bundles, in terms of the corresponding ``$G$-cocycles''.

For a given open cover  $\cal  U$  of  $M$,  we  shall  denote  by
$N^k(\cal U)$ the set of all k-tuples  $(U_1,\dots,  U_k)$,  where
$U_i\in\cal U$ are such that $U_1\cap\dots\cap U_k\neq\emptyset.$

\begin{defn} Let $\cal U$  be  a  finite  open  cover  of  $M$.  A
({\it smooth, quantum\/}) {\it $G$-cocycle} over $(M,\cal U)$ is a system
$\cal{C}=\Bigl\{\psi_{U\!V}\mid(U,V)\in N^2(\cal U)\Bigr\}$ of non-trivial
$S(U\cap V)$-linear
*-homomorphisms
$\psi_{U\!V}\colon S(U\cap V)\otimes\cal A\rightarrow S(U\cap V)\otimes
\cal A$ such that

\smallskip
\bla{i} The diagram
\begin{equation}\label{27}
\begin{CD}
S(U\cap V)\otimes\cal A @>{\mbox{$\psi_{U\!V}$}}>>
S(U\cap V)\otimes\cal A\\
@V{\mbox{$\id\otimes\phi$}}VV    @VV{\mbox{$\id\otimes\phi$}}V\\
S(U\cap V)\otimes\cal A\otimes\cal A  @>>{\mbox{$\psi_{U\!V}\otimes \id$}}>
S(U\cap V)\otimes\cal A\otimes\cal A
\end{CD}
\end{equation}
is commutative.

\smallskip
\bla{ii} We have
\begin{equation}\label{28}
\psi_{U\!V}\bigl[\psi_{V\!W}(\varphi)\bigr]=\psi_{U\!W}(\varphi),
\end{equation}
for each $(U,V,W)\in N^3(\cal U)$ and $\varphi\in S_{\!c}(U\cap V\cap W)$.
\end{defn}

Let us  observe  that  $S(U\cap  V)$-linearity  property  of  maps
$\psi_{U\!V}$ implies
$$ \psi_{U\!V}\left[S_{\!c}(W)\otimes \cal A\right]\subseteq
S_{\!c}(W)\otimes \cal A$$
for each (nonempty) open set $W\subseteq U\cap V$.  Furthermore,  maps
$\psi_{U\!V}$ are completely determined  by  their  restrictions  on
$S_{\!c}(U\cap V).$

The following proposition completely describes  $G$-cocycles.  Let
$G_{\!cl}$ be {\it the classical part} of $G$ (Appendix A). This  is
a classical group (a ``subgroup'' of $G$)
consisting of {\it points} of $G$ (formally *-characters on $\cal  A$).

\begin{pro}\label{pro:24} For each  $G$-cocycle
$\cal{C}=\Bigl\{\psi_{U\!V}\mid(U,V)\in
N^2(\cal  U)\Bigr\}$  there  exists  the  unique
collection of smooth maps $g_{V\!U}\colon U\cap V\rightarrow  G_{\!cl}$ such
that
\begin{equation}\label{29}
\psi_{U\!V}(\varphi\otimes
a)\restr_x=\varphi g_{V\!U}(x)(a^{(1)})\otimes a^{(2)}.
\end{equation}
Maps $g_{U\!V}$ form a classical $G_{\!cl}$-cocycle over $(M,\cal U).$

Conversely, if $g_{U\!V}$ form  a  classical  $G_{\!cl}$-cocycle  then
formula \eqref{29} determines a
quantum $G$-cocycle over $(M,\cal U)$.
\end{pro}
\begin{pf}  Let  $\cal{C}=\Bigl\{\psi_{U\!V}\mid
(U,V)\in   N^2(\cal   U)\Bigr\}$   be   a
$G$-cocycle. For each $(U,V)\in N^2(\cal{U})$ let us define a map $\mu_{V\!U}
\colon \cal A\rightarrow S(U\cap V)$ by
\begin{equation}\label{210}
\mu_{V\!U}(a)=(\id\otimes \e)\psi_{U\!V}(1\otimes a).
\end{equation}
Acting by  $\id\otimes  \e\otimes  \id$  on  both  wings  of  diagram
\eqref{27} we obtain
\begin{equation}\label{213}
\psi_{U\!V}(\varphi\otimes        a)=\varphi\mu_{V\!U}(a^{(1)})\otimes
a^{(2)}.
\end{equation}
Maps $\mu_{V\!U}$ are unital
*-homomorphisms.  Equivalently,  they  can  be
naturally
understood as smooth maps $g_{V\!U}\colon U\cap V\rightarrow  G_{\!cl},$  by
exchanging the order of arguments:
$$ \bigl[\mu_{V\!U}(a)\bigr](x)=\bigl[g_{V\!U}(x)\bigr](a). $$
We see that \eqref{29} holds.  Now  acting  by  $\id\otimes  \e$  on
\eqref{28}, using \eqref{29} and the definition
of the product in $G_{\!cl}$ we conclude that
\begin{equation}
g_{U\!V}g_{V\!W}=g_{U\!W}
\end{equation}
for each $(U,V,W)\in N^3(\cal U).$ In other words, maps
$g_{V\!U}$ form a classical $G_{\!cl}$-cocycle over  $(M,\cal{U})$.
The  second part  of  the  proposition   easily   follows
from the coassociativity of $\phi$ and the definition  of
the  product  in $G_{\!cl}$.
\end{pf}

Property \eqref{29} implies that maps $\psi_{U\!V}$  are  bijective.
Indeed, the inverse is explicitly given by
\begin{equation}
\psi_{U\!V}^{-1}(\varphi\otimes a)\restr_x=\varphi g_{U\!V}(x)
(a^{(1)})\otimes a^{(2)}.
\end{equation}
In particular, \eqref{28} implies
$$
\psi_{U\!U}(f)=f \qquad
\psi_{U\!V}^{-1}=\psi_{V\!U}.
$$

We see that $G$-cocycles are  in  a  natural  correspondence  with
$G_{\!cl}$-cocycles. On the other hand, $G_{\!cl}$-cocycles are  in  a
natural correspondence with classical  principal  $G_{\!cl}$-bundles
over $M$ (endowed with a trivialization system).

A similar correspondence holds between  quantum  $G$-cocycles  and
quantum principal $G$-bundles. Let $P=(\cal B,i,F)$ be  a  quantum
principal $G$-bundle over $M$. For a given (nonempty) open set
$V\subseteq M$ let us denote by  $I_V$  the  lineal  in  $\cal  B$
consisting of elements of the form  $q=i(\varphi)b$,  where  $b\in
\cal B$ and $\varphi\in S_{\!c}(V)$.  Lemma~\ref{lem:22}  \bla{ii}
implies that $I_V$ is a (two-sided) *-ideal in $\cal B$.

Let $(U,\pi_U) $ be a local trivialization of $P$.
The following lemma is a direct consequence of
properties listed in Definition~\ref{def:bun}.
\begin{lem}\label{lem:25} Let $V\subseteq U$ be a nonempty open set. Then
$$ \pi_U(I_V)\subseteq S_{\!c}(V)\otimes \cal A$$
and the restriction $(\pi_U\restr I_V)\colon I_V\rightarrow
S_{\!c}(V)\otimes\cal A$ is a *-isomorphism. \qed
\end{lem}

Let  $\psi_U\colon S_{\!c}(U)\otimes\cal  A\rightarrow \cal{B}$  be  a
*-monomorphism defined by
\begin{equation}\label{214}
\psi_U=\left(\pi_U\restr I_U\right)^{-1}.
\end{equation}
Evidently, the diagram
\begin{equation}\label{215}
\begin{CD}
S_{\!c}(U)\otimes\cal A @>{\mbox{$\psi_U$}}>> \cal B\\
@V{\mbox{$\id\otimes\phi$}}VV         @VV{\mbox{$F$}}V\\
S_{\!c}(U)\otimes\cal A\otimes\cal A @>>{\mbox{$\psi_U\otimes \id$}}>
\cal B\otimes \cal A
\end{CD}
\end{equation}
is commutative.

Let us consider a trivialization system
$\tau=(\pi_U)_{U\in\cal U}$ for $P$.
\begin{lem}\label{lem:26}
There exists the unique $G$-cocycle
$\cal C_\tau=\Bigl\{\psi_{U\!V}\mid (U,V)\in N^2(\cal U)\Bigr\}$
satisfying
\begin{equation}
\psi_{U\!V}(q)=\pi_U\psi_V(q)
\end{equation}
for each $(U,V)\in N^2(\cal U)$ and
$q\in S_{\!c}(U\cap V)\otimes \cal A.$
\end{lem}
\begin{pf} The above formula defines maps $\psi_{U\!V}$ on algebras
$S_{\!c}(U\cap  V)\otimes  \cal  A.$  These maps are $S(U\cap
V)$-linear. Because of this it  is  possible  to  extend  them
uniquely to *-homomorphisms $\psi_{U\!V}\colon S(U\cap V)\otimes
\cal{A}\rightarrow
S(U\cap V)\otimes\cal A$. Covariance property \eqref{27}  follows
from  \eqref{215}.  Cocycle  condition  \eqref{28}  is  a   direct
consequence of the definition of maps $\psi_{U\!V}$.
\end{pf}

Let us consider an arbitrary $G$-cocycle
$\cal C_\tau=\Bigl\{\psi_{U\!V}\mid(U,V)\in N^2(\cal U)
\Bigr\}$, and  let us define a *-algebra $\cal T$ as a direct sum
$$ \cal T=\sideset{}{^\oplus}\sum_{U\in\cal U} S(U)\otimes
\cal A. $$
Let $\widetilde{\cal B}$ be a set consisting of elements
$b\in\cal T$ satisfying
\begin{equation}\label{222}
(\sideset{_U}{_{U\cap V}}\mid \otimes \id)p_U(b)=
\psi_{U\!V}(\sideset{_V}{_{U\cap V}}\mid \otimes \id)p_V(b)
\end{equation}
for each $(U,V)\in N^2(\cal U)$, where $p_U$ and
$\sideset{_U}{_{U\cap V}}\mid$ are the corresponding coordinate
projections and restriction maps.

All maps  figuring  in  \eqref{222}  are  *-homomorphisms.  Hence,
$\widetilde{\cal B}$ is a *-subalgebra of $\cal T$. The formula
\begin{equation}\label{223}
(p_U\otimes \id)F_{\cal T}=(\id\otimes \phi)p_U
\end{equation}
determines a *-homomorphims $F_{\cal T}\colon \cal T\rightarrow
\cal T\otimes\cal{A}$.
Diagram \eqref{27} implies that $\widetilde{\cal  B}$  is
$F_{\cal T}$-invariant, in the sense that $F_{\cal T}(
\widetilde{\cal B})\subseteq \widetilde{\cal B}\otimes
\cal A$. Let $\widetilde{F}\colon \widetilde{\cal B}
\rightarrow \widetilde{\cal B}\otimes \cal A$ be the corresponding
restriction map. The formula
\begin{equation}
p_U\tilde{\imath}(f)=(f\restr_U)\otimes 1
\end{equation}
defines a *-homomorphism $\tilde{\imath}\colon S(M)\rightarrow
\widetilde{\cal B}$. Let $\pi_U\colon \widetilde{\cal B}
\rightarrow S(U)\otimes \cal A$ be
the restrictions of coordinate projection maps.
\begin{pro}\label{pro:27} The triplet $\widetilde{P}=
(\widetilde{\cal  B},  \tilde{\imath},  \widetilde{F})$  is a
principal $G$-bundle over $M$. The family
$\tau=(\pi_U)_{U\in\cal  U}$  is  a  trivialization   system   for
$\widetilde{P}$. The corresponding $G$-cocycle  coincides  with the
initial one. In other words $\cal{C}=\cal{C}_\tau$.\qed
\end{pro}
The above proposition directly follows from the construction of
$\widetilde{P}$.
Let  $P=(\cal  B,i,F)$  be  a  principal  $G$-bundle   over   $M$,
with a trivialization system $\tau$.
\begin{lem}\label{lem:28}
The following identities hold
\begin{equation}
(\sideset{_U}{_{U\cap V}}\mid \otimes \id)\pi_U(b)=
\psi_{U\!V}(\sideset{_V}{_{U\cap V}}\mid \otimes \id)\pi_V(b),
\end{equation}
where $\psi_{U\!V}$ are transition functions from $\cal{C}_\tau$.
\end{lem}
\begin{pf}
It is sufficient to check that above equalities hold on elements
of the form $q=i(\varphi)b$, where $\varphi\in S_{\!c}(U\cap V)$.
However, this is equivalent to
$$ \psi_{U\!V}\pi_V(q)=\pi_U(q)$$
which is the definition of $\psi_{U\!V}.$
\end{pf}
\begin{pro}\label{pro:29} Let $\widetilde{P}=(\widetilde{\cal B},
\tilde{\imath},\widetilde{F})$ be  a  principal  $G$-bundle  constructed
from  the  $G$-cocycle  $\cal  C_\tau$.  Then  the  *-homomorphism
$j_\tau\colon \cal B\rightarrow\cal T$ defined by
\begin{equation}\label{225}
p_Uj_\tau=\pi_U
\end{equation}
isomorphically maps $\cal B$ onto $\widetilde{\cal B}$. Moreover,
the following equalities hold
\begin{align}
\widetilde{F}j_\tau&=(j_\tau\otimes \id)F\label{227}\\
j_\tau i&=\tilde{\imath}.\label{226}
\end{align}
\end{pro}
\begin{pf}
According to Lemma~\ref{lem:28} we have $j_\tau(\cal B)\subseteq
\widetilde{\cal B}$. Further
$$p_Uj_\tau i(\varphi)=(\varphi\restr_U)\otimes1
=p_U\tilde{\imath}(\varphi),$$
for each $\varphi\in S(M)$  and  $U\in\cal  U$.  Thus  \eqref{226}
holds. Together with \eqref{225} this implies
\begin{equation}\label{228}
j_\tau\psi_U=\widetilde{\psi}_U
\end{equation}
where $\widetilde{\psi}_U$ are the  corresponding  right inverses  for
maps $\widetilde{\pi}_U=p_U\restr\widetilde{\cal B}$.

The    map    $j_\tau$    is    surjective,     because     spaces
$\widetilde{\psi}_U\bigl[S_{\!c}(U)\otimes\cal A\bigr]$ linearly span
$\widetilde{\cal B}$. Injectivity of $j_\tau$ is a consequence  of
Lemma~\ref{lem:21}. Hence, $j_\tau\colon\cal{B}\leftrightarrow
\widetilde{\cal{B}}$.

Finally, we have
$$(p_Uj_\tau\otimes \id)F=(\pi_U\otimes \id)F=(\id\otimes\phi)\pi_U=
(\id\otimes\phi)p_Uj_\tau=(p_U\otimes \id)\widetilde{F}j_\tau,$$
for each $U\in\cal U$. Consequently, \eqref{227} holds.
\end{pf}

In summary, the following natural correspondences hold:
\begin{equation*}
\left\{
\begin{gathered}
\mbox{quantum principal}\\
\mbox{$G$-bundles}
\end{gathered}\right\}\leftrightarrow\left\{\mbox{$G$-cocycles}
\right\}\leftrightarrow\left\{\mbox{$G_{\!cl}$-cocycles}\right\}
\leftrightarrow\left\{
\begin{gathered}
\mbox{classical principal}\\
\mbox{$G_{\!cl}$-bundles}
\end{gathered}\right\}
\end{equation*}

In this sense, each  quantum  $G$-bundle  $P$
determines  a  classical $G_{\!cl}$-bundle $P_{\!cl}$, and vice versa.

The  correspondence  $P\leftrightarrow  P_{\!cl}$   has   a   simple
geometrical explanation. Each  quantum  group  $G$  is  inherently
inhomogeneous, because it always possesses a nontrivial  classical
part $G_{\!cl}$ consisting of points of $G$ (because of $\e\in  G_{\!cl}$)
and (as far as $\cal A$ is not commutative) a  nontrivial  quantum
part, imaginable as the ``complement'' to $G_{\!cl}$  in  $G$.  It  is
clear that ``transition functions''  being  diffeomorphisms  at  the
level of spaces,  preserve  this  intrinsic  decomposition.  As  a
result, because of the right covariance, transition functions  are
completely determined by their ``restrictions'' on $G_{\!cl}$.

In  fact  the  correspondence  $P\leftrightarrow  P_{\!cl}$  can  be
formulated independently  of  trivialization  systems  $\tau$.  If
$P=(\cal B,i,F)$ is given then the elements of $P_{\!cl}$ are  in  a
natural bijection with *-characters of $\cal B$. In  other  words,
$P_{\!cl}$ is consisting of {\it classical points} of $P$.

Conversely, if $P_{\!cl}$ is given then  $P$  can  be  recovered  by
applying a variant of  the  classical  construction  of  extending
structure groups.

Let $r\colon g\mapsto r_g$ be the (left) action  of  $G_{\!cl}$
on  the  algebra
$S(P_{\!cl})$, induced by the right action of $G_{\!cl}$  on  $P_{\!cl}$.
Let $\zeta^\sstar\colon g\mapsto\zeta_g^\sstar$ be the
left action of $G_{\!cl}$ on $\cal A$. Explicitly,
\begin{gather}
r_g(\varphi)(x)=\varphi(xg)\\
\zeta^\sstar_g=(g^{-1}\otimes \id)\phi.
\end{gather}

Operators $r_g\otimes\zeta_g^\sstar$  are  automorphisms
of a *-algebra $S(P_{\!cl})\otimes \cal A$. Let $\cal B$
be   the   corresponding   fixed-point subalgebra.
It is easy to see that formulas
\begin{gather}
F(b)=(\id\otimes\phi)(b)\\
i(\varphi)=\varphi\pi^\star_M\otimes 1
\end{gather}
(where $\pi_M\colon P_{\!cl}\rightarrow M$ is the projection)
determine *-homomorphisms  $i\colon S(M)\rightarrow\cal B$
and $F\colon \cal B\rightarrow\cal B\otimes \cal A$ such  that
$P=(\cal B,i,F)$ is a principal $G$-bundle over $M$. The initial
bundle $P_{\!cl}$ is realized as the set of classical points of $P$.
\section{Differential Calculus}
Let $P=(\cal B,i,F)$ be a quantum principal $G$-bundle over $M$.
As the starting point for this section, we shall  formulate  three
basic assumptions about a differential  calculus  over $P$.
We  shall  assume  that  the  calculus  on  $P$  is  based  on   a
graded-differential algebra
$$\Omega_P=\sideset{}{^\oplus}\sum_{k\geq 0}\Omega^k_P$$
possessing the following properties:

\smallskip
\bla{diff1} The algebra $\cal B$ is realized as the  0-th  order
summand of $\Omega_P$. In other words, $\Omega_P^0=\cal B$.

\smallskip
\bla{diff2} As a differential algebra, $\Omega_P$ is generated by
$\cal B$.

\smallskip
The next (and the last) assumption  expresses  an  idea  of  local
triviality of the calculus.  It  relates  the  calculus  over  the
bundle $P$  with
differential structures over the structure quantum group $G$
and the base manifold $M$. The  calculus  over  $M$  will  be  the
classical one, based on a graded-differential algebra  $\Omega(M)$
consisting of differential forms. For each open set $U\subseteq M$
we shall denote  by  $\Omega(U)$ and  $\Omega_{\!c}(U)$  algebras  of
differential forms on $U$ (having compact supports).

Concerning the  calculus  over  $G$,  it  will  be  based  on  the
universal   differential   envelope $\Gamma^\wedge$ of   a   given
first-order
differential calculus $\Gamma$ over  $G$.  Properties  of
such  structures  are  collected   in   Appendix   B.   A   symbol
$\grten$ will be used for the graded tensor product  of
graded (differential) algebras.

\smallskip
\bla{diff3} Let $(U,\pi_U)$ be a local trivialization
for $P$ and $\psi_U\colon S_{\!c}(U)\otimes\cal  A\rightarrow  \cal  B$  the
corresponding ``right inverse''. Then $\pi_U$ and $\psi_U$ are  extendible
to homomorphisms
$\pi_U^\wedge\colon \Omega_P\rightarrow\Omega(U)\grten
\Gamma^\wedge$ and $\psi_U^\wedge\colon \Omega_{\!c}(U)\grten
\Gamma^\wedge\rightarrow\Omega_P$   of   (graded-)    differential
algebras.

\smallskip
Property {\it diff2\/} as    well    as    the     fact     that
$\Omega_{\!c}(U)\grten\Gamma^\wedge$     is      generated,
as a differential algebra, by $S_{\!c}(U)\otimes \cal A$,  imply  that
homomorphisms  $\pi_U^\wedge$  and  $\psi_U^\wedge$  are  uniquely
determined. It is easy to see that
\begin{equation}\label{31}
\pi_U^\wedge\psi_U^\wedge(w)=w
\end{equation}
for each $w\in\Omega_{\!c}(U)\grten\Gamma^\wedge$.

For a given open set $V\subseteq M$ let
$I_V^\wedge\subseteq\Omega_P$  be  the   differential   subalgebra
generated by $I_V\subseteq\cal B$.

\begin{lem}\label{lem:31}
\bla{i} Algebras $I_V^\wedge$ are ideals in $\Omega_P$.

\smallskip
\bla{ii} If $(U,\pi_U)$ is a  local  trivialization  for  $P$  and  if
$V\subseteq U$ then
\begin{gather*}
\psi_U^\wedge(\Omega_{\!c}(U)\grten\Gamma^\wedge)=I_V^\wedge\\
\pi_U^\wedge(I_V^\wedge)=\Omega_{\!c}(V)\grten
\Gamma^\wedge.
\end{gather*}
\end{lem}
\begin{pf}   The   second   statement   follows   directly    from
Lemma~\ref{lem:25}  and  definition \eqref{214}.  Concerning
\bla{i}, let us prove it first in a special case  described  in \bla{ii}.
It is sufficient to check that
$b\psi_U^\wedge(f), \, \psi_U^\wedge(f)b,  \,  db\psi_U^\wedge(f)$
and $\psi_U^\wedge(f)db$ belong to $I_V^\wedge=\psi_U^\wedge(
\Omega_{\!c}(V)\grten\Gamma^\wedge), $ for each $f\in
\Omega_{\!c}(V)\grten\Gamma^\wedge$ and $b\in\cal B$. Each
$f\in\Omega_{\!c}(V)\grten\Gamma^\wedge$ can be written  as
a sum of elements of the form $f_0df_1\dots df_k,$ where $f_i\in
S_{\!c}(V)\otimes\cal A$. We have
$
b\psi_U^\wedge(f_0df_1\dots df_k)=b\psi_U(f_0)d\psi_U(f_1)
\dots d\psi_U(f_k)\in I_V^\wedge
$
because $b\psi_U(f_0)\in\psi_U(S_{\!c}(V)\otimes\cal  A).$  Inclusions
$\psi_U^\wedge(f)b\in I^\wedge_V$  follow  in  a  similar  manner.
Further,
$db\psi_U^\wedge(f)=d(b\psi_U^\wedge(f))-b\psi_U^\wedge(df)\in
I^\wedge_V,$ and similarly
$\psi_U^\wedge(f)db\in I^\wedge_V.$

Let $V\subseteq M$ be an arbitrary open set and
$\tau=(\pi_U)_{U\in\cal U}$ an arbitrary trivialization system for
$P$. It is then easy to see that $I_V^\wedge$ is linearly  spanned
by  ideals  $I_{V\cap  U}^\wedge$,  where  $U\in  \cal  U$.  Thus,
$I_V^\wedge$ is an ideal in $\Omega_P$.
\end{pf}
\begin{lem}\label{lem:32}
 Let $\tau$ be a  trivialization  system  for  $P$.
Then  every  map
$\psi_{U\!V}$ from the corresponding $G$-cocycle $\cal{C}_\tau$
is uniquely extendible to a graded-differential automorphism
$\psi_{U\!V}^\wedge\colon \Omega(U\cap V)\grten\Gamma^\wedge
\rightarrow\Omega(U\cap V)\grten\Gamma^\wedge$.
\end{lem}
\begin{pf}   It is sufficient to
construct $\psi_{U\!V}^\wedge$ as automorphisms
of $\Omega_{\!c}(U\cap V)\grten\Gamma^\wedge$.
For   each   $(U,V)\in   N^2(\cal   U)$  let us define
$\psi_{U\!V}^\wedge$ to be the composition of the isomorphisms
$\psi_V^\wedge\colon \Omega_{\!c}(U\cap
V)\grten\Gamma^\wedge\rightarrow I_{U\cap V}^\wedge$ and $
(\psi_U^\wedge)^{-1}\colon I_{U\cap V}^\wedge\rightarrow
\Omega_{\!c}(U\cap V)\grten\Gamma^\wedge$. By  construction
$\psi_{U\!V}^\wedge$ is a grade-preserving differential automorphism
which extends the action of $\psi_{U\!V}$.
Uniqueness follows from  the  fact  that
$S_{\!c}(U\cap V)\otimes\cal A$ generates the differential algebra
$\Omega_{\!c}(U\cap V)\grten\Gamma^\wedge$.
\end{pf}

Consequently,  not  all  differential  structures  over  $G$   are
relevant for our considerations. The calculus $\Gamma$ must be
compatible with transition  functions $\psi_{U\!V}$.  This  is  a
motivation for the following

\begin{defn} A first-order differential calculus $\Gamma$ over
$G$ is called {\it admissible} iff for each $G$-cocycle  $\cal  C$
every transition function $\psi_{U\!V}\colon S(U\cap V)\otimes \cal A
\rightarrow S(U\cap V)\otimes \cal A$ is
extendible  to  a  homomorphism   $\psi_{U\!V}^\wedge\colon \Omega(U\cap
V)\grten\Gamma^\wedge\rightarrow
\Omega(U\cap   V)\grten\Gamma^\wedge$ of differential
algebras. Homomorphisms $\psi^\wedge_{U\!V}$ are grade preserving,
bijective, $\Omega(U\cap V)$-linear and uniquely determined.
\end{defn}
As we shall prove, each admissible calculus  over
$G$, together with  requirements  {\it diff1--3\/},
completely determines the corresponding calculus  $\Omega_P$  over
$P$. As first, the  notion  of
admissibility will be analyzed in more details.

As explained in Appendix A, the Lie algebra  $\lie(G_{\!cl})$  can  be
naturally understood  as  the  space  of  (hermitian)  functionals
$X\colon \cal A\rightarrow \Bbb{C}$
satisfying  $$X(ab)=\e(a)X(b)+\e(b)X(a)$$  for
each $a,b\in \cal A$. Hence, for each $X\in\lie(G_{\!cl})$ the map
\begin{equation}\label{34}
\ellG_X=-(X\otimes \id)\phi
\end{equation}
is  a derivation   on   $\cal   A$.   Further,
$\ellG \colon \lie(G_{\!cl})\rightarrow \Der(\cal A)$ is
a monomorphism of Lie
algebras.   The   image   of   $\ellG$   consists   precisely    of
right-invariant derivations on $\cal A$.

Let  $\cal  C=\Bigl\{\psi_{U\!V}\mid(U,V)\in
N^2(\cal{U})\Bigr\}$   be   a
$G$-cocycle over $(M,\cal U)$. For each $(U,V)\in N^2(\cal U)$  we
shall denote by $\partial^{U\!V}\colon \cal A\rightarrow \Omega(U\cap  V)$
a linear map defined by
\begin{equation}\label{35}
\partial^{U\!V}(a)=g_{V\!U}(a^{(1)})d\bigl(g_{U\!V}(a^{(2)})\bigr).
\end{equation}
It is easy to see that
\begin{equation}\label{36}
\partial^{U\!V}(ab)=\e(a)\partial^{U\!V}(b)+\e(b)\partial^{U\!V}(a)
\end{equation}
for  each  $a,b\in  \cal  A$.  Hence,   $\partial^{U\!V}$   can   be
understood in a natural manner as an element of the space
$\Omega(U\cap V)\otimes \lie(G_{\!cl}).$
\begin{lem}\label{lem:33}
A  first-order  calculus  $\Gamma$  over  $G$  is
admissible iff the following implications hold
\begin{align}
\biggl\{\sum_i a_idb_i=0\biggr\}&\Rightarrow
\biggl\{\sum_i\zeta_g^\sstar(a_i)d\zeta_g^\sstar(b_i)=0\biggr\}\label{37}\\
\biggl\{\sum_ia_idb_i=0\biggr\}&\Rightarrow
\biggl\{\sum_ia_i\ellG_X(b_i)=0\biggr\}\label{38}
\end{align}
for each $g\in G_{\!cl}$ and $X\in \lie(G_{\!cl})$.
\end{lem}
\begin{pf} Maps $\psi_{U\!V}^\wedge$ have the form
\begin{equation}\label{39}
\psi_{U\!V}^\wedge(\alpha\otimes\vartheta)=\alpha\varphi_{U\!V}^\wedge
(\vartheta)
\end{equation}
where $\varphi_{U\!V}^\wedge\rightarrow\Omega(U\cap
V)\grten\Gamma^\wedge$   is
are (unique) graded-differential homomorphism  extending  the  maps
\begin{equation}\label{310}
\varphi_{U\!V}(a)=g_{V\!U}(a^{(1)})\otimes a^{(2)}.
\end{equation}
If $\Sum_ia_idb_i=0$ then
\begin{equation*}
\begin{split}
0&=\varphi^\wedge_{U\!V}\biggl(\sum_ia_idb_i\biggr)=
\sum_i\bigl(g_{V\!U}(a_i^{(1)})\otimes a_i^{(2)}\bigr)d\bigl(
g_{V\!U}(b_i^{(1)})\otimes b_i^{(2)}\bigr)\\
&=\sum_i g_{V\!U}(a_i^{(1)})d\bigl(g_{V\!U}(b_i^{(1)})\bigr)\otimes
a_i^{(2)}b_i^{(2)}+
\sum_i g_{V\!U}(a_i^{(1)})g_{V\!U}(b_i^{(1)})\otimes
a_i^{(2)}db_i^{(2)}\\
&=\sum_i g_{V\!U}(a_i^{(1)}b_i^{(1)})\partial^{V\!U}(b_i^{(2)})
\otimes a_i^{(2)}b_i^{(3)}+\sum_ig_{V\!U}(a_i^{(1)}b_i^{(1)})\otimes
a^{(2)}_idb^{(2)}_i,
\end{split}
\end{equation*}
according to Definition 3.1. Comparing bidegrees we find
\begin{gather*}
\sum_i g_{V\!U}(a_i^{(1)}b_i^{(1)})\partial^{V\!U}(b_i^{(2)})
\otimes a_i^{(2)}b_i^{(3)}=0\\
\sum_i g_{V\!U}(a_i^{(1)}b_i^{(1)})\otimes
a^{(2)}_idb^{(2)}_i=0.
\end{gather*}

Because of arbitrariness of the $G$-cocycle, the  above  equations
imply \eqref{37}--\eqref{38}. Conversely, if  \eqref{37}--\eqref{38}
hold then the formula
\begin{equation}
\sharp_{U\!V}(adb)=g_{V\!U}(a^{(1)}b^{(1)})\otimes a^{(2)}db^{(2)}
+g_{V\!U}(a^{(1)}b^{(1)})\partial^{V\!U}(b^{(2)})\otimes
a^{(2)}b^{(3)}
\end{equation}
consistently defines a linear map $\sharp_{U\!V}\colon \Gamma
\rightarrow\Omega(U\cap V)\grten\Gamma^\wedge$.  It  is
easy to check that
$$
\sharp_{U\!V}(adb)=\varphi_{U\!V}(a)d\varphi_{U\!V}(b)
$$
for each $a,b\in\cal A$. According to Proposition~\ref{pro:B2}
there exists
the unique homomorphism
$\varphi_{U\!V}^\wedge\colon \Gamma^\wedge\rightarrow\Omega(U\cap
V)\grten\Gamma^\wedge$ of graded-differential  algebras
which extends both $\varphi_{U\!V}$ and  $\sharp_{U\!V}$.  Let  us
define maps  $\psi_{U\!V}^\wedge$  by  \eqref{39}.  These  maps  are
differential homomorphisms extending the cocycle maps $\psi_{U\!V}$.
\end{pf}

If implication \eqref{37} holds then the formula
\begin{equation}
\zeta_g^\sstar(adb)=\zeta_g^\sstar(a)d\zeta_g^\sstar(b)
\end{equation}
consistently determines a left action of $G_{\!cl}$ by automorphisms
of $\Gamma$.

It is easy to see that if \eqref{37} holds  then
\begin{equation}\label{315}
\biggl\{\sum_ia_idb_i=0\biggr\}\Rightarrow
\biggl\{\sum_i\ellG_X(a_i)db_i+a_id\ellG_X(b_i)=0\biggr\}
\end{equation}
for  each  $X\in\lie(G_{\!cl})$. In other words, the formula
\begin{equation}\label{316}
\ellG_X(adb)=\ellG_X(a)db+ad\ellG_X(b)
\end{equation}
consistently determines a linear map
$\ellG_X\colon \Gamma\rightarrow\Gamma$. Evidently, the
following equalities hold
\begin{equation*}
\ellG_X(da)=d\ellG_X(a)\qquad
\begin{aligned}
\ellG_X(a\xi)&=\ellG_X(a)\xi+a\ellG_X(\xi)\\
\ellG_X(\xi a)&=\ellG_X(\xi)a+\xi \ellG_X(a).
\end{aligned}
\end{equation*}

Let us now suppose that \eqref{38} holds. In this case the formula
\begin{equation}\label{317}
\iota_X(adb)=a \ellG_X(b)
\end{equation}
consistently determines a bimodule homomorphism $\iota_X\colon \Gamma
\rightarrow\cal A$.

It is worth noticing that the mentioned left action of
$G_{\!cl}$ on $\Gamma$
(and $\cal A$)
is, according to Proposition~\ref{pro:B2},
uniquely  extendible  to  the  left  action  of
$G_{\!cl}$  by  automorphisms  of  the  graded-differential  algebra
$\Gamma^\wedge$.  Moreover,  operators  $\ellG_X$  and  $\iota_X$  are
uniquely extendible to a grade-preserving derivation $\ellG_X\colon
\Gamma^\wedge\rightarrow\Gamma^\wedge$ commuting with $d$, and an
antiderivation
$\iota_X\colon \Gamma^\wedge\rightarrow\Gamma^\wedge$  of order $-1$
respectively.
Classical identities
\begin{gather*}
\iota_X\iota_Y+\iota_Y\iota_X=0\qquad[\ellG_X,\iota_Y]=\iota_{[X,Y]}\\
\ellG_X=d\iota_X+\iota_Xd\qquad\ellG_{[X,Y]}=[\ellG_X,\ellG_Y]
\end{gather*}
hold.
\begin{lem}\label{lem:34}
If $G_{\!cl}$ is connected then the admissibility property
is  equivalent  to implications \eqref{38} and \eqref{315}.
\end{lem}
\begin{pf} Let us suppose that $\Sum_ia_idb_i=0$. It is easy
to see that
\begin{equation}\label{318}
e^{t\ellG_X}\biggl(\sum_ia_idb_i\biggr)=
\sum_i\zeta_{g^t}^\sstar(a_i)d\zeta_{g^t}^\sstar(b_i)=0
\end{equation}
for each $t\in \Re$ and $X\in \lie(G_{\!cl})$, where $t\mapsto g^t$  is
the  1-parameter  subgroup   of   $G_{\!cl}$   generated   by   $X$.
Consequently, there exists an open set $N\ni \e$ such that
\begin{equation}\label{319}
\biggl\{\sum_ia_idb_i=0\biggr\}\Rightarrow
\biggl\{\sum_i\zeta_{g^{\!N}}^\sstar(a_i)d\zeta_{g^{\!N}}^\sstar
(b_i)=0\biggr\}
\end{equation}
for each $g^{\!N}\in N$. If $G_{\!cl}$ is connected then each $g\in
G_{\!cl}$ is a product of some elements from $N$. Inductively applying
\eqref{319} we find that \eqref{37} holds in the full generality.
\end{pf}

On the other hand, implications \eqref{38} and \eqref{315} are
equivalent to the possibility of constructing the maps $\iota_X\colon
\Gamma^\wedge\rightarrow\Gamma^\wedge$.

We pass to the construction of a calculus over $P$. Let us  fix  a
trivialization system $\tau=(\pi_U)_{U\in\cal U}$ for $P$, and  an
admissible first-order calculus $\Gamma$ over $G$.

For each $(U,V)\in N^2(\cal  U)$  the  corresponding  cocycle  map
$\psi_{U\!V}$ admits a natural extension
$\psi_{U\!V}^\wedge\colon \Omega(U\cap V)\grten
\Gamma^\wedge\rightarrow
\Omega(U\cap V)\grten\Gamma^\wedge$ given by
\begin{equation}\label{320}
\psi_{U\!V}^\wedge(\alpha\otimes\xi)=\alpha\varphi^\wedge_{U\!V}(\xi).
\end{equation}
The map $\psi_{U\!V}^\wedge$ can  be  characterized  as  the  unique
graded-differential   homomorphism   extending   $\psi_{U\!V}$.   By
definition, the maps $\psi_{U\!V}^\wedge$ are $\Omega(U\cap  V)$-linear.
In particular, subalgebras $\Omega_{\!c}(W)\grten
\Gamma^\wedge$ are
$\psi_{U\!V}^\wedge$-invariant for each open set $W\subseteq U\cap
V$.
\begin{lem}\label{lem:35}
\bla{i} The maps $\psi_{U\!V}^\wedge$ are bijective and
\begin{equation}
(\psi_{U\!V}^\wedge)^{-1}=\psi_{V\!U}^\wedge.
\end{equation}

\bla{ii} We have
\begin{equation}\label{322}
\psi_{U\!V}^\wedge \psi_{V\!W}^\wedge
(\varphi)=\psi_{U\!W}^\wedge(\varphi)
\end{equation}
for each $(U,V,W)\in N^3(\cal U)$ and $\varphi\in
\Omega_{\!c}(U\cap V\cap W)\grten\Gamma^\wedge$.
\end{lem}
\begin{pf}
Everything follows from similar properties of transition functions
$\psi_{U\!V}$,  and  from  the  fact  that  $\psi_{U\!V}^\wedge$   are
differential homomorphisms.
\end{pf}

Let us consider a graded-differential algebra
$${\cal T}^\wedge=\sideset{}{^\oplus}\sum_{U\in\cal U}
\Omega(U)\grten\Gamma^\wedge$$ and let
$\Omega(P,\tau,\Gamma)\subseteq{\cal T}^\wedge$  be  a  subset
consisting of all $w\in{\cal T}^\wedge$ satisfying
\begin{equation}\label{323}
\psi_{U\!V}^\wedge(\sideset{_V}{_{U\cap V}}\mid\otimes \id)
p_V(w)=(\sideset{_U}{_{U\cap V}}\mid\otimes \id)p_U(w)
\end{equation}
for each $(U,V)\in N^2(\cal U)$, where  $p_U$  are  corresponding
coordinate projections.

All  maps  figuring   in   \eqref{323}   are   graded-differential
homomorphisms. This implies that
$\Omega(P,\tau,\Gamma)$ is a
graded-differential subalgebra of ${\cal T}^\wedge$.

The $0$-th part of $\Omega(P,\tau,\Gamma)$
can be, according  to  Proposition~\ref{pro:29},  identified  with
$\cal B$. By the use of the previous analysis,  it  can  be  shown
easily  that  $\Omega(P,\tau,\Gamma)\subseteq{\cal  T}^\wedge$
satisfies requirements {\it diff2\/} and {\it diff3\/} too.

We shall now prove that
$\Omega(P,\tau,\Gamma)$  is,  up   to isomorphism, the  unique
graded-differential  algebra satisfying conditions {\it diff1--3}.

Let $\cal E$ be an arbitrary algebra possessing this property.
\begin{lem}\label{lem:36}
We have
\begin{equation}\label{324}
\psi_{U\!V}^\wedge(\sideset{_V}{_{U\cap V}}\mid\otimes \id)
\pi_V^\wedge(w)=
(\sideset{_U}{_{U\cap V}}\mid\otimes \id)\pi_U^\wedge(w)
\end{equation}
for each $(U,V)\in N^2(\cal U)$ and $w\in \cal{E}$.
\end{lem}
\begin{pf} Both sides  of  \eqref{324}  are  differential  algebra
homomorphisms coinciding on $\cal B=\cal{E}^0$, according to
Lemma~\ref{lem:28}. Property {\it diff2\/} implies that
they  coincide  on  the whole $\cal E$.
\end{pf}
\begin{lem}\label{lem:37}
The  system  of  maps  $\tau^\wedge=(\pi_U^\wedge)_{U\in\cal   U}$
distinguishes elements of $\cal E$.
\end{lem}
\begin{pf}
Let $(\varphi_U)_{U\in \cal U}$ be an arbitrary smooth
partition of unity for $\cal U$, and let us assume that $w\in \ker(
\pi_U^\wedge)$ for each $U\in \cal U$. Then
$i(\varphi_U)w\in I_U^\wedge\cap \ker(\pi_U^\wedge)$ for each
$U\in  \cal  U$. Hence, we   have
$i(\varphi_U)w=0$. Summing over $\cal U$ we obtain $w=0$.
\end{pf}
\begin{pro}\label{pro:38}
\bla{i} There  exists  the  unique homomorphism $j_\tau^\wedge\colon
\cal E\rightarrow\Omega(P,\tau,\Gamma)$
of   graded-differential
algebras extending the map $j_\tau\colon \cal{B}\rightarrow
\widetilde{\cal{B}}$.

\smallskip
\bla{ii} The map $j_\tau^\wedge$ is bijective.
\end{pro}
\begin{pf}  Let  us  define  a  graded-differential   homomorphism
$j_\tau^\wedge\colon \cal E\rightarrow{\cal T}^\wedge $ by equalities
$$ p_Uj_\tau^\wedge=\pi_U^\wedge.$$
According to Lemma~\ref{lem:36} we have $j_\tau^\wedge(\cal E)
\subseteq\Omega(P,\tau,\Gamma)$. The map
$j_\tau^\wedge\colon \cal  E\rightarrow  \Omega(P,\tau,\Gamma)$   is
injective, according to Lemma~\ref{lem:37}.  The above equality
implies
\begin{equation}\label{327}
j_\tau^\wedge\psi_U^\wedge=\widetilde{\psi}_U^\wedge
\end{equation}
where $\widetilde{\psi}_U^\wedge\colon \Omega_{\!c}(U)\grten
\Gamma^\wedge\rightarrow\Omega(P,\tau,\Gamma)$     is      the
unique graded-differential extension of $\widetilde{\psi}_U\colon
S_{\!c}(U)\otimes\cal A\rightarrow\widetilde{\cal B}$.  Surjectivity  of
$j_\tau^\wedge$ now follows from the fact that $
\Omega(P,\tau,\Gamma)$ is linearly generated by spaces $\im(
\widetilde{\psi}_U)$. Uniqueness of $j_\tau^\wedge$ directly
follows from property {\it diff2}.
\end{pf}

We see that  $\Omega(P,\tau,\Gamma)$ is essentially independent
of a trivialization  system  $\tau$.  For  this  reason  we  shall
simplify the notation and write
$\Omega(P,\Gamma)=\Omega(P,\tau,\Gamma).$  It   is   worth
noticing that the algebra
$\Omega(P,\Gamma)$  can  be  understood   as   the   universal
differential envelope of its first-order  part  (understood  as  a
first-order calculus over $\cal B$).

In the rest of this section algebraic properties of
$\Omega(P,\Gamma)$ will be analyzed in more details.  It  will
be assumed that a trivialization system $\tau$ is fixed.

Let us observe that the formula
\begin{equation}\label{328}
\pi_U^\wedge i^\wedge(\alpha)=\alpha\restr_U
\end{equation}
determines (the unique) graded-differential homomorphism
$i^\wedge\colon \Omega(M)\rightarrow\Omega(P,\Gamma)$ which  extends
the map $i$. The map $i^\wedge$ is injective and
\begin{equation}
i^\wedge(\alpha)w=(-1)^{\partial w\partial\alpha}wi^\wedge(
\alpha)
\end{equation}
for each $\alpha\in\Omega(M)$ and $w\in\Omega(P,\Gamma)$.

As we shall now  see,  it  is  possible  to  introduce  a  natural
coaction  of  $G$  on   $\Omega(P,\Gamma)$, trivialized   in
classical geometry. Let $c\colon\Gamma^\wedge\otimes\cal{A}\rightarrow
\Gamma^\wedge$ be a natural coaction map, defined in Appendix B.
\begin{lem}\label{lem:39} The diagram
\begin{equation}
\begin{CD}
\Bigl\{\Omega(U\cap V)\grten\Gamma^\wedge\Bigr\}\otimes
\cal A @>{\mbox{$\psi_{U\!V}^\wedge\otimes \id$}}>>
\Bigl\{\Omega(U\cap V)\grten\Gamma^\wedge\Bigr\}\otimes
\cal A\\
@V{\mbox{$\id\otimes c$}}VV  @VV{\mbox{$\id\otimes c$}}V\\
\Omega(U\cap V)\grten\Gamma^\wedge @>>
{\mbox{$\psi_{U\!V}^\wedge$}}> \Omega(U\cap V)\grten
\Gamma^\wedge
\end{CD}
\end{equation}
is commutative, for each $(U,V)\in N^2(\cal U)$.
\end{lem}
\begin{pf} A direct computation gives
\begin{equation*}
\begin{split}
\psi_{U\!V}^\wedge(\id\otimes c)(w\otimes a)&=
\psi_{U\!V}^\wedge\bigl(1_{U\cap V}\otimes \k(a^{(1)})\bigr)w(
1_{U\cap V}\otimes a^{(2)})\\
&=\bigl(g_{V\!U}\k(a^{(2)})\otimes \k(a^{(1)})\bigr)
\psi_{U\!V}^\wedge
(w)\bigl(g_{V\!U}(a^{(3)})\otimes a^{(4)}\bigr)\\
&=\bigl(g_{V\!U}(\e(a^{(2)})1)\otimes \k(a^{(1)})\bigr)
\psi_{U\!V}^\wedge
(w)(1_{U\cap V}\otimes a^{(3)})\\
&=\bigl(1_{U\cap V}\otimes \k(a^{(1)})\bigr)\psi_{U\!V}^\wedge
(w)(1_{U\cap V}\otimes a^{(2)})\\
&=(\id\otimes c)(\psi_{U\!V}^\wedge\otimes \id)(w\otimes a). \qed
\end{split}
\end{equation*}
\renewcommand{\qed}{}
\end{pf}
\begin{pro}\label{pro:310} \bla{i}
There exists the unique map
$\Delta\colon \Omega(P,\Gamma)\otimes \cal A\rightarrow
\Omega(P,\Gamma)$ such that the diagram
\begin{equation}\label{330}
\begin{CD}
\Omega(P,\Gamma)\otimes  \cal  A  @>{\mbox{$\pi_U^\wedge\otimes
\id$}}>> \Bigl\{\Omega(U)\grten\Gamma^\wedge\Bigr\}\otimes
\cal A\\
@V{\mbox{$\Delta$}}VV   @VV{\mbox{$\id\otimes c$}}V \\
\Omega(P,\Gamma)                   @>>{\mbox{$\pi_U^\wedge$}}>
\Omega(U)\grten\Gamma^\wedge
\end{CD}
\end{equation}
is commutative, for each $U\in\cal{U}$.

\smallskip
\bla{ii} The following identities hold
\begin{gather}
\Delta(w\otimes ab)=\Delta\bigl(\Delta(w\otimes        a)\otimes
b\bigr)\label{334}\\
\Delta(wu\otimes    a)=\Delta(w\otimes     a^{(1)})\Delta(u\otimes
a^{(2)})\label{333}\\
\Delta(w\otimes 1)=w\label{331}\\
\Delta\bigl(i^\wedge(\alpha)w\otimes a\bigr)=
i^\wedge(\alpha)\Delta(w\otimes
a).\label{332}
\end{gather}
\end{pro}
\begin{pf}
Uniqueness of $\Delta$ is a direct consequence of  the  fact  that
maps $\pi_U^\wedge$ distinguish elements of $\Omega(P,\Gamma)$.
To prove the existence, let us consider a map $\widetilde{
\Delta}\colon{\cal T}^\wedge\otimes \cal{A}
\rightarrow  {\cal  T}^\wedge$
defined by
$$ p_U\widetilde{\Delta}(w \otimes a)=(\id\otimes c)\bigl(p_U(w)\otimes
a\bigr). $$
Lemma~\ref{lem:39} implies that $\widetilde{\Delta}
\bigl(\Omega(P,\Gamma)\otimes\cal A\bigr)\subseteq\Omega(P,\Gamma)$.
The restriction of $\widetilde{\Delta}$ on  $\Omega(P,\Gamma)$
gives the desired map $\Delta\colon \Omega(P,\Gamma)\otimes\cal A
\rightarrow\Omega(P,\Gamma)$.  Evidently, diagram \eqref{330} is
commutative.

A direct computation gives
\begin{equation*}
\begin{split}
\pi_U^\wedge(\Delta(wu\otimes a))&=\sum_{ij}(-1)^
{\partial\vartheta_i\partial\beta_j}\alpha_i\beta_j
\otimes c(\vartheta_i\eta_j\otimes a)\\
&=\sum_{ij}(-1)^{\partial\vartheta_i\partial\beta_j}
\alpha_i\beta_j\otimes c(\vartheta_i\otimes a^{(1)})
c(\eta_j\otimes a^{(2)})\\
&=\pi_U^\wedge(\Delta(w\otimes a^{(1)})\Delta(u\otimes a^{(2)}),
\end{split}
\end{equation*}
Similarly
\begin{multline*}
\pi_U^\wedge\bigl(\Delta(w\otimes ab)\bigr)=\sum_i\alpha_i\otimes c(
\vartheta_i\otimes ab)=\sum_i\alpha_i\otimes c(c(\alpha_i\otimes
a)\otimes b)\\
=\pi_U^\wedge\bigl(\Delta(\Delta(w\otimes a)\otimes b)\bigr),
\end{multline*}
and finally
$$\pi_U^\wedge\Delta(i^\wedge(\alpha) w\otimes a)=
\bigl((\alpha\restr_U)\otimes 1\bigr)\sum_i\alpha_i\otimes c(\vartheta_i
\otimes a)=\pi_U^\wedge\bigl(i^\wedge(\alpha)\Delta(w\otimes a)
\bigr)$$
where $\pi_U^\wedge(w)=\Sum_i\alpha_i\otimes\vartheta_i$ and
$\pi_U^\wedge(u)=\Sum_j\beta_j\otimes\eta_j$. Hence
\eqref{334}--\eqref{332} hold.
\end{pf}

In the case when $\Gamma$ admits the *-structure, or if it  is
right-covariant \cite{W3} the algebra $\Omega(P,\Gamma)$ possesses
a similar property, too. To prove this we need a technical lemma.
\begin{lem}\label{lem:311} \bla{i} If  $\Gamma$  is  a  *-calculus
then  $\psi_{U\!V}^\wedge$  preserve  the  natural  *-structure   on
$\Omega(U\cap V)\grten\Gamma^\wedge$.

\smallskip
\bla{ii} If $\Gamma$ is right-covariant then the diagrams
\begin{equation}\label{335}
\begin{CD}
\Omega(U\cap V)\grten\Gamma^\wedge @>{\mbox{
$\id\otimes \rig_\Gamma^\wedge$}}>> \Bigl\{
\Omega(U\cap V)\grten\Gamma^\wedge\Bigr\}\otimes
\cal A\\
@V{\mbox{$\psi_{U\!V}^\wedge$}}VV @VV{\mbox{$\psi_{U\!V}^\wedge\otimes
\id$}}V\\
\Omega(U\cap V)\grten\Gamma^\wedge @>>{\mbox{
$\id\otimes \rig_\Gamma^\wedge$}}> \Bigl\{
\Omega(U\cap V)\grten\Gamma^\wedge\Bigr\}\otimes
\cal A
\end{CD}
\end{equation}
are commutative. Here,
$\rig_\Gamma^\wedge\colon \Gamma^\wedge\rightarrow\Gamma^\wedge\otimes
\cal A$ is a natural extension of the right action
$\rig_\Gamma\colon \Gamma\rightarrow\Gamma\otimes\cal A$.
\end{lem}
\begin{pf}
Elements of the form $w=\alpha\otimes a_0da_1\dots  da_n$,
where $\alpha\in\Omega(U\cap V)$ and $a_0,a_1,\dots ,a_n\in\cal A$,
linearly span $\Omega(U\cap V)\grten\Gamma^\wedge$.
If $\Gamma$ is *-covariant then
\begin{equation*}
\begin{split}
\psi_{U\!V}^\wedge(w^*)&=(-1)^{n(n-1)/2}
\psi_{U\!V}^\wedge\bigl(\alpha^*\otimes
d(a_n^*)\dots d(a_1^*)a_0^*\bigr)\\
&=(-1)^{n(n-1)/2}\alpha^*d[\varphi_{U\!V}(a_n^*)]\dots
d[\varphi_{U\!V}(a_1^*)]\varphi_{U\!V}(a_0^*)\\
&=(-1)^{n(n-1)/2}\alpha^*d[\varphi_{U\!V}(a_n)]^*\dots
d[\varphi_{U\!V}(a_1)]^*\varphi_{U\!V}(a_0)^*\\
&=\psi^\wedge_{U\!V}(w)^*,
\end{split}
\end{equation*}
according to \eqref{310} and Proposition~\ref{pro:B3}.  Similarly,
if $\Gamma$ is right-covariant  then  Proposition~\ref{pro:B7}
\bla{ii} implies
\begin{equation*}
\begin{split}
(\id\otimes \rig_\Gamma^\wedge)\psi_{U\!V}^\wedge(w)&=
(\id\otimes \rig_\Gamma^\wedge)\bigl(\alpha\varphi_{U\!V}
(a_0)d[\varphi_{U\!V}(a_1)]\dots d[\varphi_{U\!V}(a_n)]\bigr)\\
&=\alpha[(\varphi_{U\!V}\otimes
\id)\phi(a_0)][(d\varphi_{U\!V}\otimes
\id)\phi(a_1)]\dots[(d\varphi_{U\!V}\otimes \id)\phi(a_n)]\\
&=(\psi_{U\!V}^\wedge\otimes \id)(\id\otimes
\rig_\Gamma^\wedge)(w).\qed
\end{split}
\end{equation*}
\renewcommand{\qed}{}
\end{pf}
\begin{pro}\label{pro:312} If  $\Gamma$  is  *-covariant  then
there exists the unique antilinear map $*\colon \Omega(P,\Gamma)
\rightarrow\Omega(P,\Gamma)$ extending $*\colon \cal B\rightarrow
\cal B$, satisfying
$(wu)^*=(-1)^{\partial w\partial u}u^*  w^*$  and  commuting with
$d\colon\Omega(P,\Gamma)\rightarrow
\Omega(P,\Gamma)$. The following identities hold
\begin{gather}
i^\wedge(\alpha^*)=i^\wedge(\alpha)^*\label{339}\\
(w^*)^*=w\label{336}\\
\Delta(w\otimes a)^*=\Delta(w^*\otimes \k(a)^*).\label{338}
\end{gather}
\end{pro}
\begin{pf}
If  $\Gamma$  is  a  *-calculus  then  tensoring  the  natural
*-structure on $\Omega(U)$ with the corresponding  *-structure  on
$\Gamma^\wedge$ and taking the direct sum we obtain a  *-structure
on ${\cal T}^\wedge$. It is easy to see that
$$
(wu)^*=(-1)^{\partial w\partial u}u^*w^*\qquad
d(w^*)=d(w)^*\qquad
i^\wedge(\alpha^*)=i^\wedge(\alpha)^*
$$
for  each  $u,w\in {\cal T}^\wedge$   and   $\alpha\in\Omega(M)$.
According  to  Lemma~\ref{lem:311} \bla{i}, the algebra
$\Omega(P,\Gamma)\subseteq\cal{T}^\wedge$  is
*-invariant. The restriction of  the
*-operation   on   $\Omega(P,\Gamma)$   gives   the    desired
involution.

Applying the definition of $\Delta$ and elementary properties of $c$ we
obtain
\begin{multline*}
\pi_U^\wedge\bigl[\Delta(w\otimes a)^*\bigr]=
\sum_i\alpha_i^*\otimes c(\vartheta_i\otimes a)^*=
\sum_i\alpha_i^*\otimes c(\vartheta_i^*\otimes \k(a)^*)\\=
\pi_U^\wedge\bigl[\Delta(w^*\otimes \k(a)^*)\bigr]
\end{multline*}
for each $U\in\cal U.$  Uniqueness  of $*$  directly  follows  from
property {\it diff2\/} for $\Omega(P,\Gamma).$
\end{pf}
\begin{pro}\label{pro:313}
\bla{i} If $\Gamma$  is  right-covariant  then  there  exists  the
unique homomorphism $F^\wedge\colon \Omega(P,\Gamma)\rightarrow
\Omega(P,\Gamma)\otimes\cal A$ which extends $F$ and such that
\begin{equation}\label{340}
F^\wedge d=(d\otimes\id)F^\wedge.
\end{equation}
The following identities hold
\begin{gather}
F^\wedge i^\wedge(\alpha)=i^\wedge(\alpha)\otimes 1\label{344}\\
(\id\otimes \e)F^\wedge=\id\label{341}\\
(\id\otimes \phi)F^\wedge=(F^\wedge\otimes \id)F^\wedge\label{342}\\
F^\wedge\Delta(w\otimes a)=\sum_k\Delta(w_k\otimes a^{(2)})
\otimes \k(a^{(1)})c_ka^{(3)}\label{343}
\end{gather}
where $F^\wedge(w)=\Sum_k w_k\otimes c_k.$

\smallskip
\bla{ii} If $\Gamma$ is  also  a  *-calculus  then  $F^\wedge$  is
hermitian, in a natural manner.
\end{pro}
\begin{pf}  If  $\Gamma$  is  right-covariant   then   a   map
$F^\wedge_{\cal{T}}\colon \cal T^\wedge\rightarrow\cal
T^\wedge\otimes\cal A$ defined by
$$ (p_U\otimes \id) F^\wedge_{\cal{T}}=(\id\otimes \rig_\Gamma^\wedge)p_U$$
is a homomorphism  which,  according  to  Proposition~\ref{pro:B7}
\bla{ii}, satisfies the following equations
\begin{gather*}
F^\wedge_{\cal{T}}d=(d\otimes \id)F^\wedge_{\cal{T}}\\
(\id\otimes \e)F^\wedge_{\cal{T}}=\id\\
(\id\otimes \phi)F^\wedge_{\cal{T}}=(F^\wedge_{\cal{T}}\otimes
\id)F^\wedge_{\cal{T}}\\
F^\wedge_{\cal{T}}(\alpha)=\alpha\otimes \id
\end{gather*}
where $p_U(\alpha)\in\Omega(U)\otimes 1$ for each $U\in\cal{U}$.
Now Lemma~\ref{lem:311} \bla{ii} implies that $\Omega(P,\Gamma)=
\Omega(P,\tau,\Gamma)$ is a $F^\wedge_{\cal{T}}$-invariant subalgebra of
$\cal T$, in other words we have the inclusion
$F^\wedge_{\cal{T}}\bigl(\Omega(P,\Gamma)\bigr)
\subseteq
\Omega(P,\Gamma)\otimes\cal A$. The restriction of  $F^\wedge_{\cal{T}}$
on $\Omega(P,\Gamma)$ gives the desired map $F^\wedge$.

According to Lemma~\ref{lem:B8}
\begin{equation*}
\begin{split}
(\pi_U^\wedge\otimes \id)F^\wedge\Delta(w\otimes a)&=(\id\otimes
\rig_\Gamma^\wedge)\biggl[\sum_i\alpha_i\otimes c(
\vartheta_i\otimes a)\biggr]\\
&=\sum_{ij}\alpha_i\otimes c(\vartheta_{ij}\otimes a^{(2)})
\otimes \k(a^{(1)})c_{ij} a^{(3)}\\
&=(\pi_U^\wedge\otimes \id)\biggl[\sum_k\Delta(w_k\otimes a^{(2)})
\otimes \k(a^{(1)})c_ka^{(3)}\biggr]
\end{split}
\end{equation*}
for each $U\in\cal U$. Here, $\rig_\Gamma^\wedge(\vartheta_i)
=\Sum_j\vartheta_{ij}\otimes c_{ij}$. Uniqueness of $F^\wedge$  is
a direct consequence of property {\it diff2}. If  $\Gamma$  is
in  addition  *-covariant   then   $\Omega(P,\Gamma)$   is   a
*-subalgebra of  ${\cal  T}^\wedge$  and  $F^\wedge_{\cal{T}}$
is  hermitian, according to Proposition~\ref{pro:B7}.
\end{pf}

{}From  this  moment  we   shall   assume   that   $\Gamma$   is
left-covariant. The space of left-invariant elements  of  $\Gamma$
will be  denoted  by  $\Gamma_{\inv}$.  Further,  $\cal  R\subseteq
\ker(\e)$ will be the right $\cal A$-ideal  which  canonically \cite{W3}
determines this calculus.
\begin{pro}\label{pro:314} A left-covariant calculus  $\Gamma$
is admissible iff
\begin{equation}\label{345}
(X\otimes \id)\ad(\cal R)=\{0\}
\end{equation}
for each $X\in \lie(G_{\!cl})$.
\end{pro}
\begin{pf}
If  $\Gamma$  is  admissible  (and
left-covariant) then the following equality holds
\begin{equation}\label{346}
\varphi^\wedge_{U\!V}\pi(a)=\partial^{V\!U}(a^{(2)})\otimes \k(a^{(1)})
a^{(3)}+1_{U\cap V}\otimes \pi(a).
\end{equation}
Indeed,
\begin{equation*}
\begin{split}
\varphi^\wedge_{U\!V}\pi(a)&=\varphi^\wedge_{U\!V}
\bigl(\k(a^{(1)})da^{(2)}\bigr)\\
&=g_{V\!U}\k(a^{(2)})dg_{V\!U}(a^{(3)})\otimes \k(a^{(1)})a^{(4)}\\
&\phantom{=}+g_{V\!U}\bigl(\k(a^{(2)})a^{(3)}\bigr)
\otimes \k(a^{(1)})da^{(4)}\\
&=g_{U\!V}(a^{(2)})dg_{V\!U}(a^{(3)})\otimes \k(a^{(1)})a^{(4)}\\
&\phantom{=}+g_{V\!U}\bigl(\e(a^{(2)})1\bigr)
\otimes \k(a^{(1)})da^{(3)}\\
&=\partial^{V\!U}(a^{(2)})\otimes \k(a^{(1)})a^{(3)}+1_{U\cap V}\otimes
\pi(a)
\end{split}
\end{equation*}
according to \eqref{35} and \eqref{B37}.

If $a\in\cal R$ then
\begin{equation}\label{347}
\partial^{V\!U}(a^{(2)})\otimes \k(a^{(1)})a^{(3)}=0.
\end{equation}
It is easy to see that, because of arbitrariness of $\tau$,
equations \eqref{347} are equivalent to equations \eqref{345}.

Conversely, let us assume that \eqref{345} holds for each $X\in
\lie(G_{\!cl})$.  To  prove  admissibility  of  $\Gamma$  it   is
sufficient to check implication \eqref{38}, because \eqref{37}  is
satisfied   automatically    for    left-covariant    differential
structures. As a consequence of \eqref{345}, the formula
\begin{equation}\label{348}
\rho_X\bigl(\pi(a)\bigr)=X(a^{(2)})\k(a^{(1)})a^{(3)}
\end{equation}
consistently defines a linear map $\rho_X\colon \Gamma_{\inv}
\rightarrow \cal A$, for each $X\in\lie(G_{\!cl}).$

Now if $\Sum_ia_idb_i=0$ then
\begin{multline*}
0=\sum_i a_ib_i^{(1)}\rho_X\bigl(\pi(b_i^{(2)})\bigr)=
\sum_ia_ib_i^{(1)}X(b_i^{(3)})\k(b_i^{(2)})b_i^{(4)}\\
=\sum_ia_iX(b_i^{(2)})\e(b_i^{(1)})b_i^{(3)}
=\sum_ia_i\ellG_X(b_i)
\end{multline*}
because of \eqref{B43} and the fact that $\Gamma$  is  free over
$\Gamma_{\inv}$ as a left/right $\cal A$-module.
\end{pf}

There exists ``the simplest'' left-covariant admissible  calculus.
It is  based  on  the  right  $\cal  A$-ideal  $\widehat{\cal  R}$
consisting of all  elements  $a\in\ker(\e)$ anihilated by  operators
$(X\otimes  \id)\ad$.  This  calculus  is   also   bicovariant   and
*-covariant. It is analyzed in more details in Appendix C.

Now we are going to construct the  total  ``pull  back''  for  the
right action of $G$ on $P$. We shall assume that  $\Gamma$  is
bicovariant. In this case, as shown in Proposition~\ref{pro:B12},
the comultiplication map admits a natural extension
$\widehat{\phi}\colon \Gamma^\wedge\rightarrow\Gamma^\wedge\grten
\Gamma^\wedge$,  which   is   a   graded   differential   algebra
homomorphism.
\begin{lem}\label{lem:315} The diagram
\begin{equation}\label{349}
\begin{CD}
\Omega(U\cap V)\grten\Gamma^\wedge
@>{\mbox{$\psi^\wedge_{U\!V}$}}>> \Omega(U\cap V)\grten
\Gamma^\wedge\\
@V{\mbox{$\id\otimes\widehat{\phi}$}}VV
@VV{\mbox{$\id\otimes\widehat{\phi}$}}V\\
\Omega(U\cap V)\grten\Gamma^\wedge
\grten\Gamma^\wedge
@>>{\mbox{$\psi^\wedge_{U\!V}\otimes \id$}}>
\Omega(U\cap V)\grten\Gamma^\wedge\grten\Gamma^\wedge
\end{CD}
\end{equation}
is commutative.
\end{lem}
\begin{pf}
All  maps  figuring  in  this   diagram   are   homomorphisms   of
graded-differential  algebras,  and $\Omega(U\cap  V)$-linear
in  a  natural  manner.  Hence,  it  is  sufficent  to  check  the
commutativity in the $0$-th order level. However, this is  just  the
covariance condition for the cocycle maps.
\end{pf}
\begin{pro}\label{pro:316}\bla{i} There exists the  unique homomorphism
$$ \widehat{F}\colon \Omega(P,\Gamma)\rightarrow\Omega(P,\Gamma)
\grten\Gamma^\wedge$$
of graded-differential algebras which extends the map $F$.

\smallskip
\bla{ii} The diagram
\begin{equation}\label{350}
\begin{CD}
\Omega(P,\Gamma) @>{\mbox{$\widehat{F}$}}>>\Omega(P,\Gamma)
\grten\Gamma^\wedge\\
@V{\mbox{$\widehat{F}$}}VV @VV{\mbox{$\id\otimes\widehat{\phi}$}}V\\
\Omega(P,\Gamma)\grten\Gamma^\wedge
@>>{\mbox{$\widehat{F}\otimes \id$}}>
\Omega(P,\Gamma)\grten\Gamma^\wedge
\grten\Gamma^\wedge
\end{CD}
\end{equation}
is commutative and the following identities hold
\begin{gather}
F^\wedge=(\id\otimes p_0)\widehat{F}\label{352}\\
(\id\otimes \e^\wedge)\widehat{F}=\id\label{351}\\
\widehat{F}i^\wedge(\alpha)=i^\wedge(\alpha)\otimes 1.\label{353}
\end{gather}

\bla{iii}  If   $\Gamma$   is   in   addition   *-covariant   then
$\widehat{F}$ preserves canonical *-structures.
\end{pro}
\begin{pf}
Let us consider a linear map
$\widehat{F}_{\cal{T}}\colon \cal{T}^\wedge\rightarrow\cal{T}^\wedge\grten
\Gamma^\wedge$ given by
$$ (p_U\otimes \id)\widehat{F}_{\cal{T}}=(\id\otimes\widehat{\phi})p_U. $$
This map is a homomorphism of graded-differential algebras and
$\widehat{F}_{\cal{T}}(\alpha)=\alpha\otimes 1$ for each
$\alpha\in\cal{T}_M$, where
$$\cal{T}_M=\sideset{}{^\oplus}\sum_{U\in\cal U} \Omega(U).$$
According to Lemma~\ref{lem:315} the algebra
$\Omega(P,\Gamma)=\Omega(P,\tau,\Gamma)$ is
$\widehat{F}_{\cal{T}}$-invariant, in the sense that
$\widehat{F}_{\cal{T}}\bigl(\Omega(P,\Gamma)\bigr)
\subseteq\Omega(P,\Gamma)\grten\Gamma^\wedge$.

Let $\widehat{F}\colon \Omega(P,\Gamma)
\rightarrow\Omega(P,\Gamma)\grten\Gamma^\wedge$  be
the corresponding restriction.
The diagram \eqref{350} and equation \eqref{351}  directly  follow
from \eqref{B48} and \eqref{B50}.

Let us consider a map
$ (\id\otimes p_0)\widehat{F}\colon \Omega(P,\Gamma)\rightarrow
\Omega(P,\Gamma)\otimes\cal A$.
Evidently, this is a homomorphism which extends $F$. Moreover,
$$(\id\otimes p_0)\widehat{F}d(w)=(\id\otimes p_0)
(d\otimes \id+(-1)^{\partial*}\id\otimes
d)\widehat{F}(w)=(d\otimes p_0)\widehat{F}(w)$$
for each $w\in\Omega(P,\Gamma)$. Proposition~\ref{pro:313} implies that
$(\id\otimes p_0)\widehat{F}=F^\wedge$. Uniqueness of $\widehat{F}$
follows from property {\it diff2}.

Finally, if $\Gamma$  is  *-covariant  then  $\widehat{\phi}$  is
hermitian. This implies that $\widehat{F}_{\cal{T}}$ is hermitian,  too.
Hermicity of $\widehat{F}$ also directly follows from hermicity of $F$, and
hermicity of all differentials appearing in the game.
\end{pf}

Let us define the graded *-algebra of {\it horizontal forms} to
be the tensor product
\begin{equation}\label{354}
\hor(P)=\Omega(M)\otimes_M\cal B.
\end{equation}
This   algebra   can   be   understood   as   a   subalgebra    of
$\Omega(P,\Gamma)$ consisting of all $w$ satisfying
\begin{equation}\label{355}
\pi_U^\wedge(w)\in\Omega(U)\otimes\cal A
\end{equation}
for each $U\in\cal U$. By construction, $\hor(P)$ is independent of
a choice of $\Gamma$.

Let us now define a graded algebra of ``verticalized'' differential
forms to be, as a graded vector space
\begin{equation}
\ver(P,\Gamma)=\cal B\otimes\Gamma_{\inv}^\wedge
\end{equation}
while the product is specified by
\begin{equation}\label{357}
(q\otimes\eta)(b\otimes\vartheta)=\sum_k qb_k\otimes
(\eta\circ c_k)\vartheta
\end{equation}
where $\Sum_k b_k\otimes c_k=F(b)$. Here, $\circ$ is the left-invariant
restriction of the coaction map $c$. Associativity of this  product
easily follows from the main properties  of
$F$ and $\circ$.
We  see  that  $\cal  B$   and   $\Gamma_{\inv}^\wedge$   are
subalgebras of $\ver(P,\Gamma)$, in a natural manner. For  each
$U\in\cal U$ the map
$$\pi_U\otimes \id\colon  \ver(P,\Gamma)\rightarrow  S(U)\otimes
\cal A\otimes\Gamma_{\inv}^\wedge\cong S(U)\otimes\Gamma^\wedge $$
becomes a homomorphism of graded algebras. Actually this  property
characterizes the product in $\ver(P,\Gamma)$, because the maps
$\pi_U\otimes \id$ distinguish elements of this algebra.

The algebra $\ver(P,\Gamma)$ can be  equipped  with  a  natural
differential, defined by
\begin{equation}\label{358}
d_v(b\otimes\vartheta)=\sum_k b_k\otimes\pi(c_k)\vartheta+
b\otimes d\vartheta.
\end{equation}
We have
\begin{equation}\label{359}
(\pi_U\otimes \id)d_v(b\otimes\vartheta)=\sum_i\left[\alpha_i
\otimes a_i^{(1)}\otimes\pi(a_i^{(2)})\vartheta+\alpha_i
\otimes a_i\otimes d\vartheta_i\right]
\end{equation}
where $\pi_U(b)=\Sum_i\alpha_i\otimes a_i$. We see that
locally $$d_v\leftrightarrow(\id\otimes d)\colon S(U)\otimes\Gamma^\wedge
\rightarrow S(U)\otimes\Gamma^\wedge.$$

Furthermore, right actions of $G$ on $\cal B$ and
$\Gamma^\wedge_{\inv}$  naturally induce
the right action $F_v$ of $G$ on $\ver(P,\Gamma)$. More precisely,
\begin{equation}\label{360}
F_v(b\otimes\vartheta)=\sum_{kl}b_k\otimes\vartheta_l\otimes
c_kd_l
\end{equation}
where $\adj^\wedge(\vartheta)=\Sum_l\vartheta_l
\otimes d_l.$ This action can be also characterized by relations
\begin{equation}\label{361}
(\pi_U\otimes \id^2)F_v=(\id\otimes\varrho^\wedge_\Gamma)
(\pi_U\otimes \id).
\end{equation}

The differential $d_v$ is $F_v$-covariant, in the sense that
\begin{equation}\label{362}
F_vd_v=(d_v\otimes \id)F_v.
\end{equation}
Indeed, we have
\begin{multline*}
F_vd_v(b\otimes\vartheta)=\sum_{kl}\left(b_k\otimes\pi(c_k^{(3)})
\vartheta_l\otimes c_k^{(1)}\k(c_k^{(2)})c_k^{(4)}d_l+
b_k\otimes d\vartheta_l\otimes c_kd_l\right)\\
=\sum_{kl}\left(b_k\otimes\pi(c_k^{(1)})
\vartheta_l\otimes c_k^{(2)}d_l+
b_k\otimes      d\vartheta_l\otimes       c_kd_l\right)=(d_v\otimes
\id)F_v(b\otimes\vartheta).
\end{multline*}

Graded-differential  algebra  $\ver(P,\Gamma)$  can   be   also
obtained from  $\Omega(P,\Gamma)$  by  factoring through horizontal
forms. More precisely, let $H$ be the ideal in $\Omega(P,\Gamma)$
generated  $di\bigl(S(M)\bigr)$.
Then  $\ver(P,\Gamma)$   is   naturally
isomorphic   to   the   factoralgebra    $\Omega(P,\Gamma)/H.$
Moreover,  $H$  is  a  right-invariant  ideal  and,  according  to
\eqref{359} and \eqref{361}  the  factorized  $F^\wedge$  and  $d$
coincide with $F_v$ and $d_v$ respectively.  We  shall  denote  by
$\pi_v$ the factor projection map.

The homomorphism $\pi_v\colon\Omega(P,\Gamma)\rightarrow\ver(P,\Gamma)$
possesses the following properties
\begin{gather}
(\pi_v\otimes \id)F^\wedge=F_v\pi_v\label{363}\\
\pi_v d=d_v\pi_v\label{364}\\
\pi_v(b)=b\otimes 1.\label{365}
\end{gather}
The last two properties uniquely characterize $\pi_v$.

Finally if $\Gamma$ is *-covariant then $H$ is *-invariant and
there exists the unique *-structure  on  $\ver(P,\Gamma)$  such
that $\pi_v$ is hermitian. Explicitly, this *-structure is given by
\begin{equation}\label{366}
(b\otimes\vartheta)^*=\sum_k b_k^*\otimes (\vartheta^*\circ c_k^*).
\end{equation}
\section{Connections $\&$ Pseudotensorial Forms}
This section is devoted to the study of counterparts of
(pseudo)tensorial forms. In particular, we shall develop the formalism
of connections.

As first, the classical concept of pseudotensoriality will be
translated into the noncommutative context. Let us assume for a
moment that the  bundle is classical. Let us consider a
representation $\widetilde{\rho}\colon G\rightarrow\lin(V)$ in a
vector space $V$. Then a $V$-valued $k$-form $\widetilde{w}$ on $P$ is
called {\it pseudotensorial} of $(\widetilde{\rho},V)$-type \cite{KN} iff
$$ g^*(\widetilde{w})=\rho(g^{-1})\widetilde{w}$$
for each $g\in G$, where $g^*$ is the pull back of the corresponding
right action. The form $\widetilde{w}$ is called {\it tensorial},
if it vanishes whenever at least one argument is vertical.

Pseudotensoriality property can be equivalently formulated in terms of
the map $w\colon V^*\rightarrow\Omega(P)$, where $w(\vartheta)=\vartheta
\widetilde{w}$, via the following diagram
\begin{equation}\label{41}
\begin{CD}
V^* @>{\mbox{$w$}}>>\Omega(P)\\
@V{\mbox{$\rho(g)$}}VV @VV{\mbox{$g^*$}}V\\
V^* @>>{\mbox{$w$}}>\Omega(P)
\end{CD}
\end{equation}
where $\rho$ is the contragradient representation of $\widetilde{\rho}$.
Moreover, $\widetilde{w}$ is tensorial iff
$w(\vartheta)$ is horizontal for each $\vartheta\in V^*$.

Let us turn back to the noncommutative context. Let $P=(\cal{B},i,F)$ be
a quantum principal $G$-bundle over $M$ and $\rho\colon L\rightarrow
L\otimes\cal{A}$ a (nonsingular) representation \cite{W2} of $G$ in a
complex vector space $L$. Let $\Gamma$ be an admissible right-covariant
calculus over $G$. The above diagram naturally suggests to define
pseudotensorial forms as linear maps $w\colon
L\rightarrow\Omega(P,\Gamma)$ such that the diagram
\begin{equation}\label{42}
\begin{CD}
L@>{\mbox{$w$}}>>\Omega(P,\Gamma)\\
@V{\mbox{$\rho$}}VV @VV{\mbox{$F^\wedge$}}V\\
L\otimes\cal{A} @>>{\mbox{$w\otimes\id$}}>
\Omega(P,\Gamma)\otimes\cal{A}
\end{CD}
\end{equation}
is commutative.

Let us denote by $\psi(P,\rho,\Gamma)$ the space of corresponding
pseudotensorial forms. This space is naturally graded
\begin{equation}\label{43}
\psi(P,\rho,\Gamma)=\sideset{}{^\oplus}\sum_{i\geq 0}\psi^i(P,\rho,\Gamma)
\end{equation}
where the grading is induced from $\Omega(P,\Gamma)$. Strictly speaking the
above decomposition holds if $L$ is finite-dimensional. The space
$\psi(P,\rho,\Gamma)$ is a bimodule over $\Omega(M)$, in a natural manner.
According to \eqref{340}, the space of pseudotensorial forms is
invariant under compositions with $d\colon\Omega(P,\Gamma)\rightarrow
\Omega(P,\Gamma)$.

We shall denote by $\tau(P,\rho)$ the subspace consisting of {\it
tensorial forms} $w$, characterized by
\begin{equation}\label{44}
w(L)\subseteq\hor(P).
\end{equation}
Actually $\tau(P,\rho)$ is a graded $\Omega(M)$-submodule of
$\psi(P,\rho,\Gamma)$. Let us observe that
$\tau(P,\rho)$ is independent of a specification of $\Gamma$.

If $L$ is endowed with an antilinear involution $*\colon L\rightarrow L$
such that $\rho$ is hermitian, in a natural manner, and if $\Gamma$ is a
*-calculus then the formula
$$ w^*(\vartheta)=(w(\vartheta^*))^*$$
defines a *-structure on $\psi(P,\rho,\Gamma)$. The space $\tau(P,\rho)$
is *-invariant.

Tensorial forms possess a simple local representation.
\begin{pro}\label{pro:41} \bla{i} For each $w\in\tau(P,\rho)$ and
$U\in\cal{U}$ there exists the unique linear map $\varphi_U\colon
L\rightarrow\Omega(U)$ such that
\begin{equation}\label{45}
\pi_U^\wedge w=(\varphi_U\otimes\id)\rho.
\end{equation}
We have
\begin{equation}\label{46}
\bigl(\varphi_V(\vartheta){\restr}_{U\cap V}\bigr)=
\sum_k\bigl(\varphi_U(\vartheta_k){\restr}_{U\cap V}\bigr)g_{U\!V}(c_k)
\end{equation}
for each $\vartheta\in L$ and $(U,V)\in N^2(\cal{U})$, where
$\Sum_k\vartheta_k\otimes c_k=\rho(\vartheta)$.

\smallskip
\bla{ii} Conversely, if maps $\varphi_U$ satisfy equalities \eqref{46}
then there exists the unique $w\in\tau(P,\rho)$ such that
\eqref{45} holds.
\end{pro}
\begin{pf}
 We have
$$\pi_U^{\wedge}w(L)\subseteq\Omega(U)\otimes\cal{A}$$
for each $w\in\tau(P,\rho)$ and $U\in\cal{U}$.
On the other hand \eqref{42} is equivalent to the following equations
\begin{equation}\label{48}
(\id\otimes\phi)\bigl[\pi_U^\wedge w(\vartheta)\bigr]=
\bigl(\pi_U^\wedge w\otimes\id\bigr)\rho(\vartheta).
\end{equation}
Acting by $\id\otimes\e\otimes\id$ on both sides of this equation we obtain
\eqref{45} with $\varphi_U=(\id\otimes\e)\pi_U^\wedge w$. Conversely,
a direct verification shows that \eqref{48} follows from \eqref{45}.

Let us now analyze how $\varphi_U$ and $\varphi_V$ are related on the
overlaping of regions $U$ and $V$.

For an arbitrary system of linear maps $\varphi_U\colon L\rightarrow
\Omega(U)$, the formula \eqref{45} determines a linear map $w\colon
L\rightarrow\cal{T}^\wedge$. According to \eqref{323} a
necessary and sufficient condition for the inclusion $w(L)\subseteq
\Omega(P,\tau,\Gamma)$ can be written in the form
\begin{equation}\label{49}
(_U\!{\vert}_{U\cap V}\otimes\id)(\varphi_U\otimes\id)\rho(\vartheta)
=\sum_k{}_V\!{\vert}_{U\cap V}\bigl(\varphi_V(\vartheta_k)\bigr)g_{V\!U}
(c_k^{(1)})\otimes c_k^{(2)},
\end{equation}
which is equivalent to \eqref{46}.
\end{pf}

{}From this moment it will be assumed that $\Gamma$ is the simplest
left-covariant admissible calculus. Explicitly, $\Gamma$ is a
first-order calculus based on the right ideal $\widehat{\cal{R}}$
consisting of all $a\in\ker(\e)$ such that $(X\otimes\id)\ad(a)=0$ for
each $X\in\lie(G_{\!cl})$. As explained in Appendix C, this is a
bicovariant *-calculus.

Furthermore, we shall restrict the consideration to the case
$$ L=\Gamma_{\inv}\qquad\rho=\adj.$$
In this case we shall simplify the notation and write
$\Omega(P)$, $\ver(P)$, $\tau(P)$ and $\psi(P)$ for the corresponding
algebras and modules.

Finally, we shall fix a section $\eta\colon \L^*\rightarrow\Gamma_{\inv}$
of $\nu\colon\Gamma_{\inv}\rightarrow \L^*$ (Appendix C) which intertwines
*-structures and adjoint actions of $G_{\!cl}$. Hence we can write
\begin{equation}\label{411}
\Gamma_{\inv}=\L^*\oplus\ker(\nu),
\end{equation}
with $\eta\nu$ playing the role of the projection on the first factor.

If $\varphi_U$ are local representatives of $w\in\tau(P)$ then maps
$$\varphi_{U}^{cl}=\varphi\eta\nu\qquad\varphi_U^{\perp}=\varphi_U(1-\eta\nu)
$$
satisfy \eqref{46}, too. This, together with Proposition~\ref{pro:41},
enables us to introduce the ``classical'' and the ``quantum'' component
of $w$, by
$$\pi_U^\wedge w_{\!cl}=(\varphi_U^{cl}\otimes\id)\adj
\qquad\pi_U^\wedge w_{\perp}=(\varphi_U^{\perp}\otimes\id)\adj.$$
By construction, $$w=w_{\!cl}+w_{\perp}.$$

We shall denote by $\tau_{\!cl}(P)$ and $\tau_{\perp}(P)$ corresponding
mutually complementary graded *-$\Omega(M)$-submodules of $\tau(P)$.
Elements of $\tau_{\!cl}(P)$ will be called {\it classical} tensorial forms.
\begin{pro}\label{pro:42} A tensorial form $w$ is classical
iff the diagram
\begin{equation}\label{412}
\begin{CD}
\Gamma_{\inv}\otimes\cal{A}
@>{\mbox{$w\otimes\id$}}>>\Omega(P)\otimes\cal{A}\\
@V{\mbox{$\circ$}}VV @VV{\mbox{$\Delta$}}V\\
\Gamma_{\inv} @>>{\mbox{$w$}}> \Omega(P)
\end{CD}
\end{equation}
is commutative.
\end{pro}
\begin{pf} Let us suppose that $w$ is classical. In local trivialization
terms, this means
\begin{equation}\label{413}
\varphi_U(\vartheta\circ a)=\e(a)\varphi_U(\vartheta),
\end{equation}
for each $\vartheta\in\Gamma_{\inv}$, $U\in\cal{U}$ and $a\in\cal{A}$. On
the other hand, according to \eqref{330} and \eqref{B36}, commutativity
of \eqref{412} is equivalent to equalities
\begin{multline}\label{414}
\pi_U^\wedge\bigl[w(\vartheta\circ
a)\bigr]=(\varphi_U\otimes\id)\adj(\vartheta\circ a)=\sum_k\varphi_U(
\vartheta_k\circ a^{(2)})\otimes\k(a^{(1)})c_ka^{(3)}\\
=\sum_k\varphi_U(\vartheta_k)\otimes\k(a^{(1)})c_k a^{(2)}=\pi_U^\wedge
\Delta\bigl[w(\vartheta)\otimes a\bigr],
\end{multline}
where $\adj(\vartheta)=\Sum_k\vartheta_k\otimes c_k$.

If \eqref{413} holds then, evidently, \eqref{414} holds. Conversely, if
\eqref{414} holds then acting by $\id\otimes\e$ on both sides of the
third equality we obtain \eqref{413}.
\end{pf}

We pass to the study of connection forms.
\begin{defn}\label{def:41}
A {\it connection} on $P$ is every pseudotensorial 1-form $\omega$
satisfying
\begin{align}
\omega(\vartheta^*)&=\omega(\vartheta)^*\label{415a}\\
\pi_v\omega(\vartheta)&=1\otimes\vartheta\label{415b}
\end{align}
for each $\vartheta\in\Gamma_{\inv}$.
\end{defn}

Condition \eqref{415b} plays the role of the classical requirement that
connections map fundamental vector fields into their generators.
Connections naturally form an infinite-dimensional affine space (as far as
$\Gamma_{\inv}$ is non-trivial).
\begin{lem}\label{lem:43} \bla{i}
Each quantum principal bundle $P$ admits a connection.

\smallskip
\bla{ii} For an arbitrary connection $\omega$ on $P$, and a linear map
$\alpha\colon\Gamma_{\inv}\rightarrow\Omega(P)$, the map $\alpha+\omega$
is a connection iff $\alpha$ is a hermitian 1-order tensorial form.
\end{lem}
\begin{pf}
Let us consider an arbitrary smooth partition of unity
$\bigl(\rho_U\bigr)_{U\in\cal{U}}$ for $\cal{U}$, and define a
map $\omega\colon\Gamma_{\inv}\rightarrow\Omega(P)$ by
\begin{equation}\label{416}
\omega(\vartheta)=\sum_{U\in\cal{U}}\psi^\wedge_U(\rho_U\otimes
\vartheta).
\end{equation}
This map is a connection on $P$. The second statement easily follows
from Definition~\ref{def:41}.
\end{pf}

Let $\con(P)$ be the affine space of all connections on $P$. The
following proposition describes connections in terms of gauge
potentials.
\begin{pro}\label{pro:44} \bla{i} For each $\omega\in\con(P)$ there exist
the unique linear maps $A_U\colon\Gamma_{\inv}\rightarrow\Omega(U)$ such
that
\begin{equation}\label{417}\pi_U^\wedge\omega(\vartheta)=
\sum_kA_U(\vartheta_k)\otimes c_k+1_U\otimes\vartheta
\end{equation}
for each $U\in\cal{U}$, where $\Sum_k\vartheta_k\otimes c_k
=\adj(\vartheta)$. These maps are hermitian and
\begin{equation}\label{418}
\bigl(A_V(\vartheta){\restr}_{U\cap V}\bigr)=
\sum_k\bigl(A_U(\vartheta_k){\restr}_{U\cap V}\bigr)g_{U\!V}(c_k)+
\partial_{U\!V}(\vartheta)
\end{equation}
for each $(U,V)\in N^2(\cal{U})$, where $\partial_{U\!V}\pi=\partial^{U\!V}$.

\smallskip
\bla{ii} Conversely, if hermitian maps $A_U\colon\Gamma_{\inv}\rightarrow
\Omega(U)$ are given such that \eqref{418} holds, then the formula
\eqref{417} determines a connection on $P$.
\end{pro}
\begin{pf}
The proof is essentially the
same as for Proposition~\ref{pro:41}.
\end{pf}
\begin{defn}\label{def:42}
A connection $\omega$ is called {\it classical} iff the diagram
\begin{equation}\label{419}
\begin{CD}
\Gamma_{\inv}\otimes\cal{A}@>{\mbox{$\circ$}}>>\Gamma_{\inv}\\
@V{\mbox{$\omega\otimes\id$}}VV @VV{\mbox{$\omega$}}V\\
\Omega(P)\otimes\cal{A} @>>{\mbox{$\Delta$}}>\Omega(P)
\end{CD}
\end{equation}
is commutative.
\end{defn}
\begin{pro}\label{pro:45}
A connection $\omega$ is classical iff $$A_U\eta\nu=A_U\iff
A_U(\vartheta\circ a)=\e(a)A_U(\vartheta),$$
for each $U\in\cal{U}$.
\end{pro}
\begin{pf}
A similar reasoning as in the proof of Proposition~\ref{pro:42}.
\end{pf}

Every connection can be written as a sum of a classical connection, and
a ``purely quantum'' part.
\begin{pro}\label{pro:46}
For each $\omega\in\con(P)$ there exist the unique classical connection
$\omega_{\!cl}$ and hermitian tensorial 1-form
$\omega_\perp\in\tau_\perp(P)$ such that
\begin{equation}\label{w=w+w}
\omega=\omega_{\!cl}+\omega_\perp.
\end{equation}
\end{pro}
\begin{pf}
Let us start from the corresponding gauge potentials $A_U$ and define
$$A_U^{cl}=A_U\eta\nu\qquad A_U^{\perp}=A_U-A_U^{cl}.$$
{}From \eqref{418} it follows that
\begin{align*}
\bigl(A_V^{cl}(\vartheta){\restr}_{U\cap V}\bigr)&=\sum_k
\bigl(A_U^{cl}(\vartheta_k){\restr}_{U\cap V}\bigr)g_{U\!V}(c_k)
+\partial_{U\!V}(\vartheta)\\
\bigl(A_V^{\perp}(\vartheta){\restr}_{U\cap V}\bigr)&=
\sum_k\bigl(A_U^{\perp}(\vartheta_k){\restr}_{U\cap V}\bigr)g_{U\!V}(c_k).
\end{align*}
It is easy to see that $A_{U}^{cl}$ and $A_U^{\perp}$ are hermitian.
Hence, there exist a classical connection $\omega_{\!cl}$ and a
hermitian element
$\omega_\perp\in\tau_\perp^1(P)$ such that
\begin{align*}
\pi_U^\wedge\omega_{\!cl}(\vartheta)&=(A_U^{cl}\otimes\id)\adj(\vartheta)+
1_U\otimes\vartheta\\
\pi_U^\wedge\omega_\perp(\vartheta)&=(A_U^\perp\otimes\id)\adj(\vartheta)
\end{align*}
for each $\vartheta\in\Gamma_{\inv}$. Evidently, \eqref{w=w+w} holds.
This decomposition is unique, because of mutual complementarity between
$\tau_{\!cl}(P)$ and $\tau_\perp(P)$.
\end{pf}

{}From this moment it will be assumed that the subalgebra
$\Gamma_{\inv}^\wedge$ of left-invariant elements is realized as a
complement to the space $S_{\inv}^\wedge
\subseteq\Gamma_{\inv}^\otimes$, with the
help of a linear section $\iota\colon\Gamma_{\inv}^\wedge\rightarrow
\Gamma_{\inv}^\otimes$ of the factorization map,
which intertwines *-structures and adjoint actions of $G$.
Here $S^\wedge_{\inv}$ is the left-invariant part of the ideal
$S^\wedge\subseteq\Gamma^{\otimes}$ and $\Gamma_{\inv}^\otimes$
is the tensor
algebra over $\Gamma_{\inv}$ (Appendix B).

It is easy to see (for example, applying a quantum analog of the method
of group projectors) that $\iota$ always exists. If $\Gamma_{\inv}$ is
finite-dimensional then $\iota$ can be constructed by identifying
$\Gamma_{\inv}^\wedge$ with the orthocomplement of $S_{\inv}^\wedge$, with
respect to an appropriate scalar product.

However, it is important to mention that in various interesting
situations (for example, if $G=S_\mu U(2)$ and $\mu\in(-1,1)\setminus
\{0\}$) the space $\Gamma_{\inv}$ will be infinite-dimensional.

For each connection $\omega$, let us denote by
$\omega^\otimes\colon\Gamma_{\inv}^\otimes\rightarrow\Omega(P)$ the
corresponding unital multiplicative extension. Let
$\omega^\wedge\colon\Gamma^\wedge_{\inv}\rightarrow\Omega(P)$ be the
composition of maps $\iota$ and $\omega^{\otimes}$.
\begin{pro}\label{pro:47}
\bla{i} The diagram
\begin{equation}\label{424}
\begin{CD}
\Gamma_{\inv}^\wedge @>{\mbox{$\omega^\wedge$}}>>\Omega(P)\\
@V{\mbox{$\adj^\wedge$}}VV @VV{\mbox{$F^\wedge$}}V\\
\Gamma_{\inv}^\wedge\otimes\cal{A} @>>{\mbox{$\omega^\wedge
\otimes\id$}}> \Omega(P)\otimes\cal{A}
\end{CD}
\end{equation}
is commutative.

\smallskip
\bla{ii} We have
\begin{equation}
\pi_v\omega^\wedge(\vartheta)=1\otimes\vartheta
\end{equation}
for each $\vartheta\in\Gamma_{\inv}^\wedge$.

\smallskip
\bla{iii} The map $\omega^\wedge$ is *-preserving.

\smallskip
\bla{iv} If $\omega$ is classical then $\omega^\wedge$ is multiplicative and
the diagram
\begin{equation}\label{425}
\begin{CD}
\Gamma_{\inv}^\wedge\otimes\cal{A} @>{\mbox{$\omega^\wedge\otimes\id$}}>>
\Omega(P)\otimes\cal{A}\\
@V{\mbox{$\circ$}}VV @VV{\mbox{$\Delta$}}V\\
\Gamma^\wedge_{\inv} @>>{\mbox{$\omega^\wedge$}}> \Omega(P)
\end{CD}
\end{equation}
is commutative.
\end{pro}
\begin{pf}
Property \bla{i} is a simple consequence of the pseudotensoriality of
$\omega$ and of the $\adj^\otimes$-invariance of
$\iota\bigl(\Gamma_{\inv}^\wedge\bigr)$. Property \bla{ii} follows from
\eqref{415b}, and the multiplicativity of $\pi_v$.

To prove \bla{iii}, it is sufficient to observe that $\omega^\otimes$
intertwines *-structures on $\Gamma_{\inv}^\otimes$ and $\Omega(P)$.

Let us assume that $\omega$ is classical. We shall prove that
$\omega^\otimes$ vanishes on the ideal
$S_{\inv}^\wedge\subseteq\Gamma_{\inv}^\otimes$. In accordance with
considerations performed in Appendix B, it is sufficient to check that
\begin{equation}\label{426}
\omega^\otimes\bigl[\pi(a^{(1)})\otimes\pi(a^{(2)})\bigr]=0
\end{equation}
for each $a\in\widehat{\cal{R}}$. In the local trivialization system,
this is equivalent to the following equalities
\begin{multline*}
\bigl[(A_U\otimes\id)\adj\pi(a^{(1)})\bigr]\pi(a^{(2)})
+\pi(a^{(1)})\bigl[(A_U\otimes\id)\adj\pi(a^{(2)})\bigr]\\
{}+\bigl[(A_U\otimes\id)\adj\pi(a^{(1)})\bigr]
\bigl[(A_U\otimes\id)\adj\pi(a^{(2)})\bigr]=0.
\end{multline*}
A direct calculation shows that the last term, as well as the
sum of the first two, vanishes. Consequently $\omega^\wedge$ is
multiplicative.

Commutativity of \eqref{425} is a direct consequence of \eqref{333},
\eqref{B24} and \eqref{412}, and
the multiplicativity of $\omega^\wedge$.
\end{pf}

With the help of $\omega^\wedge$ the space $\Omega(P)$ can be
naturally decomposed into a tensor product of $\hor(P)$ and
$\Gamma_{\inv}^\wedge$.

Let us suppose that $\vh(P)=\hor(P)\otimes\Gamma_{\inv}^\wedge$ is
endowed with a graded *-algebra structure, via the natural
indentification
\begin{equation}\label{429}
\hor(P)\otimes\Gamma_{\inv}^\wedge\leftrightarrow\Omega(M)
\otimes_M\ver(P).
\end{equation}
The algebra $\vh(P)$ represents ``vertically-horizontally'' decomposed
forms on the bundle. We shall denote by $F_{v\!h}$ the natural right
action of $G$ on $\vh(P)$.

For each $\omega\in\con(P)$ the formula
\begin{equation}\label{427}
m_\omega(\varphi\otimes\vartheta)=\varphi\omega^\wedge(\vartheta)
\end{equation}
defines a linear grade-preserving map $m_\omega\colon\vh(P)\rightarrow
\Omega(P)$.
\begin{pro}\label{pro:48} \bla{i} The map $m_\omega$ is bijective.

\smallskip
\bla{ii} The diagram
\begin{equation}\label{428}
\begin{CD}
\vh(P)@>{\mbox{$m_\omega$}}>> \Omega(P)\\
@V{\mbox{$F_{v\!h}$}}VV @VV{\mbox{$F^\wedge$}}V\\
\vh(P)\otimes\cal{A} @>>{\mbox{$m_\omega\otimes\id$}}>
\Omega(P)\otimes\cal{A}
\end{CD}
\end{equation}
is commutative.

\smallskip
\bla{iii} If $\omega$ is classical then $m_\omega$ is an isomorphism of
graded *-algebras.
\end{pro}
\begin{pf}
As first we prove that $m_\omega$ is injective. Each
$\alpha\in\vh(P)\setminus\{0\}$ can be written in the form
$\alpha=\Sum_i w_i\otimes\vartheta_i+\psi$, where $\vartheta_i\in
\Gamma_{\inv}^{\wedge k}$
are homogeneous linearly independent elements and
$w_i\neq 0$, while $\psi$ is the element having the second degrees
less then $k$. If $m_\omega(\alpha)=0$ then
$$\sum_i\pi_U(w_i)\vartheta_i=0$$
for each $U\in\cal{U}$. This implies $\Sum_i w_i\otimes\vartheta_i=0$,
which is a contradiction.

In order to prove that $m_\omega$ is surjective, it is sufficient to
check that
$$\psi_U^\wedge\bigl(\Omega_{\!c}(U)\otimes\Gamma^{\wedge k}\bigr)
\subseteq m_\omega\bigl(\vh(P)\bigr)$$
for each $U\in\cal{U}$ and $k\geq 0$.

For $k=0$ the statement is obvious. Let us suppose that the above
inclusion holds for degrees up to some fixed $k$. Equation
\eqref{417} together with the definition of $\omega^\wedge$ gives
\begin{equation}\label{431}
\pi_U^\wedge\bigl[m_\omega(w\otimes\vartheta)\bigr]
=\sum_i\alpha_i\otimes a_i\vartheta+\beta
\end{equation}
where $\vartheta\in(\Gamma_{\inv}^\wedge)^{k+1}$ and
$w=\psi_U^\wedge\Bigl(\Sum_i\alpha_i\otimes a_i\Bigr)$, while
$\beta\in\Omega_{\!c}(U)\otimes\Gamma^\wedge$, with the second degrees less
then $k+1$.

Acting by $\psi_U^\wedge$ on both sides of \eqref{431} we get
$$\psi_U^\wedge\Bigl(\sum_i\alpha_i\otimes a_i\vartheta\Bigr)
=m_\omega(w\otimes\vartheta)-\psi_U^\wedge(\beta).$$
By the inductive assumption,
the right-hand side of the above equality belongs to $\im(m_\omega)$.
Hence $m_\omega$ is bijective.

The commutativity of \eqref{428} is a direct consequence of \eqref{427},
and Proposition~\ref{pro:47} \bla{i}.

Finally, let us suppose that $\omega$ is classical. According to
Proposition~\ref{pro:47} \bla{iv}
and definition \eqref{357} of the product in $\ver(P)$, we have
\begin{equation}\label{prod-vh}
(u\otimes\vartheta)(w\otimes\eta)=(-1)^{\partial\vartheta\partial w}
\sum_kuw_k\otimes (\vartheta\circ c_k)\eta
\end{equation}
and hence
\begin{equation*}
\begin{split}
m_\omega\bigl[(u\otimes\vartheta)(w\otimes\eta)\bigr]&=
(-1)^{\partial w\partial\vartheta}\sum_kuw_k\omega^{\wedge}
(\vartheta\circ c_k)\omega^\wedge(\eta)\\
&=(-1)^{\partial w\partial\vartheta}\sum_kuw_k\Delta\bigl(\omega^{\wedge}
(\vartheta)\otimes c_k\bigr)\omega^\wedge(\eta)\\
&=u\omega^\wedge(\vartheta)w\omega^\wedge(\eta)=
m_\omega(u\otimes\vartheta)m_\omega(w\otimes\eta).
\end{split}
\end{equation*}
Here $F^\wedge(w)=\Sum_kw_k\otimes c_k$ and we have used the identity
\begin{equation}\label{432}
\sum_k w_k\Delta(\alpha\otimes c_k)=(-1)^{\partial \alpha\partial w}
\alpha w,
\end{equation}
where $\alpha$ is arbitrary (and $\omega$ is horizontal).
Similarly, the *-structure on $\vh(P)$ is given by
\begin{equation}\label{*-vh}
(w\otimes\vartheta)^*=\sum_k w_k^*\otimes(\vartheta^*\circ c_k^*)
\end{equation}
and hence
\begin{multline*}
m_\omega\bigl[(w\otimes\vartheta)^*\bigr]=
\sum_kw_k^*\omega^\wedge(\vartheta^*\circ c_k^*)=
\sum_kw_k^*\Delta\bigl(\omega^\wedge(\vartheta)^*\otimes c_k^*\bigr)\\
=(-1)^{\partial w\partial\vartheta}\omega^\wedge(\vartheta)^* w^*
=\bigl[m_\omega(w\otimes\vartheta)\bigr]^*.\qed
\end{multline*}
\renewcommand{\qed}{}
\end{pf}

It is of some interest to analyze in more details the question of the
multiplicativity of $\omega^\wedge$.

\begin{defn}\label{def:43}
A connection $\omega$ is called {\it multiplicative} iff
$$\omega^\otimes(S_{\inv}^\wedge)=\{0\}.$$
\end{defn}

Equivalently, $\omega$  is multiplicative iff $\omega^\wedge$ is a
multiplicative map. In this case $\omega^\wedge$ is independent of the
embedding $\iota$, and coincides with $\omega^\otimes/S_{\inv}^\wedge$. As
already mentioned, the multiplicativity of $\omega^\wedge$ is equivalent
to \eqref{426}. This gives a quadratic constraint in $\con(P)$.
In the general case the left-hand side of \eqref{426} determines a
linear map $r_\omega\colon\widehat{\cal{R}}\rightarrow\Omega(P)$. This
map ``measures'' a lack of multiplicativity of $\omega$.
\begin{pro}\label{pro:49} We have
\begin{equation}\label{433}
r_\omega=m_\Omega(w_\perp\pi\otimes\omega_\perp\pi)\phi,
\end{equation}
where $m_\Omega$ is the product map in $\Omega(P)$. In local terms
\begin{equation}\label{434}
\pi_U^\wedge
r_\omega=(r_\omega^U\otimes\id)\bigl(\ad{\restr}\widehat{\cal{R}}\bigr)
\end{equation}
where $r_\omega^U(a)=A_U^\perp\pi(a^{(1)})A_U^\perp\pi(a^{(2)})$.
In particular $r_\omega$ is a horizontally-valued map.
\end{pro}
\begin{pf} Using local expressions for $\omega_{\!cl}$ and $\omega_\perp$,
equations \eqref{426} and \eqref{B40}, and Proposition~\ref{pro:45} we obtain
\begin{equation*}
\begin{split}
\pi_U^\wedge r_\omega(a)-\pi_U^\wedge m_\Omega(\omega_\perp\pi
\otimes\omega_\perp\pi)\phi(a)&=\pi_U^\wedge m_\Omega(\omega_\perp\pi
\otimes\omega_{\!cl}\pi)\phi(a)\\
&\phantom{=}+\pi_U^\wedge m_\Omega(\omega_{\!cl}\pi
\otimes\omega_{\perp}\pi)\phi(a)\\
&=A_U^\perp\pi(a^{(2)})A_U^{cl}\pi(a^{(3)})\otimes\k(a^{(1)})a^{(4)}\\
&\phantom{=}+A_U^{cl}\pi(a^{(2)})A_U^{\perp}\pi(a^{(3)})
\otimes\k(a^{(1)})a^{(4)}\\
&\phantom{=}+A_U^\perp\pi(a^{(2)})\otimes\k(a^{(1)})a^{(3)}\pi(a^{(4)})\\
&\phantom{=}-A_U^\perp\pi(a^{(3)})\otimes\pi(a^{(1)})\k(a^{(2)})a^{(4)}\\
&=A_U^\perp\pi(a^{(3)})A_U^{cl}\pi\bigl(\k(a^{(2)})a^{(4)}\bigr)\otimes
\k(a^{(1)})a^{(5)}\\
&\phantom{=}+A_U^\perp\pi(a^{(3)})\otimes\k(a^{(2)})a^{(4)}\pi\bigl(
\k(a^{(1)})a^{(5)}\bigr),
\end{split}
\end{equation*}
for each $a\in\widehat{\cal{R}}$. Remembering that $\widehat{\cal{R}}$
is $\ad$-invariant we conclude that the above terms vanish. Hence
\eqref{433} holds. Property \eqref{434} simply follows from
\eqref{433}.\end{pf}

\section{Horizontal Projection, Covariant Derivative $\&$ Curvature}
For each $\omega\in\con(P)$ let $h_\omega\colon\Omega(P)\rightarrow
\Omega(P)$ be a linear map given by
\begin{equation}\label{def-hw}
h_\omega=m_\omega(\id\otimes p_{\inv}^0)m_\omega^{-1}.
\end{equation}
Let $D_\omega\colon\Omega(P)\rightarrow\Omega(P)$ be a linear map
defined as a composition
\begin{equation}\label{def-Dw}
D_\omega=h_\omega d.
\end{equation}
Evidently, both maps are $\hor(P)$-valued.
\begin{defn}\label{def:51}
Operators $h_\omega$ and $D_\omega$ are called {\it the horizontal
projection} and {\it the covariant derivative} associated to $\omega$.
\end{defn}
The following
statement easily follows from the analysis of the previous section.
\begin{pro}\label{pro:51}
\bla{i} The map $h_\omega$ is $\Omega(M)$-linear and
projects $\Omega(P)$ onto $\hor(P)$.

\smallskip
\bla{ii} We have
\begin{equation}\label{51}
(D_\omega-d)\bigl(\Omega(M)\bigr)=\{0\}\qquad
D_\omega(w\varphi)=(dw)h_\omega(\varphi)+(-1)^{\partial w}wD_\omega(
\varphi)
\end{equation}
for each $w\in\Omega(M)$ and $\varphi\in\Omega(P)$.

\smallskip
\bla{iii} Maps $h_\omega$ and $D_\omega$ are
invariant under the action of $G$. In other words, the diagrams
\begin{equation}\label{52}
\begin{CD}
\Omega(P)@>{\mbox{$F^\wedge$}}>> \Omega(P)\otimes\cal{A}\\
@V{\mbox{$h_\omega$}}VV @VV{\mbox{$h_\omega\otimes\id$}}V\\
\Omega(P) @>>{\mbox{$F^\wedge$}}> \Omega(P)\otimes\cal{A}
\end{CD}\qquad\qquad
\begin{CD}
\Omega(P)@>{\mbox{$F^\wedge$}}>> \Omega(P)\otimes\cal{A}\\
@V{\mbox{$D_\omega$}}VV @VV{\mbox{$D_\omega\otimes\id$}}V\\
\Omega(P) @>>{\mbox{$F^\wedge$}}> \Omega(P)\otimes\cal{A}
\end{CD}
\end{equation}
are commutative.

\smallskip
\bla{iv} If $\omega$ is classical then $h_\omega$ is a *-homomorphism.
Furthermore
\begin{equation}\label{53}
D_\omega(\psi\varphi)=D_\omega(\psi)h_\omega(\varphi)+(-1)^{\partial\psi}
h_{\omega}(\psi)D_\omega(\varphi)
\end{equation}
for each $\psi,\varphi\in\Omega(P)$.\qed
\end{pro}

By construction, the space $\hor(P)$ is $D_\omega$-invariant. The
corresponding restriction is described by the following
\begin{pro}\label{pro:52} If $\varphi\in\hor(P)$ then
\begin{equation}\label{54}
D_\omega(\varphi)=d(\varphi)-(-1)^{\partial\varphi}m_\Omega(\id
\otimes\omega\pi)F^\wedge(\varphi).
\end{equation}
In local terms,
\begin{equation}\label{55}
\pi_U^\wedge D_\omega(\varphi)=\sum_i\Bigl\{d(\alpha_i)\otimes a_i
-(-1)^{\partial\alpha}\alpha_iA_U\pi(a_i^{(1)})\otimes a_i^{(2)}\Bigr\},
\end{equation}
where $\Sum_i\alpha_i\otimes a_i=\pi_U^\wedge(\varphi)$.
\end{pro}
\begin{pf}
We have
$$ \pi_U^\wedge d(\varphi)=\sum_id(\alpha_i)\otimes a_i
+(-1)^{\partial\alpha}\alpha_i\otimes a_i^{(1)}\pi(a_i^{(2)}). $$
and hence
\begin{multline*}
\pi_U^\wedge D_\omega(\varphi)=\sum_id(\alpha_i)\otimes a_i
-(-1)^{\partial\alpha}\sum_i\alpha_iA_U\pi(a_i^{(3)})\otimes a_i^{(1)}
\k(a_i^{(2)})a_i^{(4)}\\
=\sum_i\Bigl\{d(\alpha_i)\otimes a_i
-(-1)^{\partial\alpha}\alpha_iA_U\pi(a_i^{(1)})\otimes a_i^{(2)}\Bigr\}
\end{multline*}
according to Definition~\ref{def:51}.
This proves \eqref{55}. Let us compute the right-hand side of
\eqref{54}. We have
\begin{equation*}
\begin{split}
\pi_U^\wedge\bigl[d(\varphi)&-(-1)^{\partial\varphi}m_\Omega(\id
\otimes\omega\pi)F^\wedge(\varphi)\bigr]=\sum_id(\alpha_i)
\otimes a_i\\ &{}+(-1)^{\partial\alpha}\sum_i
\alpha_i\otimes a_i^{(1)}\pi(a_i^{(2)})\\
&{}-(-1)^{\partial\alpha}\sum_i
\alpha_iA_U\pi(a_i^{(1)})\otimes a_i^{(2)}\\
&{}-(-1)^{\partial\alpha}
\sum_i\alpha_i\otimes a_i^{(1)}\pi(a_i^{(2)})\\
=&\sum_i\Bigl\{d(\alpha_i)\otimes a_i
-(-1)^{\partial\alpha}\alpha_iA_U\pi(a_i^{(1)})\otimes a_i^{(2)}\Bigr\}
=\pi_U^\wedge D_\omega(\varphi).\qed
\end{split}
\end{equation*}
\renewcommand{\qed}{}
\end{pf}
For given linear maps
$\alpha,\beta\colon\Gamma_{\inv}\rightarrow\Omega(P)$ we shall denote by
$[\alpha,\beta]$ and $\langle\alpha,\beta\rangle$ linear maps defined by
\begin{align}
[\alpha,\beta]&=m_\Omega(\alpha\otimes\beta)c^\top\label{56}\\
\langle\alpha,\beta\rangle&=m_\Omega(\alpha\otimes\beta)\delta\label{57}
\end{align}
where
$c^\top\colon\Gamma_{\inv}\rightarrow\Gamma_{\inv}\otimes\Gamma_{\inv}$ is
the ``transposed commutator'' map \cite{W3} explicitly given by
\eqref{C15} and $\delta\colon\Gamma_{\inv}\rightarrow
\Gamma_{\inv}\otimes\Gamma_{\inv}$ is the ``embedded differential''
defined by
\begin{equation}\label{emb-d}
\delta(\vartheta)=\iota d(\vartheta).
\end{equation}
If $\alpha,\beta\in\psi(P)$ then
$\langle\alpha,\beta\rangle,[\alpha,\beta]\in\psi(P)$, according to
Lemma~\ref{lem:C6}. In particular these brackets map
$\tau(P)\times\tau(P)$ into $\tau(P)$. Similar brackets can be
introduced for maps valued in an arbitrary algebra.

According to Proposition~\ref{pro:51} the space $\psi(P)$ is mapped,
via compositions with $h_\omega$ and $D_\omega$, into $\tau(P)$. In
particular $\tau(P)$ is $D_\omega$-invariant.
\begin{pro}\label{pro:53} \bla{i} We have
\begin{equation}\label{loc-D}
(\pi_U^\wedge D_\omega\varphi)(\vartheta)=(\Bigl\{d\varphi_U
-(-1)^{\partial\varphi}[\varphi_U,A_U]\Bigr\}\otimes\id)\adj(\vartheta)
\end{equation}
where $\varphi_U$ are local representatives of $\varphi\in\tau(P)$.

\smallskip
\bla{ii} The following identity describes the action of $D_\omega$ on
tensorial forms
\begin{equation}\label{58}
D_\omega\varphi=d\varphi-(-1)^{\partial\varphi}[\varphi,\omega].
\end{equation}
\end{pro}
\begin{pf}
We have
$$ \pi_U^\wedge d\varphi(\vartheta)=\sum_kd\varphi_U(\vartheta_k)\otimes
c_k+(-1)^{\partial\varphi}\varphi_U(\vartheta_k)\otimes
c_k^{(1)}\pi(c_k^{(2)})$$
where $\Sum_k\vartheta_k\otimes c_k=\adj(\vartheta)$. Taking the
horizontal projection we obtain
\begin{equation*}
\begin{split}
(\pi_U^\wedge D_\omega\varphi)(\vartheta)&=\sum_k\Bigl\{
d\varphi_U(\vartheta_k)\otimes c_k-(-1)^{\partial\varphi}
\varphi_U(\vartheta_k)A_U\pi(c_k^{(3)})\otimes c_k^{(1)}\k(c_k^{(2)})
c_k^{(4)}\Bigr\}\\
&=\sum_k
d\varphi_U(\vartheta_k)\otimes c_k-(-1)^{\partial\varphi}
\varphi_U(\vartheta_k)A_U\pi(c_k^{(1)})\otimes c_k^{(2)}\\
&=\sum_k \bigl(d\varphi_U
-(-1)^{\partial\varphi}[\varphi_U,A_U]\bigr)(\vartheta_k)\otimes c_k.
\end{split}
\end{equation*}
A computation of the right-hand side of \eqref{58} gives
\begin{multline*}
\pi_U^\wedge\Bigl\{d\varphi-(-1)^{\partial\varphi}[\varphi,\omega]\Bigr\}
=\sum_k
d\varphi_U(\vartheta_k)\otimes c_k-(-1)^{\partial\varphi}
\varphi_U(\vartheta_k)A_U\pi(c_k^{(1)})\otimes c_k^{(2)}\\
=(\pi_U^\wedge D_\omega\varphi)(\vartheta).\qed
\end{multline*}
\renewcommand{\qed}{}
\end{pf}
Let $q_\omega\colon\psi(P)\rightarrow\psi(P)$ be a linear map defined by
\begin{equation}\label{59}
q_\omega(\varphi)=\langle\omega,\varphi\rangle-(-1)^{\partial\varphi}
\langle\varphi,\omega\rangle-(-1)^{\partial\varphi}[\varphi,\omega].
\end{equation}
By definition, this map is $\Omega(M)$-linear  from the right.
\begin{pro}\label{pro:54} The space $\tau(P)$ is $q_\omega$-invariant.
\end{pro}
\begin{pf}
For a given $\vartheta\in\Gamma_{\inv}$ let us choose $a\in\ker(\e)$
satisfying conditions listed
in Lemma~\ref{lem:C7} ({\it i\/}). We have then
\begin{equation*}
\begin{split}
-(-1)^{\partial\varphi}\bigl(\pi_U^\wedge q_\omega(\varphi)\bigr)(\vartheta)&=
\sum_k\Bigl\{\varphi_U(\vartheta_k)A_U\pi(c_k^{(1)})\otimes c_k^{(2)}+
\varphi_U(\vartheta_k)\otimes c_k^{(1)}\pi(c_k^{(2)})\Bigr\}\\
&\phantom{=}-\varphi_U\pi(a^{(2)})A_U\pi(a^{(3)})\otimes\k(a^{(1)})a^{(4)}\\
&\phantom{=}-\varphi_U\pi(a^{(2)})\otimes\k(a^{(1)})a^{(3)}\pi(a^{(4)})\\
&\phantom{=}+(-1)^{\partial\varphi}A_U\pi(a^{(2)})\varphi_U\pi(a^{(3)})\otimes
\k(a^{(1)})a^{(4)}\\
&\phantom{=}+\varphi_U\pi(a^{(3)})\otimes\pi(a^{(1)})\k(a^{(2)})a^{(4)},
\end{split}
\end{equation*}
for each $\varphi\in\tau(P)$.

On the other hand, applying \eqref{B40} and
\eqref{B21} we find
\begin{multline*}
\varphi_U\pi(a^{(2)})\otimes\k(a^{(1)})a^{(3)}\pi(a^{(4)})-\varphi_U
\pi(a^{(3)})\otimes\pi(a^{(1)})\k(a^{(2)})a^{(4)}\\
=\sum_k\varphi_U(\vartheta_k)\otimes c_k^{(1)}\pi(c_k^{(2)}).
\end{multline*}
Combining the above equalities we obtain finally
\begin{equation}\label{loc-qw}
\bigl(\pi_U^\wedge q_\omega(\varphi)\bigr)(\vartheta)=\bigl(
q_\omega^U(\varphi)\otimes\id\bigr)\adj(\vartheta)
\end{equation}
where
\begin{equation}\label{512}
q_\omega^U(\varphi)=\langle A_U,\varphi_U\rangle
-(-1)^{\partial\varphi}\langle\varphi_U,A_U\rangle-(-1)^{\partial\varphi}
[\varphi_U,A_U].
\end{equation}
We see that $q_\omega(\varphi)$ is tensorial.
\end{pf}

If $\omega$ is classical then the operator $q_\omega$ vanishes on
tensorial forms. Indeed, in this case
$$ A_U\pi(ab)=\e(a)A_U\pi(b)+\e(b)A_U\pi(a)$$
which, together with \eqref{56}--\eqref{57}, implies
\begin{equation*}\begin{split}
[\varphi_U,A_U](\vartheta)&=\varphi_U\pi(a^{(2)}) A_U\pi
\bigl(\k(a^{(1)})a^{(3)}\bigr)\\
&=\varphi_U\pi(a^{(1)})A_U\pi(a^{(2)})-(-1)^{\partial\varphi}
A_U\pi(a^{(1)})\varphi_U\pi(a^{(2)})\\
&=-\bigl(\langle\varphi_U,A_U\rangle-(-1)^{\partial\varphi}\langle
A_U,\varphi_U\rangle\bigr)(\vartheta).
\end{split}\end{equation*}
Consequently, in the general case the operator
$q_\omega{\restr}\tau(P)$ depends only
on the quantum part $\omega_\perp$ of $\omega$, and can be written in an
explicitly tensorial form
\begin{equation}\label{513}\begin{split}
q_\omega(\varphi)&=\langle\omega_\perp,\varphi\rangle-(-1)^{\partial\varphi}
\langle\varphi,\omega_\perp\rangle-(-1)^{\partial\varphi}[\varphi,
\omega_\perp]\\
q_\omega^U(\varphi)&=\langle A_U^\perp,\varphi_U\rangle
-(-1)^{\partial\varphi}\langle\varphi_U,A_U^\perp\rangle-
(-1)^{\partial\varphi}
[\varphi_U,A_U^\perp].
\end{split}\end{equation}

The rest of the section is devoted to the introduction and the analysis
of the curvature form.
\begin{defn}\label{def:52}
A tensorial 2-form
\begin{equation}\label{def-Rw}
R_\omega=D_\omega\omega
\end{equation}
is called {\it the curvature} of $\omega$.
\end{defn}

This definition directly follows classical differential geometry.
However, in contrast to the classical case, the curvature is generally
$\delta$-dependent.
\begin{pro}\label{pro:55}
We have
\begin{equation}\label{516}
\pi_U^\wedge R_\omega(\vartheta)=(F_U\otimes\id)\adj(\vartheta)
\end{equation}
where
\begin{equation}\label{517}
F_U=dA_U-\langle A_U,A_U\rangle.
\end{equation}
\end{pro}
\begin{pf}
A direct calculation gives
\begin{equation*}
\begin{split}
-\bigl(\pi_U^\wedge\omega^\wedge\bigr)(d\vartheta)&=1_U\otimes\pi(a^{(1)})
\pi(a^{(2)})
+A_U\pi(a^{(2)})\otimes\k(a^{(1)})a^{(3)}\pi(a^{(4)})\\
&\phantom{=}-A_U\pi(a^{(3)})\otimes\pi(a^{(1)})\k(a^{(2)})a^{(4)}\\
&\phantom{=}+A_U\pi(a^{(2)})A_U\pi(a^{(3)})\otimes\k(a^{(1)})a^{(4)}\\
=1_U\otimes\pi(a^{(1)})\pi(a^{(2)})&+\sum_k\Bigl\{A_U(\vartheta_k)
\otimes c_k^{(1)}\pi(c_k^{(2)})-\langle A_U,A_U\rangle(\vartheta_k)
\otimes c_k\Bigr\}.
\end{split}
\end{equation*}
On the other hand
$$\bigl(\pi_U^\wedge d\omega\bigr)(\vartheta)= -1_U\otimes\pi(a^{(1)})
\pi(a^{(2)})+\sum_k\Bigl\{dA_U(\vartheta_k)\otimes c_k -A_U(\vartheta_k)
\otimes c_k^{(1)}\pi(c_k^{(2)})\Bigr\}.$$
Here $\adj(\vartheta)=\Sum_k\vartheta_k\otimes c_k$ and $a\in\ker(\e)$ is
chosen as explained in Lemma~\ref{lem:C7}.

Combining the above expressions we find
\begin{equation}\label{519}
\pi_U\bigl(d\omega(\vartheta)-
\omega^\wedge(d\vartheta)\bigr)=
\sum_k\Bigl\{dA_U(\vartheta_k)\otimes c_k-\langle
A_U,A_U\rangle(\vartheta_k)\otimes c_k\Bigr\}.
\end{equation}

To complete the proof it is sufficient to observe that last two summands
in the right-hand side of the
above equation are horizontal while the first one is completely
``vertical''.
\end{pf}

Now, the analogs of classical Structure Equation and Bianchi identity
will be derived.
\begin{pro}\label{pro:56} The following identities hold
\begin{align}
R_\omega&=d\omega-\langle\omega,\omega\rangle\label{520}\\
D_\omega R_\omega-q_\omega(R_\omega)&=
\langle\omega_\perp,\langle\omega_\perp,\omega_\perp\rangle\rangle
-\langle\langle\omega_\perp,\omega_\perp,\rangle,\omega_\perp\rangle.
\label{521}
\end{align}
\end{pro}
\begin{pf} The previous proposition and equation \eqref{519} imply
$$(\pi_U^\wedge d\omega)(\vartheta)=\bigl(\pi_U^\wedge
\omega^\wedge\bigr)(d\vartheta)+\bigl(\pi_U^\wedge
R_\omega\bigr)(\vartheta)=\bigl\{\pi_U^\wedge(R_\omega+\langle\omega,\omega
\rangle)\bigr\}(\vartheta),$$
for each $\vartheta\in\Gamma_{\inv}$ and $U\in\cal{U}$.
Hence \eqref{520} holds.

Equation \eqref{512} and Proposition~\ref{pro:53} \bla{ii} imply
\begin{multline*}
\bigl[\pi_U^\wedge\bigl(D_\omega
R_\omega -q_\omega(R_\omega)\bigr)\bigr](\vartheta)=\sum_kddA_U
(\vartheta_k)\otimes c_k-\sum_k\langle dA_U,A_U\rangle(\vartheta_k)\otimes
c_k\\ {}+\sum_k\langle A_U,dA_U\rangle(\vartheta_k)\otimes
c_k+\sum_k\Bigl\{
\langle F_U,A_U\rangle(\vartheta_k)\otimes c_k-\langle A_U,F_U\rangle(
\vartheta_k)\otimes c_k\Bigr\}\\
=\sum_k\bigl(\langle A_U,\langle A_U,A_U\rangle\rangle-
\langle\langle A_U,A_U\rangle,A_U\rangle\bigr)(\vartheta_k)\otimes c_k.
\end{multline*}

On the other hand, using Lemma~\ref{lem:C8} we conclude that
$$
\langle A_U,\langle A_U,A_U\rangle\rangle-
\langle\langle A_U,A_U\rangle,A_U\rangle
=\langle A_U^\perp,\langle A_U^\perp,A_U^\perp\rangle\rangle-
\langle\langle A_U^\perp,A_U^\perp\rangle,A_U^\perp\rangle.
$$
This is the local expression for the right-hand side of \eqref{521}.
\end{pf}

If $\omega$ is classical then \eqref{520}--\eqref{521} are equivalent to
classical Structure Equation and Bianchi identity for $\omega$, if
$\omega$ is understood as a (standard) connection on $P_{\!cl}$.

More generally, if $\omega$ is multiplicative then the right-hand side
of \eqref{521} vanishes. Indeed in this case we have
$$\langle\omega_\perp,\omega_\perp\rangle\pi(a)=
-\omega_\perp\pi(a^{(1)})\omega_\perp\pi(a^{(2)})$$
for each $a\in\cal{A}$.

It is important to mention that the proofs of identities contained in
Propositions~\ref{pro:54}--\ref{pro:56}, the choose of an embedding
$\iota$ figures only via its restriction on $d(\Gamma_{\inv})$, which
determines the embedded differential map $\delta$.

Generally, a map $\delta$ can be constructed by fixing a $*\k$-invariant
$\ad$-invariant complement $\cal{L}\subseteq\ker(\e)$ of
$\widehat{\cal{R}}$, and defining
\begin{equation}\label{522}
-\delta=(\pi\otimes\pi)\phi(\pi{\restr}\cal{L})^{-1}.
\end{equation}
If, in addition, $\phi(\cal{L})\subseteq 1\otimes\cal{L}+
\cal{L}\otimes 1+\cal{L}\otimes\cal{L}$ then the above $\delta$
satisfies
$$(\delta\otimes\id)\delta=(\id\otimes\delta)\delta$$
and right-hand side of \eqref{521} vanishes identically.

Our restriction to the minimal admissible left-covariant calculus
$\Gamma$ is not essential. All considerations can be performed using an
arbitrary admissible bicovariant *-calculus. Moreover, if the bundle is
trivial we can abandon the assumption of admissibility, and work in a
fixed global trivialization.

For example if we take $\cal{R}=\{0\}$ then $\Gamma$ becomes the
``maximal'' calculus. In this case $\Gamma_{\inv}=\ker(\e)$ and
$\Gamma^\wedge=\Gamma^\otimes$ is the universal differential envelope of
$\cal{A}$ (modulo the relation $d1=0$). Because of $S^\wedge=\{0\}$, every
connection is multiplicative and $\delta$ is uniquely determined.
\section{Examples}
In this section we consider some illustrative examples related to the
presented theory. We shall discuss ``nonclassical'' phenomena appearing
in the formalism of connections, as well as interesting properties of
appropriate differential caluli over the structure group $G$.

Two types of $G$ will be considered. The case of a classical Lie group
$G$, and the quantum case $G=S_\mu U(2)$.

As a possible application in theoretical physics, we shall briefly
describe a ``gauge theory'' based on quantum principal bundles.

\subsection{Classical Structure Groups}
Let us assume that $G$ is a classical compact Lie group ($\cal{A}$ is
commutative and $G_{\!cl}=G$). The corresponding
principal bundles are objects of classical differential geometry.

The minimal admissible calculus over $G$ coincides with the classical
one, based on standard 1-forms. The corresponding universal differential
envelope gives the classical higher-order calculus on $G$, based on
standard differential forms.

The classical calculus on $G$, together with the classical calculus on
the base manifold $M$, induces the classical differential calculus on
corresponding principal bundles. The whole theory presented in this
paper is equivalent to the classical theory.

However, if we start from a {\it nonstandard\/} differential calculus on
$G$ then, generally, ``quantum phenomena'' will enter the game.

Let $\Gamma$ be an arbitrary admissible bicovariant *-calculus over $G$,
and let $\cal{R}\subseteq\ker(\e)$ be the corresponding $\cal{A}$-ideal.
We have
$$\cal{R}\subseteq\ker(\e)^2$$
because of the admissibility of $\Gamma$.

For example, if $\cal{R}=\ker(\e)^k$ with $k\geq 2$,
then $\Gamma_{\inv}$ is naturally isomorphic
to the space of $(k-1)$-jets in the neutral element $\e\in G$.

Let $P$ be a principal $G$-bundle over $M$ and $\omega\in\con(P)$.
After choosing a splitting \eqref{411} the ``classical-quantum''
decomposition of $\omega$ can be performed. Components of the field
$\omega_\perp$ are ``labeled'' by elements of the space $\ker(\nu)$. The
field $\omega_\perp$ figures in ``quantum terms'' introduced in previous
two sections. Generally these terms do not vanish. Moreover they
already figure in the case of a {\it finite\/} group $G$.

\subsection{The Minimal Admissible Calculus For Quantum
{\it S\!U}{\rm (}{\it 2\/}{\rm )}}

This subsection is devoted to the analysis of the minimal admissible
left-covariant calculus $\Gamma$ over the group $G=S_\mu U(2)$. We shall
also briefly discuss certain features of corresponding principal
bundles.

As  first, let us assume that $\mu\in(-1,1)\setminus\{0\}$. As explained
in Appendix A, $G_{\!cl}=U(1)$ in a natural manner. The (complex)
Lie algebra of
$G_{\!cl}$ is spanned by a single element
$X\colon\cal{A}\rightarrow\Bbb{C}$ determined by
\begin{equation}\label{61}
X(\alpha)=-X(\alpha^*)=\frac{1}{2}\qquad X(\gamma)=X(\gamma^*)=0.
\end{equation}
The correspondence $X\leftrightarrow 1$ enables us to identify
$\lie(G_{\!cl})=\Bbb{C}$. In particular, the space $\Gamma_{\inv}$ can be
viewed (via the map $\rho$) as a certain subspace of $\cal{A}$.
\begin{pro}\label{pro:61}
The map $\rho\colon\Gamma_{\inv}\rightarrow\cal{A}$ is a bijection onto
the subalgebra $\cal{Q}\subseteq\cal{A}$ consisting of left
$U(1)$-invariant elements. A natural basis in $\cal{Q}$ is given by elements
$\xi_{n,k}$ where $n\in\Bbb{Z}$ and $k\in\Bbb{N}\cup\{0\}$ and
\begin{equation}\label{62}
\xi_{n,k}=\begin{cases}(-\mu)^n(\gamma\gamma^*)^k\gamma^n\alpha^n &
\mbox{if $n\geq 0$}\\
(\alpha^*)^{-n}(\gamma^*)^{-n}(\gamma^*\gamma)^k& \mbox{if $n\leq 0$}
\end{cases}
\end{equation}
\end{pro}
\begin{pf} According to \cite{W1} the elements
$\alpha^n\gamma^k\gamma^{*r}$ form a basis in $\cal{A}$ (by definition
$\alpha^{-n}=\alpha^{*n}$). It is easy to see that $g\in
U(1)$ acts on the left by multiplying these elements by $z^{n-k+r}$, where
$z=g(\alpha)$. Hence, $\cal{Q}$ is spanned by basis elements satisfying
$n-k+r=0$. Equivalently, elements \eqref{62} form a basis in $\cal{Q}$.

We have to verify that $\cal{Q}=\rho\bigl(\Gamma_{\inv}\bigr)$. According
to Lemma~\ref{lem:C10} ({\it i\/}) the image of $\rho$ is contained in
$\cal{Q}$. It is easy to see that
\begin{equation}\label{64}
\begin{aligned}
\rho\pi(\alpha)&=\frac{1}{2}-\gamma\gamma^*\\
\rho\pi(\alpha^*)&=\mu^2\gamma\gamma^*-\frac{1}{2}
\end{aligned}\qquad
\begin{gathered}
\rho\pi(\gamma^*)=\alpha^*\gamma^*\\
\rho\pi(\gamma)=-\alpha\gamma.
\end{gathered}
\end{equation}
Furthermore, a straightforward calculation gives
\begin{gather}
\xi_{n,k}\circ\alpha^{\phantom{*}}=\mu^{-2k-\vert n\vert}\xi_{n,k}+
(\mu^{\vert n\vert}-\mu^{-2k-\vert n\vert})\xi_{n,k+1}\label{mod1}\\
\xi_{n,k}\circ\alpha^*=\mu^{2k+\vert n\vert}\xi_{n,k}+
\mu^2(\mu^{\vert n\vert}-\mu^{2k+3\vert n\vert})\xi_{n,k+1}\label{mod2}\\
\xi_{n,k}\circ\gamma^{\phantom{*}}=(1-\mu^{2(k+n)})\xi_{n+1,k},\quad
n\geq 0\label{mod3}\\
\xi_{n,k}\circ\gamma^*=(1-\mu^{2(k-n)})\xi_{n-1,k},\quad n\leq
0\label{mod4}\\
\xi_{n,k}\circ\gamma^{\phantom{*}}=(1-\mu^{-2k})\xi_{n+1,k+1}+
\mu^{-2k}(1-\mu^{2(k-n)})\xi_{n+1,k+2}, \quad n<0\label{mod5}\\
\xi_{n,k}\circ\gamma^*=(1-\mu^{-2k})\xi_{n-1,k+1}+
\mu^{-2k}(1-\mu^{2(k+n)})\xi_{n-1,k+2},\quad n>0\label{mod6}
\end{gather}
The $\circ$ operation is given by
$\xi\circ a=\k(a^{(1)})\xi a^{(2)}$. We see that $\cal{Q}$ is invariant
under $\circ$. Above formulas imply that $\cal{Q}$ is generated, as a
right $\cal{A}$-module, by elements \eqref{64}. Having in mind that
$\rho(\Gamma_{\inv})$ is a right $\cal{A}$-submodule of $\cal{Q}$ (as
follows from \eqref{C5}) we conclude that $\rho$ is surjective.
\end{pf}

The following proposition describes the right $\cal{A}$-ideal
$\widehat{\cal{R}}$ corresponding to the calculus $\Gamma$.
\begin{pro}\label{pro:62}
We have
\begin{equation}\label{67}
\widehat{\cal{R}}=\bigl(\mu^2\alpha+\alpha^*-(1+\mu^2)1\bigr)\ker(\e).
\end{equation}
\end{pro}
\begin{pf} Let $\cal{R}$ be the right-hand side of \eqref{67}. According
to Lemma~\ref{lem:C10} ({\it ii\/}) the space $\cal{R}$ is contained in
$\widehat{\cal{R}}$.

On the other hand, the space of
$\ad$-invariant elements of $\cal{A}$ consists precisely
of polynomials of $\mu^2\alpha+\alpha^*$ and we have
$$\ad(ba)=b\ad(a)$$
for each $a\in\cal{A}$ and an $\ad$-invariant element $b\in\cal{A}$. In
particular, corresponding
multiple irreducible subspaces are closed under
the left multiplication by $\ad$-invariant elements. Furthermore,
primitive elements for nonsinglet multiple irreducible subspaces of
$\ad$ are of the form $p(\mu^2\alpha+\alpha^*)\gamma^k$ and
$p(\mu^2\alpha+\alpha^*)\gamma^{*k}$, corresponding to spin $k$ highest
and lowest weights respectively. Hence, in the decomposition of the
factorized adjoint action on $\ker(\e)/\cal{R}$ each irreducible
multiplet appears no more than once. On the other hand, elements
$\rho\pi(\gamma^n)$, $\rho\pi(\gamma^{*n})$ and $\rho\pi(\mu^2\alpha+
\alpha^*)$ are all non-zero (as follows from \eqref{64},
\eqref{mod3}--\eqref{mod4} and
\eqref{C5}). Therefore, for each spin value, the representation $\ad$
contains at least one irreducible multiplet. Consequently $\cal{R}=
\widehat{\cal{R}}$.
\end{pf}

We pass to the detailed analysis of the adjoint action $\adj$. In terms
of the identification $\Gamma_{\inv}=\cal{Q}$ we have
$$\adj=(\phi{\restr}\cal{Q}).$$ Let us assume that $\Gamma_{\inv}$ is
endowed with a natural $\adj$-invariant scalar product, induced by the
Haar measure (as explained in Appendix C). We are going to decompose
$\adj$ into irreducible multiplets. Let us consider operators
\begin{equation}\label{614}
K_\pm=(\id\otimes X_\pm)\adj\qquad K_3=(\id\otimes X)\adj
\end{equation}
which are counterparts for the ``creation'' and ``anihilation'', as well
as the ``third spin component'' operator. Here
$X_\pm\colon\cal{A}\rightarrow\Bbb{C}$ are linear functionals satisfying
\begin{equation}\label{615}
X_\pm(ab)=X_\pm(a)\chi(b)+\e(a)X_\pm(b)
\end{equation}
where $\chi\colon\cal{A}\rightarrow\Bbb{C}$ is a multiplicative functional
determined by
$$\chi(\alpha)=\frac{1}{\mu}\qquad \chi(\alpha^*)=\mu\quad \chi(\gamma)=
\chi(\gamma^*)=0.$$
We shall addopt the following normalization
$$ X_\pm(\alpha)=X_\pm(\alpha^*)=X_+(\gamma)=X_-(\gamma^*)=0
\quad -\mu X_+(\gamma^*)=X_-(\gamma)=1. $$

It turns out that the following identities hold
\begin{gather}
K_+K_--\mu^2K_-K_+=\frac{1-\mu^{-4K_3}}{1-\mu^{-2\phantom{K_3}}}\label{619}\\
K_3K_+-K_+K_3=K_+\qquad K_3K_--K_-K_3=-K_-\label{617}\\
K_3(\vartheta\eta)=K_3(\vartheta)\eta+\vartheta K_3(\eta)\label{620}\\
K_\pm(\vartheta\eta)=K_\pm(\vartheta)\chi_\adj(\eta)+\vartheta
K_\pm(\eta)\label{621}
\end{gather}
where $\chi_\adj=(\id\otimes \chi)\adj$. Furthermore, we have
\begin{gather}
\chi_\adj(\xi_{n,k})=\mu^{-2n}\xi_{n,k}\qquad
K_3(\xi_{n,k})=n\xi_{n,k}\label{623}\\
\begin{aligned}
K_+(\xi_{n,k})&=\frac{1-\mu^{2k}}{\mu^{n+2}(1-\mu^2)}\xi_{n+1,k-1}\\
K_-(\xi_{n,k})&=\frac{\mu^{1-n}(1-\mu^{2k})}{1-\mu^2}\xi_{n-1,k-1}
\end{aligned}\qquad
\begin{aligned}
\phantom{\frac{1}{\mu}}n\geq& 0\\
\phantom{\frac{1}{\mu}}n\leq& 0.
\end{aligned}\label{625}
\end{gather}

Now \eqref{623}--\eqref{625} imply that
$$\cal{Q}=\sideset{}{^\oplus}\sum_{k\geq 0}\cal{Q}_k,$$
where $\cal{Q}_k$ are irreducible subspaces for the $k$-spin
representation. In particular
\begin{equation}\label{627}
\cal{Q}_k=\sideset{}{^\oplus}\sum_{\vert m\vert\leq k}\cal{Q}_{k,m}
\end{equation}
where $\cal{Q}_{k,m}=\ker(mI-K_3)\cap\cal{Q}_k$. The spaces
$\cal{Q}_{k,m}$ are 1-dimensional. Hence it is possible
to construct an orthonormal basis
in $\cal{Q}$ by choosing unit vectors $\zeta_{k,m}\in\cal{Q}_{k,m}$. A
priori, there exists an ambiguity for this choice, one phase factor for
each $\zeta_{k,m}$. However, requiring that non-vanishing matrix
elements of $K_\pm$ are positive, the ambiguity is reduced to one phase
factor for each multiplet. According to \cite{W1}, we have
\begin{equation}
K_+\zeta_{k,m}=v_{k,m+1}\zeta_{k,m+1}\qquad
K_-\zeta_{k,m}=v_{k,m}\zeta_{k,m-1}\label{628}
\end{equation}
where
$$ v_{k,m}=\mu^{1-m-k}\bigl((k+m)_\mu(k-m+1)_\mu\bigr)^{1/2}
\qquad n_\mu=\frac{1-\mu^{2n}}{1-\mu^{2\phantom{n}}}.$$

Let $\cal{P}$ be the space of one-variable polynomials. It is easy to
see that
\begin{equation}\label{630}
\zeta_{0,k}=p_k(\gamma^*\gamma)
\end{equation}
where  $p_k\in\cal{P}$ are $k$-th order polinomials
orthonormal with respect to a scalar
product given by
\begin{equation}\label{631}
(p,q)=\int\!p^*q.
\end{equation}
Here $\int\colon\cal{P}\rightarrow\Bbb{C}$ is a linear functional
given by
\begin{equation}\label{637}
\int \! x^n=(n+1)_\mu^{-1}.
\end{equation}
We shall assume that leading coefficients of polynomials $p_k$ are
positive. This completely fixes vectors $\zeta_{k,m}$.
\begin{pro}\label{pro:63} \bla{i} Polynomials $p_k$ are given by
\begin{equation}\label{632}
p_k(x)=(-1)^kc_k\partial^k\Bigl[x^k\prod_{j=1}^k(1-\mu^{1-j}x)\Bigr]
\end{equation}
where $c_k>0$ are normalization constants and $\partial\colon\cal{P}
\rightarrow\cal{P}$ is a linear map specified by
\begin{equation}\label{633}
\partial(x^n)=n_\mu x^{n-1}.
\end{equation}

\bla{ii} The following identities hold
\begin{equation}\label{634}
\begin{aligned}
\zeta_{k,m}=(-1)^m\mu^{km-m}&\biggl(\frac{(k-m)_\mu!}{(k+m)_\mu!}
\biggr)^{1/2}(\partial^mp_k)(\gamma\gamma^*)\gamma^m\alpha^m\\
\zeta_{k,-m}=\mu^{km}\alpha^{*m}\gamma^{*m}&\biggl(
\frac{(k-m)_\mu!}{(k+m)_\mu!}\biggr)^{1/2}
(\partial^mp_k)(\gamma\gamma^*)
\end{aligned}
\end{equation}
where $m\in\{0,\dots,k\}$ and $n_\mu!=\displaystyle{\prod}_{j=1}^nj_\mu$.
\end{pro}
\begin{pf}
The map $\partial$ satisfies the following ``Leibniz rule''
\begin{equation}\label{635}
\partial(pq)(x)=(\partial p)(x)q(x)+p(\mu^2 x)(\partial q)(x),
\end{equation}
as directly follows from \eqref{633}. More generally
\begin{equation}\label{636}
\partial^n(pq)(x)=\sum_{k=0}^n \binom{n}{k}_{\!\!\mu}(\partial^{n-k}p)
(\mu^{2k}x)(\partial^kq)(x)
\end{equation}
for each $n\in\Bbb{N}$. In the above formula
$$\binom{n}{k}_{\!\!\mu}=\frac{n_\mu!}{k_\mu!(n-k)_\mu!}.$$

It is easy to see that
\begin{equation}\label{639}
\int\partial(p)=p(1)-p(0)
\end{equation}
for each $p\in\cal{P}$. Inductively using \eqref{635} and \eqref{639} we
obtain the following ``partial integration'' rule
\begin{multline*}
\int \!q\partial^n(p)=\sum_{k=1}^n(-1)^{k-1}\mu^{-k(k-1)}
\Bigl\{(\partial^{n-k}p)
(\mu^{2k-2}x)(\partial^{k-1}q)(x)\Big\vert_0^1\Bigr\}\\
{}+(-1)^n\mu^{-n(n-1)}\int\!(\partial^nq)p(\mu^{2n}x).
\end{multline*}

It is now easy to prove that polynomials $p_k$ given by \eqref{632} are
mutually orthogonal. Furthermore, leading coefficients of these
polynomials are positive. Having in mind that $p_k$ are normed we
conclude that \eqref{630} holds.

To prove ({\it ii\/}) it is sufficient to act by $K_\pm^m$ on both sides
of \eqref{630}, and to apply \eqref{625} and \eqref{628}.
\end{pf}

It is worth noticing that $\cal{Q}$ is *-invariant. The map
$*\colon\cal{Q}\rightarrow\cal{Q}$ corresponds to the canonical
*-structure on $\Gamma_{\inv}$. We have
$$\zeta_{k,-m}^*=(-\mu)^m\zeta_{k,m}.$$

In the classical limit the algebra $\cal{A}$ consists of polynomial
functions on the group $SU(2)$. The subalgebra $\cal{Q}$ then consists
of polynomial functions invariant under left translations by diagonal
matrices from $U(1)$. Equivalently, $\cal{Q}$ can be described as the
algebra of polynomial functions on the 2-sphere $S^2$, because the above
mentioned action defines the Hopf fibering $S^3\rightarrow S^2$. In this
picture $\zeta_{k,m}$ become spherical harmonics, and $K_3$, $K_\pm$
correspond to standard angular momentum operators.

Of course, for $\mu=1$ the minimal admissible calculus is just the classical
$3$-dimensional one. As we shall see later, a similar situation holds for
$\mu=-1$.

In the general case the algebra $\cal{Q}$ represents polynomial
functions on a ``quantum 2-sphere'' \cite{P}. At the level of spaces,
the inclusion $\cal{Q}\hookrightarrow\cal{A}$ is interpretable as the
``quantum Hopf fibering''.
\begin{pro}\label{pro:64}
The space $S_{\inv}^{\wedge 2}$ consists precisely of elements of the form
\begin{multline*}
q=1\otimes\Bigl[\pi(b)+\frac{2\mu^2}{1-\mu^2}(\gamma^*\gamma)\circ
b\Bigr]+\Bigl[\pi(b)+\frac{2\mu^2}{1-\mu^2}(\gamma^*\gamma)\circ
b\Bigr]\otimes 1\\
{}-\frac{2\mu^2}{1-\mu^2}\Bigl[(1+\mu^2)\gamma\gamma^*\otimes\gamma\gamma^*
+\mu\alpha^*\gamma^*\otimes\alpha\gamma+\frac{1}{\mu}\alpha\gamma
\otimes\alpha^*\gamma^*\Bigr]({\circ}\otimes{\circ})\phi(b),
\end{multline*}
where $b\in\ker(\e)$.
\end{pro}
\begin{pf}
The statement follows from Lemma~\ref{lem:B11}, Proposition~\ref{pro:62},
and properties \eqref{B40} and \eqref{64}.
\end{pf}

Let us now consider a quantum principal $G$-bundle $P$ over a compact
manifold $M$. According to the results of Section 2, the structure of
$P$ is completely determined by its classical part $P_{\!cl}$, which is a
classical $U(1)$-bundle over $M$. Let us consider a connection $\omega$,
and describe its components $\omega_{\!cl}$ and $\omega_\perp$. As first,
we have to specify a splitting \eqref{411}. Modulo the identification
$\Gamma_{\inv}=\cal{Q}$ we have $\nu=(\e{\restr}\cal{Q})$. With the help of
$\nu$, let us identify $\L^*$ with the 1-dimensional subspace in
$\Gamma_{\inv}$ generated by $1$. The elements of the
subspace $\L^*$ are characterized by $\xi\circ a=\e(a)\xi$.

Therefore, the classical
component $\omega_{\!cl}$ is locally determined by 1-form
$A_U(1)$. From the point of view of classical geometry, this
1-form is a gauge potential of $\omega_{\!cl}$, understood as a connection
on $P_{\!cl}$. On the other hand, the quantum component $\omega_\perp$ is
locally determined by a collection of 1-forms $A_U(\xi_{n,k})$, where
$(n,k)\neq (0,0)$. Globally, we have a collection of
tensorial 1-forms on $P_{\!cl}$.

It is important to mention that such a classical reinterpretation of
connections destroys the information about irreducible multiplets
structure of corersponding gauge potentials. Because of mutual {\it
incompatibility} of decompositions \eqref{411} and \eqref{627}.

Let us now describe  a construction of the embedded differential map
$\delta$. In the context of this example, $\delta$ can be naturally
introduced with the help of a splitting
$\ker(\e)=\widehat{\cal{R}}\oplus\cal{L}$, where
$\cal{L}\subseteq\ker(\e)$ is the minimal $\ad$-invariant lineal which
contains $\mu^2\alpha+\alpha^*-(1+\mu^2)1$ and $\gamma^k$, for each
$k\in\Bbb{N}$. Explicitly, this lineal can be constructed by extracting
irreducible multiplets from $\ad(\gamma^k)$. The map $\delta$ is given
by \eqref{522}.

According to \eqref{517} the local expression for the curvature is given
by
$$ F_U\pi(a)=A_U\pi(a)+A_U\pi(a^{(1)})A_U\pi(a^{(2)}),$$
where $a\in\cal{L}$.

Let us consider the case $\mu=-1$. As explained in Appendix A, the
classical part of $G$ is isomorphic to a semidirect product of groups
$U(1)$ and $\Bbb{Z}_2=\{-1,1\}$. The corresponding
Lie algebra is generated by a
single element $X$, as in the previous example. Let $\Gamma$ be the
minimal admissible left-covariant calculus. Equations \eqref{64}
reduce to
\begin{equation}\label{644}
\begin{gathered}
\rho\pi(\gamma^*)=\alpha^*\gamma^*\quad\rho\pi(\gamma)=\gamma\alpha\\
\rho\pi(\alpha)=-\rho\pi(\alpha^*)=\frac{1}{2}-\gamma\gamma^*.
\end{gathered}
\end{equation}

The $\circ$-structure is given by
\begin{equation}\label{645}
\begin{gathered}
\pi(\gamma)\circ\{\alpha,\alpha^*\}=-\pi(\gamma)\quad\pi(\gamma^*)\circ
\{\alpha,\alpha^*\}=-\pi(\gamma^*)\\
\pi\{\gamma,\gamma^*\}\circ\{\gamma,\gamma^*\}=\{0\}\quad
\pi(\alpha)\circ a=\e(a)\pi(\alpha).
\end{gathered}
\end{equation}
Consequently, elements
$$\eta_+=\pi(\gamma)\quad\eta_3=\pi(\alpha-\alpha^*)\quad\eta_-=\pi
(\gamma^*)$$
form a basis in $\Gamma_{\inv}$.

{}From \eqref{645} and Lemma~\ref{lem:B15} ({\it i\/}) it follows that the
flip-over operator $\sigma$ is just the standard transposition.
Furthermore, the space $S_{\inv}^\wedge$ is consisting precisely of symmetric
elements of $\Gamma_{\inv}^{\otimes 2}$.

It is worth noticing that the map $\delta$ is uniquely determined,
because $\Gamma_{\inv}^{\otimes 2}$ contains only one irreducible
triplet. Explicitly,
\begin{equation}\label{647}
\begin{aligned}
\delta(\eta_+)&=(\eta_3\otimes\eta_+-\eta_+\otimes\eta_3)/2\\
\delta(\eta_-)&=(\eta_-\otimes\eta_3-\eta_3\otimes\eta_-)/2
\end{aligned}
\qquad
\delta(\eta_3)=\eta_+\otimes\eta_--\eta_-\otimes\eta_+
\end{equation}
and hence
\begin{equation}\label{648}
\delta=-\frac{1}{2}c^\top
\end{equation}
in accordance with Lemma~\ref{lem:C7} ({\it ii\/}). Furthermore, we have
\begin{equation}\label{649}
\widehat{\cal{R}}=\ker(\e)^2
\end{equation}
as in the classical case.

The formalism of connections, based on this calculus $\Gamma$, becomes
essentially the same as in the classical $SU(2)$ case. In particular,
because of the symmetricity of $S_{\inv}^{\wedge2}$, every connection is
multiplicative. Hence, the right-hand side of  Bianchi identity
vanishes. Further, the ``perturbation'' $q_\omega$ also vanishes, as
follows directly from \eqref{647}--\eqref{648} and \eqref{513}. The
presence of the decomposition $\omega=\omega_{\!cl}+\omega_\perp$ is the
only nonclassical phenomena appearing at the level of connections.
\subsection{Trivial Bundles and Non-Admissible Structures}

According to the previous example, compatibility conditions between a
left-covariant differential calculus $\Gamma$ over $G=S_\mu U(2)$, and
``transition functions'' of an appropriate principal bundle can be
fulfilled only in the infinite-dimensional case. This automatically
rules out various interesting finite-dimensional differential
structures.

Such obstructions can be avoided if we restrict the formalism on {\it
trivial principal bundles}. In this case $\cal{B}=S(M)\otimes\cal{A}$,
and a differential calculus on $P$ can be constructed by taking the
product $\Omega(M)\grten\Gamma^\wedge=\Omega(P)$.

Of course, such a calculus over $P$ does not satisfy property
{\it diff3}. On the other hand, if $\Gamma$ is an arbitrary
bicovariant *-calculus then essentially all considerations of Sections 4
and 5 can be repeated in this ``trivial'' framework. The only exception
is that there exist no analogs for classical connections. Because it is
not longer possible to construct the restriction map
$\nu\colon\Gamma_{\inv}\rightarrow \L^*$.

Each connection $\omega$ possesses a global gauge potential $A^\omega\colon
\Gamma_{\inv}\rightarrow\Omega(M)$, given by
\begin{equation}\label{650}
\omega(\vartheta)=(A^\omega\otimes\id)\adj(\vartheta)+1_M\otimes\vartheta.
\end{equation}
The curvature is of the form
\begin{equation}\label{651}
R_\omega=(F^\omega\otimes\id)\adj\qquad
F^\omega=dA^\omega-\langle A^\omega,A^\omega\rangle.
\end{equation}

As a concrete illustration, let us consider the case $G=S_\mu U(2)$
where
$\mu\in(-1,1)\setminus\{0\}$, and let $\Gamma$ be a $4$-dimensional
calculus described in \cite{W3}. By definition, the corresponding right
$\cal{A}$-ideal $\cal{R}$ is generated by multiplets
\begin{gather*}
\mbox{\it 1}=\Bigl\{\,a\bigl(\mu^2\alpha+\alpha^*-(1+\mu^2)1\bigr)\,\Bigr\}
\qquad\mbox{\it
3}=\Bigl\{\,a\gamma,\,a(\alpha-\alpha^*),\,a\gamma^*\,\Bigr\}\\
\mbox{\it 5}=\Bigl\{\,
\gamma^2,\,\gamma(\alpha-\alpha^*),\,\mu^2\alpha^{*2}-
(1+\mu^2)(\alpha\alpha^*-\gamma\gamma^*)+\alpha^2,\,
\gamma^*(\alpha-\alpha^*),\,\gamma^{*2}\,\Bigr\}
\end{gather*}
where $a=\mu^2\alpha+\alpha^*-(\mu^3+1/\mu)1$. It turns out that the
elements
\begin{equation}\label{basis}
\tau=\pi(\mu^2\alpha+\alpha^*)\quad\eta_+=\pi(\gamma)\quad\eta_3=\pi(\alpha-
\alpha^*)\quad\eta_-=\pi(\gamma^*)
\end{equation}
form a basis in $\Gamma_{\inv}$. The canonical right $\cal{A}$-module
structure on $\Gamma_{\inv}$ is given by
\begin{gather}
\begin{aligned}
\tau\circ\gamma^{\phantom{*}}&=\frac{(1-\mu)(1-\mu^3)}{\mu}\eta_+\\
\tau\circ\gamma^*&=\frac{(1-\mu)(1-\mu^3)}{\mu}\eta_-
\end{aligned}\qquad
\begin{aligned}
\tau\circ\alpha^*&=\frac{1+\mu^4}{\mu(1+\mu^2)}\tau-
\frac{\mu(1-\mu)(1-\mu^3)}{1+\mu^2}\eta_3\\
\tau\circ\alpha^{\phantom{*}}&=\frac{1+\mu^4}{\mu(1+\mu^2)}\tau+
\frac{(1-\mu)(1-\mu^3)}{\mu(1+\mu^2)}\eta_3
\end{aligned}\notag\\
\eta_+\circ\gamma^*=\eta_-\circ\gamma=-\frac{(1+\mu)(1-\mu^2)}{\mu(1+
\mu^2)(1-\mu^3)}\tau-\frac{1-\mu^2}{\mu(1+\mu^2)}\eta_3\label{652}\\
\eta_3\circ\gamma=-\frac{1-\mu^2}{\mu}\eta_+\quad
\eta_+\circ\gamma=\eta_-\circ\gamma^*=0\quad
\eta_3\circ\gamma^*=-\frac{1-\mu^2}{\mu}\eta_-\notag\\
\begin{aligned}
-\eta_3\circ\alpha^*&=\frac{(1+\mu)(1-\mu^2)}{\mu(1+\mu^2)
(1-\mu^3)}\tau-\frac{2\mu}{1+\mu^2}\eta_3\\
\eta_3\circ\alpha^{\phantom{*}}&=\frac{\mu(1+\mu)(1-\mu^2)}{(1+\mu^2)
(1-\mu^3)}\tau+\frac{2\mu}{1+\mu^2}\eta_3
\end{aligned}\qquad\quad
\begin{aligned}
\eta_+\circ\alpha&=\eta_+=\eta_+\circ\alpha^*\\
\eta_-\circ\alpha&=\eta_-=\eta_-\circ\alpha^*
\end{aligned}\notag
\end{gather}

The ideal $\cal{R}$ is $\ad, *\k$-invariant. This means \cite{W3} that
$\Gamma$ is a bicovariant *-calculus. By the use of
\eqref{B38} and \eqref{B39} it is easy to determine the *-involution and
the adjoint action $\adj$. We have
\begin{equation}\label{653}
\begin{gathered}
\eta_+^*=\mu\eta_-\qquad\eta_3^*=-\eta_3\qquad\mu\eta_-^*=\eta_+\\
\tau^*=-\tau\qquad\adj(\tau)=\tau\otimes 1\\
\adj(\eta_+)=\eta_+\otimes\alpha^2-\eta_3\otimes\alpha\gamma+
\mu^2\eta_-\otimes
\gamma^2\\
\adj(\eta_3)=(1+\mu^2)\eta_+\otimes\gamma^*\alpha+\eta_3\otimes(
\alpha\alpha^*-\gamma\gamma^*)-(1+\mu^2)\eta_-\otimes\gamma\alpha^*\\
\adj(\eta_-)=\eta_+\otimes\gamma^{*2}+\eta_3\otimes\alpha^*\gamma^*+
\eta_-\otimes
\alpha^{*2}.
\end{gathered}
\end{equation}
We see that $\tau$ form a singlet, while $\{\eta_+,\eta_3,\eta_-\}$
form a triplet, relative to $\adj$.

We are now going to compute the space
$S_{\inv}^{\wedge2}\subseteq\Gamma_{\inv}^{\otimes2}$. Acting by
$(\pi\otimes\pi)\phi$ on the generating elements of $\cal{R}$, using
\eqref{652} and \eqref{B40}, and taking linear combinations we obtain a
lineal spanned by
\begin{gather}
\mbox{\it 5}=\left\{
\begin{gathered}
\eta_+\otimes\eta_+\qquad\mu\eta_3\otimes\eta_3+\mu^4\eta_+\otimes\eta_-
+\eta_-\otimes\eta_+\qquad\eta_-\otimes\eta_-\\
\mu^2\eta_+\otimes\eta_3+\eta_3\otimes\eta_+\qquad\eta_-\otimes\eta_3+\mu^2
\eta_3\otimes\eta_-
\end{gathered}\right\}\notag\\
\mbox{\it 3}=\biggl\{
\frac{1+\mu^4}{1-\mu^3}(\tau\otimes\eta_j+
\eta_j\otimes\tau)+(1-\mu)
\vrkpa_j\vert j\in\{+,-,3\}\biggr\}\label{655}\\
\mbox{\it 1}=\biggl\{
\frac{\phantom{\mu}(1+\mu)(1+\mu^3)}{\mu(1-\mu)(1-\mu^3)}\tau\otimes\tau
+\mu\eta_3\otimes\eta_3-(1+\mu^2)(\eta_+\otimes\eta_-+
\mu^2\eta_-\otimes\eta_+)
\biggr\}\notag
\end{gather}
where we have used the following abbreviations
\begin{gather*}
\vrkpa_+=\eta_+\otimes\eta_3-\mu^2\eta_3\otimes\eta_+\qquad
\vrkpa_-=\eta_3\otimes\eta_--\mu^2\eta_-\otimes\eta_3\\
\vrkpa_3=(1-\mu^2)\eta_3\otimes\eta_3+\mu(1+\mu^2)(\eta_+
\otimes\eta_--\eta_-\otimes\eta_+).
\end{gather*}
\begin{lem}
It turns out that $S_{\inv}^{\wedge2}$ coincides with the lineal generated by
the above elements.
\end{lem}
\begin{pf}
According to Lemma~\ref{B15} ({\it ii\/}) elements
of $S_{\inv}^{\wedge2}$ are $\sigma$-invariant, where $\sigma$
is the canonical
flip-over operator. On the other hand, the space $\ker(I-\sigma)$ is
10-dimensional, spanned by the above elements and $\tau\otimes\tau$.
Consequently, in order to determine $S_{\inv}^{\wedge2}$, it is sufficient to
analyze elements of the form $(\pi\otimes\pi)\phi(a)$, where
$a\in\cal{R}$ is $\ad$-invariant. This follows from the fact that
$(\pi\otimes\pi)\phi$ intertwines $\ad$ and $\adj^{\otimes 2}$. However,
$\ad$-invariant elements of $\cal{R}$ are just linear combinations of
terms of the form
$$r_n=\bigl(\mu^2\alpha+\alpha^*-(\mu^3+1/\mu)1\bigr)\bigl(
\mu^2\alpha+\alpha^*-(1+\mu^2)1\bigr)(\mu^2\alpha+\alpha^*)^n.$$
Inductively using \eqref{B40} and \eqref{652} we find
$$(\pi\otimes\pi)\phi(r_n)=\mu^{-2n}(1+\mu^6)^n(\pi\otimes\pi)\phi(r_0).$$

On the other hand, the last (singlet) term in \eqref{655} coincides with
the element
$\bigl(\mu(1+\mu^2)/(1-\mu)(1-\mu^5)\bigr)(\pi\otimes\pi)\phi(r_0)$.
Hence, elements \eqref{655} generate $S_{\inv}^{\wedge 2}$.
\end{pf}

Let us compute the differential
$d\colon\Gamma^\wedge\rightarrow\Gamma^\wedge$. As first, let us observe
that
\begin{equation}\label{656}
\pi(a)=\frac{\mu}{(1-\mu)(1-\mu^3)}(\tau\circ a-\e(a)\tau),
\end{equation}
for each $a\in\cal{A}$. Indeed, it is evident that \eqref{656} holds for
$a=1$, and from \eqref{652} we conclude that it holds for
$a\in\{\alpha,\alpha^*,\gamma,\gamma^*\}$. Remembering that
$\{\alpha,\alpha^*,\gamma,\gamma^*\}$ generate $\cal{A}$ and using
\eqref{B40} and linearity of both sides of \eqref{656} we conclude that
the above equality holds for all $a\in\cal{A}$.

As a consequence of identities \eqref{656} and \eqref{B43} we find
\begin{equation}\label{657}
d\vartheta=\frac{\mu}{(1-\mu)(1-\mu^3)}(\tau\vartheta-(-1)^{\partial\vartheta}
\vartheta\tau)
\end{equation}
for each $\vartheta\in\Gamma^\wedge$.

Now we shall compute the braid operator
$\sigma\colon\Gamma_{\inv}^{\otimes 2}\rightarrow\Gamma_{\inv}^{\otimes
2}$. Applying Lemma~\ref{B15} ({\it i\/}), and properties
\eqref{652}--\eqref{653} we obtain the following expressions
\begin{gather*}
\sigma(\eta_+\otimes\eta_+)
=\eta_+\otimes\eta_+\qquad\sigma(\eta_-\otimes\eta_-)=\eta_-\otimes\eta_-
\qquad\sigma(\vartheta\otimes\tau)=\tau\otimes\vartheta\\
\begin{aligned}
\sigma(\tau\otimes\eta_-)&=\frac{1+\mu^6}{\mu^2(1+\mu^2)}\eta_-\otimes\tau
+\frac{(1-\mu)(1-\mu^3)}{\mu^2}\vrkpa_-\\
\sigma(\tau\otimes\eta_3)&=\frac{1+\mu^6}{\mu^2(1+\mu^2)}\eta_3\otimes\tau
+\frac{(1-\mu)(1-\mu^3)}{\mu^2}\vrkpa_3\\
\sigma(\tau\otimes\eta_+)&=\frac{1+\mu^6}{\mu^2(1+\mu^2)}\eta_+\otimes\tau
+\frac{(1-\mu)(1-\mu^3)}{\mu^2}\vrkpa_+
\end{aligned}\\
\begin{aligned}
\sigma(\eta_3\otimes\eta_3)&=(3-\mu^2-\frac{1}{\mu^2})\eta_3\otimes\eta_3\\
&\phantom{=}+\frac{1-\mu^4}{\mu}
(\eta_-\otimes\eta_+-\eta_+\otimes\eta_-)-\frac{(1-\mu)(1-\mu^4)}
{\mu^2(1-\mu^3)}\eta_3\otimes\tau
\end{aligned}\\
\begin{aligned}
\sigma(\eta_+\otimes\eta_3)&=\eta_3\otimes\eta_+-\frac{(1+\mu)
(1-\mu^2)}{\mu^2(
1-\mu^3)}\eta_+\otimes\tau+(1-\frac{1}{\mu^2})\eta_+\otimes\eta_3\\
\sigma(\eta_-\otimes\eta_3)&=\eta_3\otimes\eta_-+\frac{(1+\mu)(1-\mu^2)}{
1-\mu^3}\eta_-\otimes\tau+(1-\mu^2)\eta_-\otimes\eta_3\\
\sigma(\eta_3\otimes\eta_+)&=\eta_+\otimes\eta_3+\frac{(1+\mu)(1-\mu^2)}{
1-\mu^3}\eta_+\otimes\tau+(1-\mu^2)\eta_3\otimes\eta_+\\
\sigma(\eta_3\otimes\eta_-)&=\eta_-\otimes\eta_3-\frac{(1+\mu)
(1-\mu^2)}{\mu^2(
1-\mu^3)}\eta_-\otimes\tau+(1-\frac{1}{\mu^2})\eta_3\otimes\eta_-
\end{aligned}\\
\begin{aligned}
\sigma(\eta_+\otimes\eta_-)&=\eta_-\otimes\eta_+-\frac{1-\mu^2}{
\mu(1+\mu^2)}\eta_3\otimes\eta_3-\frac{(1+\mu)(1-\mu^2)}
{\mu(1+\mu^2)(1-\mu^3)}\eta_3\otimes\tau\\
\sigma(\eta_-\otimes\eta_+)&=\eta_+\otimes\eta_-+\frac{1-\mu^2}{
\mu(1+\mu^2)}\eta_3\otimes\eta_3+\frac{(1+\mu)(1-\mu^2)}
{\mu(1+\mu^2)(1-\mu^3)}\eta_3\otimes\tau
\end{aligned}
\end{gather*}

Furthermore, $\mbox{sp}(\sigma)=\bigl\{1,-\mu^2,
-1/\mu^2\bigr\}$. The operator $\sigma$ is diagonalized in the basis
consisting of vectors \eqref{655}, $\tau\otimes\tau$, and the following
two $\adj^{\otimes 2}$-triplets
\begin{equation*}
\begin{gathered}
\tau\otimes\eta_+-\mu^2\eta_+\otimes\tau+\frac{1-\mu^3}{1+\mu^{\phantom{3}}}
\vrkpa_+\\
\tau\otimes\eta_3-\mu^2\eta_3\otimes\tau+\frac{1-\mu^3}{1+\mu^{\phantom{3}}}
\vrkpa_3\\
\tau\otimes\eta_--\mu^2\eta_-\otimes\tau+\frac{1-\mu^3}{1+\mu^{\phantom{3}}}
\vrkpa_-
\end{gathered}\qquad
\begin{gathered}
\mu^2\tau\otimes\eta_+-\eta_+\otimes\tau-\frac{1-\mu^3}{1+\mu^{\phantom{3}}}
\vrkpa_+\\
\mu^2\tau\otimes\eta_3-\eta_3\otimes\tau-\frac{1-\mu^3}{1+\mu^{\phantom{3}}}
\vrkpa_3\\
\mu^2\tau\otimes\eta_--\eta_-\otimes\tau-\frac{1-\mu^3}{1+\mu^{\phantom{3}}}
\vrkpa_-
\end{gathered}
\end{equation*}
corresponding to values $-\mu^2$ and $-1/\mu^2$ respectively.

It is interesting to observe that there exists an indefinite
$\adj$-invariant scalar product on $\Gamma_{\inv}$,  such that $\sigma$
is unitary, relative to the induced product in $\Gamma_{\inv}^{\otimes
2}$. Such a product is given by
\begin{equation}\label{660}
\begin{gathered}
(\tau,\tau)=-\frac{(1-\mu^3)^2(1+\mu^2)}{(1+\mu)^2}\\
(\eta_+,\eta_+)=\mu^2\qquad(\eta_3,\eta_3)=1+\mu^2\qquad(\eta_-,\eta_-)=
1/\mu^2,
\end{gathered}
\end{equation}
while $\eta_\pm,\eta_3,\tau$ are assumed to be
mutually orthogonal. The unitarity  of
$\sigma$ easily follows from the $\adj$-invariance of the introduced
scalar product, and the identity
\begin{equation}\label{661}
(\vartheta\circ a,\eta)=(\vartheta,\eta\circ\k^2(a)^*).
\end{equation}
The product is uniquely determined by the above conditions, up to a
scalar multiple.
There exists a natural splitting $\ker(\e)=\cal{R}\oplus\cal{L}$, where
$\cal{L}$ is the lineal spanned by elements
$\{\gamma,\gamma^*,\alpha-1,\alpha^*-1\}$. This splitting enables us to
introduce the embedded differential map $\delta$. A
direct calculation gives
\begin{equation*}
\begin{aligned}
-(1+\mu^2)\delta(\tau_{\phantom{\zeta}})=&
\tau\otimes\tau+\mu^2\eta_3\otimes\eta_3-\mu(1+\mu^2)(\eta_+\otimes
\eta_-+\mu^2\eta_-\otimes\eta_+)\\
-(1+\mu^2)\delta(\eta_\zeta)=
&\tau\otimes\eta_\zeta+\eta_\zeta\otimes\tau+\vrkpa_\zeta,\quad\zeta
\in\{{+},3,{-}\}.
\end{aligned}
\end{equation*}
The map $\delta$ is coassociative, by construction.

According to \eqref{650}--\eqref{651} the curvature has the form
\begin{gather*}
F^\omega(\tau)=dA^\omega(\tau)+\mu(1-\mu^2)A^\omega(\eta_-)A^\omega(\eta_+)\\
F^\omega(\eta_3)=dA^\omega(\eta_3)+2\mu A^\omega(\eta_+)A^\omega(\eta_-)\\
F^\omega(\eta_-)=dA^\omega(\eta_-)+A^\omega(\eta_3)A^\omega(\eta_-)\qquad
F^\omega(\eta_+)=dA^\omega(\eta_+)+A^\omega(\eta_+)A^\omega(\eta_3).
\end{gather*}

It is worth noticing that essentially the same expressions for {\it
singlet} and {\it triplet} components of $\delta$ and $F^\omega$ can be
obtained in the framework of the previous example.
\subsection{Gauge Theories}
Classical principal bundles provide a natural mathematical framework for
the study of gauge theories \cite{DV}. It is therefore interesting to
see what will be the counterparts of these theories, in the context of
quantum principal bundles \cite{D}.

In analogy with the classical case, the simplest possibility is to
consider lagrangians of the form
\begin{equation}\label{665}
L(\omega)=\sum_\vartheta\bigl(F^\omega(\vartheta),F^\omega(\vartheta)
\bigr)_{\!M}
\end{equation}
where elements $\vartheta$ form an orthonormal system in $\Gamma_{\inv}$
with respect to an $\ad$-invariant scalar product, and $(\,)_M$ is the
scalar product in $\Omega(M)$, induced by a metric on $M$ (here $M$
plays the role of space-time).

Properties of such ``quantum gauge'' theories essentially depend,
besides on the ``symmetry group'' $G$, on the following two
prespecifications.

As first, it is necessary to fix a bicovariant *-calculus $\Gamma$. This
determines kinematical degrees of freedom, as well as ``infinitezimal
gauge transformations''.

Secondly, we have to choose a map $\delta$. This influences dynamical
properties of the theory, because $\delta$ implicitly figures in the
self-interacting part  of \eqref{665}. In the classical case
the curvature is $\delta$-independent.

For instance, in the context of the previous example, we find a
four-component gauge field consisting of mutually interacting singlet
and triplet fields. However if we change $\delta$ and define
$$\delta(\vartheta)=\frac{\mu}{(1-\mu)(1-\mu^3)}(\tau\otimes\vartheta
+\vartheta\otimes\tau)$$
then \eqref{665} will describe non-interacting fields. On the other
hand, in the context of the second example, we find a self-interacting
infinite-component gauge field with all integer spin multiplets in the
game.

Closely related with this line of thinking is a question of ``gauge
transformations''. The most direct way of introducing gauge
trannsformations as vertical automorphisms of $P$ gives nothing new.
Every such automorphism of $P$ preserves the classical part $P_{\!cl}$, and
moreover it is completely determined by the corresponding
``restriction'', which is a classical gauge transformation of $P_{\!cl}$.
In such a way we obtain an isomorphism between gauge groups for $P$ and
$P_{\!cl}$. However, a proper quantum generalization of gauge
transformations can be introduced via the concepts of quantum
(infinitezimal) gauge bundles. These are the bundles associated to $P$,
relative to the adjoint actions $\{\ad,\adj\}$ respectively. It turns
out that operators $h_\omega$, $D_\omega$ and $R_\omega$ are covariant
with respect to natural actions of these bundles on $P$. Moreover, the
lagrangian \eqref{665} is gauge-invariant, in the appropriate sense.
\appendix
\section{Classical Points}
Let $G$ be a compact matrix quantum group. We have denoted by $G_{\!cl}$
the set of *-characters of $\cal{A}$. The elements of $G_{\!cl}$ are
interpretable as classical points of $G$.

The quantum group structure on $G$ induces a classical group structure
on $G_{\!cl}$, in a natural manner. The product and the inverse are given
by
\begin{gather}
gf=(g\otimes f)\phi\label{A1}\\
g^{-1}=g\k.\label{A2}
\end{gather}
The counit $\e\colon\cal{A}\rightarrow\Bbb{C}$ is the neutral element of
$G_{\!cl}$.
\begin{lem}\label{lem:A1}
\bla{i} The formula $$\iota_u(g)_{ij}=g(u_{ij})$$
defines a monomorphism $\iota_u\colon G_{\!cl}\rightarrow GL(n)$.

\smallskip
\bla{ii} The image $\iota_u(G_{\!cl})$ is compact.
\end{lem}
\begin{pf}
Without a lack of generality we can assume \cite{W2} that $u$ is a
unitary matrix. In this case matrices $\iota_u(g)$ belong to $U(n)$.
We have
\begin{equation*}
\begin{split}
\iota_u(gf)_{ij}&=(gf)(u_{ij})=(g\otimes f)\phi(u_{ij})\\
&=\sum_{k=1}^ng(u_{ik})f(u_{kj})=\sum_{k=1}^n\iota_u(g)_{ik}\iota_u(f)_{
kj}.\end{split}
\end{equation*}
Hence $\iota_u$ is a group homomorphism. This map is injective, because
$\cal{A}$ is generated, as a *-algebra, by the matrix elements $u_{ij}$.

Because of the compactness of $U(n)$, it is sufficient to prove that the
image of $\iota_u$ is closed. Let us suppose that a sequence of matrices
$\iota_u(g_k)$ converges to $T\in U(n)$. This means that the sequence of
numbers $g_k(u_{ij})$ converges to $T_{ij}$ for each
$i,j\in\{1,\dots,n\}$. It follows that a sequence $g_k(a)$ is
convergent, for each $a\in\cal{A}$. Now the formula
\begin{equation}
g(a)=\lim_k g_k(a)\label{A3}
\end{equation}
consistently defines a *-character $g\colon\cal{A}\rightarrow\Bbb{C}$ with
the property $\iota_u(g)=T$.
\end{pf}

The monomorphism $\iota_u$ enables us to interpret $G_{\!cl}$ as a compact
group of matrices. In particular, $G_{\!cl}$ is a Lie group in a natural
manner. Furthermore the space $G_{\!cl}$ is an algebraic submanifold of
$U(n)$. The Hopf *-algebra $\cal{A}_{\!cl}$ of polynomial
functions on $G_{\!cl}$ is generated by elements $u_{ij}^{cl}(g)=g(u_{ij})$.
Let ${\restr}_{cl}\colon\cal{A}\rightarrow\cal{A}_{\!cl}$ be the
restriction homomorphism. Let $\lie(G_{\!cl})$ be the (complex) Lie algebra of
$G_{\!cl}$, understood as
the tangent space to $G_{\!cl}$, in the point $\e$.

The formula
\begin{equation}\label{A4}
X(a)=d\bigl({\restr}_{\!cl}(a)\bigr)_{\!\e}(X)
\end{equation}
enables us to interpret elements $X\in\lie(G_{\!cl})$ as certain linear
functionals on $\cal{A}$.
\begin{lem}\label{lem:A2}
\bla{i} We have
\begin{equation}\label{A5}
X(ab)=\e(a)X(b)+\e(b)X(a)
\end{equation}
for each $a,b\in\cal{A}$. Conversely, if
$X\colon\cal{A}\rightarrow\Bbb{C}$ is a hermitian linear functional such
that \eqref{A5} holds then $X$ is interpretable via
\eqref{A4} as a real element of $\lie(G_{\!cl})$.

\smallskip
\bla{ii} In terms of the above identification, the
Lie brackets are given by
\begin{equation}\label{A6}
[X,Y](a)=X(a^{(1)})Y(a^{(2)})-Y(a^{(1)})X(a^{(2)}).
\end{equation}
\end{lem}
\begin{pf}
It is clear that functionals $X$ given by \eqref{A4} satisfy \eqref{A5}.
If $X$ is a hermitian functional satisfying \eqref{A5} then the formula
\begin{equation}\label{A7}
g^t(a)=\e\Bigl[\sum_{k=0}^\infty\frac{1}{k!}
\bigl((\id\otimes X)\phi\bigr)^k(a)t^k\Bigr]
\end{equation}
determines a 1-parameter subgroup of $G_{\!cl}$. The corresponding
generator coincides with $X$, in the sense of \eqref{A4}. Finally,
\eqref{A6} directly follows from \eqref{A4}, and the definition of Lie
brackets.
\end{pf}
In terms of the identification \eqref{A4} the conjugation in
$\lie(G_{\!cl})$ is given by
$$X^*(a)=X(a^*)^*.$$

Let $F\in M_n(\Bbb{C})$ be the canonical intertwiner \cite{W2} between
$u$ and its second contragradient $u^{cc}$. Then
\begin{lem}\label{lem:A3} We have
$$\iota_u(g)F=F\iota_u(g),$$
for each $g\in G_{\!cl}$.
\end{lem}
\begin{pf}
According to definitions of $F$ and $u^{cc}$, we have
$$FuF^{-1}=u^{cc}=(\id\otimes\k^2)u.$$
Acting by $g\in G_{\!cl}$ on this equality, and remembering that
$g\k^2=g$,
we conclude that $F$ and $\iota_u(g)$ commute.
\end{pf}

In a generic case when all eigenvalues of $F$ are mutually different,
the group $G_{\!cl}$ will be very small, because every element
$U\in\iota_u(G_{\!cl})$ is a function of $F$. In particular $G_{\!cl}$ will
be Abelian.

Furthermore, a rough information about the minimal size of $G_{\!cl}$ is
cointained in $F$. According to the results of \cite{W2} we have
$F^{it}\in\iota_u(G_{\!cl})$, for each $t\in\Re$. Hence, the closure of
this 1-parameter subgroup is contained in $\iota_u(G_{\!cl})$. This
closure is isomorphic to a torus the dimension of which is equal to the
number of rationally linearly independent elements of the spectrum of
$\log(F)$.

In the rest of this appendix classical parts of some concrete quantum
groups will be computed.

\medskip
{\it The Classical Case}
\medskip

Let us assume that $\cal{A}$ is commutative. Then so is $A$ and
according to \cite{W2}, $G$ is an ordinary compact matrix group
consisting of characters of $A$. Since every compact matrix group is an
algebraic manifold in the corresponding matrix space, the   restriction
map $g\mapsto g{\restr}\cal{A}$ is an isomorphism between $G$ and
$G_{\!cl}$.

\medskip
{\it Quantum $SU(2)$ groups}
\medskip

By definition \cite{W1}, the $C^*$-algebra representing
continuous functions on the group
$G=S_\mu U(2)$ is generated by elements $\alpha$
and $\gamma$, and relations
\begin{equation}\label{A9}
\begin{gathered}
\alpha\alpha^*+\mu^2\gamma\gamma^*=1\qquad\alpha^*\alpha+\gamma^*\gamma=1\\
\alpha\gamma=\mu\gamma\alpha\qquad\alpha\gamma^*=\mu\gamma^*\alpha
\qquad\gamma\gamma^*=\gamma^*\gamma
\end{gathered}
\end{equation}
while
\begin{equation*}
u=\begin{pmatrix}
\alpha&-\mu\gamma^*\\
\gamma&\phantom{-\mu}\alpha^*\end{pmatrix}.
\end{equation*}
Let us consider the case $\mu\in(-1,1)\setminus\{0\}$. Relations
\eqref{A9} imply that every $g\in G_{\!cl}$ satisfies
$$ \vert g(\alpha)\vert=1\quad g(\gamma)=g(\gamma^*)=0.$$
Consequently $g$ is completely determined by the number $g(\alpha)\in
U(1)$. Moreover, the correspondence $G_{\!cl}\ni g\mapsto g(\alpha)\in
U(1)$ is a group isomorphism.

If $\mu=-1$ relations \eqref{A9} give the following constraints
\begin{equation*}
\begin{split}
\vert g(\alpha)\vert&=1\quad g(\gamma)=g(\gamma^*)=0, \mbox{ or}\\
\vert g(\gamma)\vert&=1\quad g(\alpha)=g(\alpha^*)=0.
\end{split}
\end{equation*}
In this case
$$G_{\!cl}=U(1)\wedge\Bbb{Z}_2$$
in a natural manner.

\medskip
{\it Quantum $SU(n)$ groups}
\medskip

Let us assume that $\mu\in(-1,1)\setminus\{0\}$. By definition \cite{W4}
the $C^*$-algebra $A$ representing continuous functions on
$G=S_\mu U(n)$ groups is generated by elements $u_{ij}$, where
$i,j\in\{1,\dots,n\}$, and relations
\begin{equation}\label{A10}
\begin{gathered}
\sum_{j=1}^n u_{ij}u_{kj}^*=\delta_{ik}I\qquad\sum_{j=1}^n u_{ji}^*
u_{jk}=\delta_{ik}I\\
{\sum}_* u_{i_1j_1}\dots u_{i_nj_n}E_{j_1\dots j_n}=E_{i_1\dots i_n}I.
\end{gathered}
\end{equation}
The last summation is performed over indexes $j$, and
$$E_{i_1\dots i_n}=(-\mu)^{I(i)}$$
where $I(i)$ is the number of inversions in the sequence
$i=(i_1,\dots,i_n)$, if the sequence is a permutation. Other components
of $E$ vanish, by definition.

The fundamental representation of $G$ is irreducible. Let us compute the
canonical intertwiner $F$. The conjugate representation $u^c$ can be
naturally realized as a subrepresentation of the $(n-1)$-th tensor power
of $u$. The carrier space $H$ is spanned by vectors
$$x_k={\sum}_* E_{kj_1\dots j_{n-1}}e_{j_1}\otimes\dots\otimes
e_{j_{n-1}}.$$
Here $e_i$ are absolute basis vectors in $\Bbb{C}^n$, and the summation is
performed over indexes $j$. We have
$$F=cj^\dagger j$$
where $c>0$ and $j\colon\Bbb{C}\rightarrow H$ is the canonical
antilinear map defined by $j(e_k)=x_k$. Now, a direct computation gives
\begin{equation}\label{A11}
Fe_k=\mu^{2k-n-1}e_k
\end{equation}
for each $k\in\{1,\dots,n\}$.

According to Lemma~\ref{lem:A3}, matrices $\iota_u(g)$ are diagonal.
Relations \eqref{A10} imply that corresponding diagonal elements
$\iota_{ii}(g)$ are complex units, and that
$$\prod_i\iota_{ii}(g)=1.$$
The same relations imply that conversely for any sequence of numbers
$z_1,\dots z_n\in U(1)$ satisfying $\prod_i z_i=1$ there exists the
unique $g\in G_{\!cl}$ such that $\iota_{ii}(g)=z_i$. In summary, $G_{\!cl}$
is isomorphic to the $(n-1)$-dimensional torus.

\medskip
{\it Abelian Quantum Groups}
\medskip

If $G$ is Abelian then every subgroup of $G$ is Abelian, too. In
particular $G_{\!cl}$ is an Abelian compact matrix group, and as such it
is isomorphic to a product of a torus with a finite Abelian group.

According to \cite{W2} there exist a discrete finitely generated group
$\Gamma$, Hilbert space $H$ and a unitary representation
$U\colon\Gamma\rightarrow U(H)$ (the square of which is contained in its
multiple) such that $\cal{A}$ is isomorphic to the *-algebra generated
by the image of $U$. Furthermore
$$\phi\bigl(U(\gamma)\bigr)=U(\gamma)\otimes
U(\gamma)\quad\e\bigl(U(\gamma)\bigr)=1\quad\k\bigl(U(\gamma)\bigr)=
U(\gamma)^{-1}$$
for each $\gamma\in\Gamma$.
Since operators $U(\gamma)$ are mutually linearly independent \cite{W2},
every character $g\in G_{\!cl}$ can be viewed as a character on $\Gamma$,
via
$$g(\gamma)=g\bigl(U(\gamma)\bigr),$$
and vice versa. In other words $G_{\!cl}$ is isomorphic to the group of
characters of $\Gamma$.

\medskip
{\it Universal Unitary Quantum Matrix Groups}
\medskip

Let us consider a positive matrix $F\in M_n(\Bbb{C})$ such that
$$\mbox{\rm tr}(F)=\mbox{\rm tr}(F^{-1}).$$
 Let $A_F$ be a $C^*$-algebra
generated by elements $u_{ij}$, where $i,j\in\{1,\dots,n\}$, and
relations
\begin{equation}\label{A12}
\begin{aligned}
\sum_{j=1}^n u_{ij}u_{kj}^*&=\delta_{ik}I\\
\sum_{j=1}^n u_{ij}^*u_{kj}^F&=\delta_{ik}I
\end{aligned}\qquad\quad
\begin{aligned}\sum_{j=1}^n u_{ji}^*u_{jk}&=\delta_{ik}I\\
\sum_{j=1}^n u_{ji}^Fu_{jk}^*&=\delta_{ik}I,
\end{aligned}
\end{equation}
where $u^F=FuF^{-1}$.

The pair $G_F=(A_F,u)$ is a compact matrix quantum group. We are going
to describe the category of unitary representations of $G_F$. Let
$\cal{T}$ be a concrete monoidal $W^*$-category \cite{W4} generated by
elements $u$ and $u^c$, with carrier Hilbert spaces
$H_u=\Bbb{C}^n$ and $H_{u^{\!c}}=H_u^*$. It will be assumed that
$H_u$ is endowed with the standard scalar product, while the product in
$H_u^*$ is specified by $(x,y)=(x,Fy)$. The objects
of $\cal{T}$ are just the words of $u$ and $u^c$ (including the unit
object). By definition, morphisms between objects of $\cal{T}$ are
generated by ``elementary morphisms'' $t\colon\Bbb{C}\rightarrow
H_u\otimes H_u^*$ and $\bar{t}\colon H_u^*\otimes
H_u\rightarrow\Bbb{C}$,
which are given by
$$ t(1)=\sum_{i=1}^ne_i\otimes j(e_i)\qquad \bar{t}(x\otimes y)
=(j^{-1}x,y) $$
where $j\colon H_u\rightarrow H_u^*$ is the complex conjugation map.
By construction $u$ and $u^c$ are mutually conjugate objects.

Then $G_F=(A_F,u)$ is the universal $\cal{T}$-admissible pair ($u$ is a
distinguished object). In other words $G_F$ is a compact matrix quantum
group  corresponding to $\cal{T}$, in the framework of Tannaka-Krein duality
\cite{W4}. The antipode acts as follows
$$\k(u_{ij})=u^*_{ji}\qquad\k(u_{ij}^*)=u^F_{ji}.$$
The map $F=j^\dagger j$  is just the canonical intertwinner between $u$
and $u^{\!cc}$. According to Lemma~\ref{lem:A3} and relations \eqref{A12},
the elements of $\iota_u\bigl(G_F^{cl}\bigr)$ are precisely unitary
matrices commuting with $F$. Hence, $$G_F^{cl}=U(n_1)\times
\dots\times  U(n_k)$$
where $n_i$ are multiplicities of eigenvalues of $F$.

\section{Universal Differential Envelopes}
Let $\cal{A}$ be a complex unital associative algebra and $\Gamma$ a
first-order calculus \cite{W3} over $\cal{A}$. Let $\Gamma^\otimes$ be
the corresponding ``tensor bundle'' algebra, and let
$S^\wedge$ be the ideal in $\Gamma^\otimes$ generated by elements of the form
\begin{equation}\label{B1}
Q=\sum_i da_i\otimes_{\cal{A}}db_i,\quad\mbox{where}\quad\sum_i a_idb_i=0.
\end{equation}
By definition, $S^\wedge$ is a graded ideal in
$\Gamma^\otimes$ and its first (generally) nontrivial component
coincides with the set of elements $Q$ of the form \eqref{B1}.

Let $\Gamma^\wedge=\Gamma^\otimes/S^\wedge$ be the
corresponding factoralgebra.
\begin{pro}\label{pro:B1}
There exists the unique linear map $d\colon\Gamma^\wedge\rightarrow
\Gamma^\wedge$ extending the derivation $d\colon\cal{A}\rightarrow
\Gamma$ such that
\begin{align*}
d^2&=0\\
d(\vartheta\eta)&=d(\vartheta)\eta+(-1)^{\partial\vartheta}\vartheta
d(\eta)
\end{align*}
for each $\vartheta,\eta\in\Gamma^\wedge$.
\end{pro}
\begin{pf}
The formula
\begin{equation}\label{B2}
d\Bigl(\sum_i a_idb_i\Bigr)=\sum_ida_idb_i
\end{equation}
consistently defines a linear map $d\colon\Gamma\rightarrow
\Gamma^\wedge$. We have
\begin{align}
dd(a)&=0\label{B3}\\
d(a\vartheta)&=(da)\vartheta+ad(\vartheta)\label{B4}\\
d(\vartheta a)&=d(\vartheta)a-\vartheta (da)\label{B5}
\end{align}
for each $a\in\cal{A}$ and $\vartheta\in\Gamma$.
Equalities \eqref{B4}--\eqref{B5} imply that maps $d$
admit the unique
extension $d\colon\Gamma^\otimes\rightarrow\Gamma^\wedge$ satisfying
\begin{equation}\label{B6}
d(w\otimes_{\cal{A}}u)=d(w)\Pi(u)+(-1)^{\partial
w}\Pi(w)d(u)
\end{equation}
where $\Pi\colon\Gamma^\otimes\rightarrow\Gamma^\wedge$ is the
projection map. Equations \eqref{B3} and \eqref{B6} imply that
$S^\wedge\subseteq\ker(d)$. Consequently, there exists the unique
map $d\colon\Gamma^\wedge\rightarrow\Gamma^\wedge$ defined as a
factorization of the previous $d$ through $\Pi$.
This map possesses all desired properties.\end{pf}

The differential algebra $\Gamma^\wedge$ possesses the following {\it
universality property}.
\begin{pro}\label{pro:B2}
Let $\Omega$ be a differential algebra with a differential
$d_\Omega\colon\Omega\rightarrow\Omega$.

\smallskip
\bla{i} Let $\varphi\colon\cal{A}\rightarrow\Omega$ be a homomorphism
admitting the extension $\sharp_\varphi\colon\Gamma\rightarrow\Gamma$,
given by
$$\sharp_\varphi\bigl(ad(b)\bigr)=\varphi(a)d_\Omega\varphi(b).$$
Then there exists the unique differential algebra homomorphism
$\varphi^\wedge\colon\Gamma^\wedge\rightarrow\Omega$
extending both $\varphi$ and $\sharp_\varphi$.

\smallskip
\bla{ii} Similarly, if $\varphi\colon\cal{A}\rightarrow\Omega$  is an
antimultiplicative linear map and if there exists $\sharp_\varphi\colon
\Gamma\rightarrow\Omega$ satisfying
$$\sharp_\varphi\bigl(ad(b)\bigr)=d_\Omega\varphi(b)\varphi(a)$$
then $\varphi$ and $\sharp_\varphi$ admit the unique extension
$\varphi^\wedge\colon\Gamma^\wedge\rightarrow\Omega$ satisfying
\begin{align*}
\varphi^\wedge d&=d_\Omega\varphi^\wedge\\
\varphi^\wedge(\vartheta\eta)&=(-1)^{\partial\vartheta\partial\eta}
\varphi^\wedge(\eta)\varphi^\wedge(\vartheta)
\end{align*}
for each $\vartheta,\eta\in\Gamma^\wedge$.
\end{pro}
\begin{pf}
We shall check the statement ({\it i\/}). The maps $\varphi$ and
$\sharp_\varphi$ admit the unique common multiplicative extension
$\varphi^\otimes\colon\Gamma^\otimes\rightarrow\Omega$. It is easy to
see that $\varphi^\otimes(Q)=0$, for each $Q$ given by \eqref{B1}. In
other words, $S^\wedge\subseteq\ker\bigl(\varphi^\otimes\bigr)$ and hence
$\varphi^\wedge$ can be factorized through $\Pi$. In such a way we
obtain the desired map $\varphi^\wedge$. The uniqueness follows from the
fact that $\Gamma^\wedge$ is generated by $\cal{A}$, as a differential
algebra.
\end{pf}

A similar statement can be formulated for antilinear maps $\varphi$. As
a simple corollary we obtain
\begin{pro}\label{pro:B3} Let us assume that $\cal{A}$ is a *-algebra
and that $\Gamma$ is a *-calculus. There exists the unique
antilinear involution $*\colon\Gamma^\wedge\rightarrow\Gamma^\wedge$
extending *-involutions on $\cal{A}$ and $\Gamma$ and satisfying
\begin{equation}\label{B14}
\begin{aligned}
d(\vartheta^*)&=d(\vartheta)^*\\
(\vartheta\eta)^*&=(-1)^{\partial\eta\partial\vartheta}\eta^*\vartheta^*
\end{aligned}
\end{equation}
for each $\vartheta,\eta\in\Gamma^\wedge.\qed$
\end{pro}
Let us consider some examples of universal envelopes, interesting from
the point of view of quantum principal bundles.
\begin{pro}\label{pro:B16} \bla{i} Let $M$ be a compact manifold. Then
$$\Omega(M)=[\Omega^1(M)]^\wedge.$$

\smallskip
\bla{ii} If $P$ is a quantum principal
bundle over $M$ and $\Gamma$ an arbitrary
admissible calculus over $G$ then
$$\Omega(P,\Gamma)=[\Omega^1(P,\Gamma)]^\wedge.$$
In other words $\Omega(M)$ and $\Omega(P,\Gamma)$ are understandable as
universal envelopes.
\end{pro}
\begin{pf} We shall prove the statement \bla{i}. The proof of \bla{ii}
is based on \bla{i} and the universality of $\Gamma^\wedge$.

The space $\Omega^1(M)\otimes_M\Omega^1(M)$ is naturally isomorphic to a
$S(M)$-module of covariant $2$-tensors. To prove \bla{i} it is
sufficient to check that $S^{\wedge 2}$ coincides with the space
$\Sigma$ of symmetric $2$-tensors. According to universality of
$\Omega^1(M)^\wedge$ we have $S^{\wedge 2}\subseteq \Sigma$. Conversely,
elements of the form $q=df\otimes_M df$, where $f\in S(M)$, generate the
module $\Sigma$. Every such element belongs to $S^{\wedge 2}$, because
of the identity $fdf-d(f^2)/2=0$. Hence, $\Sigma\subseteq S^{\wedge 2}$.
\end{pf}

The algebra $\Gamma^\wedge$ can be alternatively constructed by applying
a method of extended bimodules \cite{C,W1,W3}.

Let $\Gamma^{\wedge}\{X\}$ be the graded differential algebra generated by
$\Gamma^\wedge$, a first-order element $X$, and the following
relations
\begin{equation}\label{B15}
\begin{gathered}
X^2=0\qquad d(X)=0\\
X\vartheta-(-1)^{\partial\vartheta}\vartheta X=d(\vartheta).
\end{gathered}
\end{equation}
On the other hand, let $\widetilde{\Gamma}$ be the extended bimodule
$$\widetilde{\Gamma}=\cal{A}\widetilde{X} \oplus\Gamma$$
with a right $\cal{A}$-module structure specified by
\begin{equation}\label{B18}
\widetilde{X}a=a\widetilde{X}+d(a).
\end{equation}
\begin{pro}\label{pro:B4}
There exists the unique homomorphism
$\Pi^\sstar\colon\widetilde{\Gamma}^\otimes\rightarrow\Gamma^{\wedge}\{X\}$
satisfying $\Pi^\sstar
(\widetilde{X})=X$ and extending the factorization map $\Pi$.
The kernel of $\Pi^\sstar$ coincides with the ideal in
$\widetilde{\Gamma}^\otimes$ generated by
$\widetilde{X}\otimes_{\cal{A}}\widetilde{X}.\qed$
\end{pro}

In other words, $\Gamma^\wedge$ can be viewed as a differential
subalgebra of $\widetilde{\Gamma}^\otimes/\ker(\Pi^\sstar)$ generated
by $\cal{A}$.

Let us turn to the quantum group context, and assume that $\cal{A}$
represents polynomial functions on a compact matrix quantum group $G$.
The following statement is a direct corrolary of
Proposition~\ref{pro:B2}.
\begin{pro}\label{pro:B7}
\bla{i} Let $\Gamma$ be a left-covariant calculus over $G$, with the
corresponding left action
$\ell_\Gamma\colon\Gamma\rightarrow\cal{A}\otimes\Gamma$. Then there
exists the unique
homomorphism $\ell_\Gamma^\wedge\colon\Gamma^\wedge\rightarrow
\cal{A}\otimes\Gamma^\wedge$ which extends $\phi$ and such that
\begin{equation}\label{B27}
\ell_\Gamma^\wedge d=(\id\otimes d)\ell_\Gamma^\wedge.
\end{equation}
This map also extends $\ell_\Gamma$ and satisfies
\begin{align}
(\e\otimes\id)\ell_\Gamma^\wedge&=\id\label{B28}\\
(\phi\otimes\id)\ell_\Gamma^\wedge&=(\id\otimes\ell_\Gamma^\wedge)
\ell_\Gamma^\wedge.\label{B29}
\end{align}
If $\Gamma$ is also a *-calculus then $\ell_\Gamma^\wedge$ is
hermitian, in a natural manner.

\smallskip
\bla{ii} Similarly, if $\Gamma$ is right-covariant then there exists the
unique homomorphism  $\wp_\Gamma^\wedge\colon\Gamma^\wedge\rightarrow
\Gamma^\wedge\otimes\cal{A}$ extending $\phi$ and satisfying
\begin{equation}\label{B31}
\wp_\Gamma^\wedge d=(d\otimes\id)\wp_\Gamma^\wedge.
\end{equation}
This homomorphism also extends the right action map $\wp_\Gamma\colon
\Gamma\rightarrow\Gamma\otimes\cal{A}$ and satisfies
\begin{align}
(\id\otimes\e)\wp_\Gamma^\wedge&=\id\label{B32}\\
(\wp_\Gamma^\wedge\otimes\id)\wp_\Gamma^\wedge&=(\id\otimes\phi)
\wp_\Gamma^\wedge.\label{B33}
\end{align}
If, in addition, the calculus
$\Gamma$ is *-covariant then $\wp_\Gamma^\wedge$ preserves corresponding
*-structures.

\smallskip
\bla{iii}
If $\Gamma$ is bicovariant then so is $\Gamma^\wedge$, that is
\begin{equation}\label{B35}
(\id\otimes\wp_\Gamma^\wedge)\ell_\Gamma^\wedge=
(\ell_\Gamma^\wedge\otimes\id)\wp_\Gamma^\wedge.\qed
\end{equation}
\end{pro}

There exists a natural grade-preserving coaction map
$c\colon\Gamma^\wedge\otimes\cal{A}\rightarrow\Gamma^\wedge$, given by
\begin{equation}
c(\vartheta\otimes a)=\k(a^{(1)})\vartheta a^{(2)}.
\end{equation}
The same formula determines the coaction of $G$ on $\Gamma^\otimes$. We
have
\begin{equation}
c(\vartheta\otimes 1)=\vartheta\qquad c\bigl(c(\vartheta\otimes
a)\otimes b\bigr)=c\bigl(\vartheta\otimes(ab)\bigr).
\end{equation}
If $\Gamma$ is *-covariant then
\begin{equation}
c(\vartheta\otimes a)^*=c\bigl(\vartheta^*\otimes\k(a)^*\bigr)
\end{equation}
for each $\vartheta\in\Gamma^{\wedge}$ and $a\in\cal{A}$.

\begin{lem}\label{lem:B8} Let us assume that $\Gamma$ is
right-covariant. Then the following identity holds
\begin{equation}\label{B36}
\rig_\Gamma^\wedge c(\vartheta\otimes a)
=\sum_kc(\vartheta_k\otimes a^{(2)})\otimes
\k(a^{(1)})c_k a^{(3)},
\end{equation}
where $\Sum_k\vartheta_k\otimes c_k=\rig_\Gamma^\wedge(\vartheta)$.
\end{lem}
\begin{pf}
We compute
\begin{multline*}
\rig_\Gamma^\wedge c
(\vartheta\otimes a)=\rig_\Gamma^\wedge\bigl(\k(a^{(1)})\vartheta
a^{(2)}\bigr)=\sum_k\k(a^{(2)})\vartheta_k a^{(3)}\otimes\k(a^{(1)})c_k
a^{(4)}\\=\sum_kc(\vartheta_k\otimes a^{(2)})\otimes\k(a^{(1)})c_k
a^{(3)}.\qed
\end{multline*}
\renewcommand{\qed}{}
\end{pf}
\begin{defn} A first-order calculus $\Gamma$ over $G$ is called
{\it $\k$-covariant} iff there exists a linear map $\sharp_\k\colon
\Gamma\rightarrow\Gamma$ such that
\begin{align}
d\k(a)&=\sharp_\k d(a)\\
\sharp_\k(a\vartheta)&=\sharp_\k(\vartheta)\k(a)
\end{align}
for each $a\in\cal{A}$ and $\vartheta\in\Gamma$.
\end{defn}

The map $\sharp_\k$ is uniquely determined by the above conditions.
Furthermore it is bijective and
\begin{equation}
\sharp_\k(\vartheta a)=\k(a)\sharp_\k(\vartheta).
\end{equation}
According to Proposition~\ref{pro:B2} the map $\sharp_\k$ can be
extended to a $d$-preserving graded-antiautomorphism $\k^\wedge\colon
\Gamma^\wedge\rightarrow\Gamma^\wedge$. If $\Gamma$ is *-covariant then
\begin{equation}\label{B56}
\k^\wedge(\k^\wedge(\vartheta^*)^*)=\vartheta
\end{equation}
for each $\vartheta\in\Gamma^\wedge$.

\begin{pro}\label{pro:B13} If the calculus $\Gamma$ is left-covariant then
$\k$-covariance is equivalent to bicovariance.\qed
\end{pro}

{}From this moment we assume that $\Gamma$ is left-covariant. Let us
denote by $\Gamma_{\inv}^\sstar$ the space of left-invariant elements
of $\Gamma^\sstar$, for $\star\in\{{\otimes},{\wedge}\}$.
The space $\Gamma^\otimes_{\inv}$ is naturally
identificable with the tensor algebra over $\Gamma_{\inv}$.
Proposition~\ref{pro:B7} ({\it i\/}) implies that $\Gamma_{\inv}^\wedge$
is a graded-differential subalgebra of $\Gamma^\wedge$. This algebra is
generated by $\Gamma_{\inv}$.

Let $\ell_\Gamma^\otimes\colon\Gamma^\otimes
\rightarrow\cal{A}\otimes\Gamma^\otimes$ be the left action of $G$ on
$\Gamma^\otimes$. The ideal $S^\wedge$ is $\ell_\Gamma^\otimes$-invariant and
$\ell_\Gamma^\wedge$ coincides with the factorized $\ell_\Gamma^\otimes$
through $\Pi$. The ideal $S^\wedge$ is decomposable as
$$S^\wedge\leftrightarrow \cal{A}\otimes S^\wedge_{\inv}.$$
It is easy to see that $\Pi(\Gamma^\otimes_{\inv})=\Gamma_{\inv}^\wedge$.
In other words $$\Gamma_{\inv}^\otimes/ S_{\inv}^\wedge
=\Gamma^\wedge_{\inv}.$$

The spaces $\Gamma_{\inv}^\sstar$ are $c$-invariant, and hence
the formula
\begin{equation}\label{B21}
\vartheta\circ a=c(\vartheta\otimes a)
\end{equation}
determines a right $\cal{A}$-module structure on them.
The following identities hold
\begin{align}
1\circ a&=\e(a)1\label{B23}\\
(\vartheta\eta)\circ a&=(\vartheta\circ a^{(1)})(\eta\circ
a^{(2)}).\label{B24}
\end{align}
If $\Gamma$ is *-covariant then the spaces $\Gamma^\sstar_{\inv}$
are *-invariant and we can write
\begin{equation}\label{B26}
(\vartheta\circ a)^*=\vartheta^*\circ \k(a)^*.
\end{equation}

Let $\pi\colon\cal{A}\rightarrow\Gamma_{\inv}$ be a linear map given by
\begin{equation}\label{B37}
\pi(a)=\k(a^{(1)})d(a^{(2)}).
\end{equation}
The map $\pi$ is surjective, and $\pi(1)=0$.

\begin{lem}\label{lem:B10} The following identities hold
\begin{gather}
\pi(a)\circ b=\pi\bigl(ab-\e(a)b\bigr)\label{B40}\\
d(a)=a^{(1)}\pi(a^{(2)})\label{B43}\\
d\pi(a)=-\pi(a^{(1)})\pi(a^{(2)})\label{B42}.
\end{gather}
\end{lem}
\begin{pf}
All these equalities follow by straightforward transformations, applying
the definition of $\pi$.
\end{pf}

We can write
$$\Gamma_{\inv}=\ker(\e)/\cal{R}$$
where $\cal{R}=\ker(\e)\cap\ker(\pi)$ is the right $\cal{A}$-ideal
which, in the sense of \cite{W3}, canonically determines the  structure
of $\Gamma$. According to \cite{W3}, the calculus is *-covariant iff
$\k(\cal{R})^*=\cal{R}$. In this case
\begin{equation}\label{B38}
\pi(a)^*=-\pi\bigl(\k(a)^*\bigr)
\end{equation}
for each $a\in\cal{A}$.

\begin{lem}\label{lem:B11} The space $S_{\inv}^{\wedge2}\subseteq
\Gamma_{\inv}\otimes\Gamma_{\inv}$ is consisting precisely of elements of
the form
\begin{equation}\label{B44}
q=\pi(a^{(1)})\otimes\pi(a^{(2)}),
\end{equation}
where $a\in\cal{R}$.
\end{lem}
\begin{pf}
The space $S_{\inv}^{\wedge 2}$ is consisting of left-invariant
projections of elements $Q$ given by \eqref{B1}. In terms of the
identification $\Gamma^\otimes\leftrightarrow\cal{A}\otimes
\Gamma_{\inv}^{\otimes}$ we have
$$Q=\sum_ia_i^{(1)}b_i^{(1)}\otimes\Bigl
\{\bigl(\pi(a_i^{(2)})\circ
b_i^{(2)}\bigr)\otimes\pi(b_i^{(3)})\Bigr\}$$
and hence
$$(\e\otimes\id)(Q)=\sum_i\pi(a_ib_i^{(1)})\otimes\pi(b_i^{(2)})-
\sum_i\e(a_i)\pi(b_i^{(1)})\otimes\pi(b_i^{(2)}).$$
The first summand on the right-hand side of the above equality vanishes,
because of $\Sum_i a_idb_i=\Sum_ia_ib_i^{(1)}\otimes\pi(b_i^{(2)})=0$. On the
other hand, the elements of the form $r=\Sum_i\e(a_i)b_i$ cover the
whole space $\ker(\pi)=\Bbb{C}1\oplus\cal{R}$.
\end{pf}

Actually the space $S_{\inv}^{\wedge 2}$ generates the whole ideal
$S_{\inv}^\wedge$ in
$\Gamma_{\inv}^\otimes$. In other words, $\Gamma_{\inv}^\wedge$ is a
quadratic algebra.

\begin{pro}\label{pro:B12} The following conditions are equivalent\par
\bla{i} The calculus $\Gamma$ is bicovariant.

\smallskip
\bla{ii} The coproduct map $\phi$ is (necessarily uniquely) extendible to
the homomorphism $\widehat{\phi}\colon\Gamma^\wedge
\rightarrow\Gamma^\wedge\grten\Gamma^\wedge$ of differential algebras.
\end{pro}
\begin{pf}
Let us suppose that ({\it i\/}) holds. Let $\widehat{\phi}\colon\Gamma
\rightarrow\Gamma^\wedge\grten\Gamma^\wedge$ be a map
given by
\begin{equation}\label{B47}
\widehat{\phi}(\vartheta)=\ell_\Gamma(\vartheta)\oplus\wp_\Gamma(\vartheta).
\end{equation}
Proposition~\ref{pro:B3} implies that this map, together with $\phi$,
can be further extended
to a differential homomorphism $\widehat{\phi}\colon\Gamma^\wedge
\rightarrow\Gamma^\wedge\grten\Gamma^\wedge$. Conversely, if ({\it
ii\/}) holds then formula \eqref{B47} defines the left and the right
actions of $G$ on $\Gamma$. In other words the calculus is bicovariant.
\end{pf}

Let us assume that $\Gamma$ is bicovariant. This is equivalent
\cite{W3} to
$\ad(\cal{R})\subseteq\cal{R}\otimes\cal{A}$.
The spaces
$\Gamma_{\inv}^\sstar$ are invariant under the right
action of $G$.
Let $\adj^\star\colon\Gamma_{\inv}^\sstar
\rightarrow\Gamma_{\inv}^\sstar\otimes\cal{A}$ be the corresponding
restriction maps. The following identity holds
\begin{equation}
\adj^\star(\vartheta\circ a)=\sum_k\vartheta_k\circ a^{(2)}\otimes
\k(a^{(1)})c_k a^{(3)},
\end{equation}
where $\Sum_k\vartheta_k\otimes c_k=\adj^\star(\vartheta)$.

Explicitly, the map
$\adj\colon\Gamma_{\inv}\rightarrow\Gamma_{\inv}\otimes\cal{A}$ is given
by\begin{equation}\label{B39}\adj\pi=(\pi\otimes\id)\ad.\end{equation}

The map $\widehat{\phi}$ possesses the property
\begin{equation}\label{B48}
(\id\otimes\widehat{\phi})\widehat{\phi}=
(\widehat{\phi}\otimes\id)\widehat{\phi}
\end{equation}
as follows from the coassociativity of $\phi$. Let $\e^\wedge\colon
\Gamma^\wedge\rightarrow\Bbb{C}$ be a homomorphism acting as $\e$ on
$\cal{A}$, and vanishing on higher-order components. Then
\begin{equation}\label{B50}
(\e^\wedge\otimes\id)\widehat{\phi}=(\id\otimes\e^\wedge)\widehat{\phi}=\id.
\end{equation}
If in addition $\Gamma$ admits a *-structure then $\widehat{\phi}$ is
a hermitian map.
Let us denote by $m^\wedge$ the multiplication map in $\Gamma^\wedge$.
\begin{pro}\label{pro:B14} The following identity holds
\begin{equation}\label{B64}
m^\wedge(\k^\wedge\otimes\id)\widehat{\phi}=
m^\wedge(\id\otimes\k^\wedge)\widehat{\phi}=1\e^\wedge.
\end{equation}
\end{pro}
\begin{pf}
It follows from the definition of $\k^\wedge$, $\e^\wedge$ and
$\widehat{\phi}$.
\end{pf}

Let $\sigma\colon\Gamma\otimes_{\cal{A}}\Gamma\rightarrow
\Gamma\otimes_{\cal{A}}\Gamma$ be the canonical braid operator
\cite{W3}. This map intertwines the corresponding left and right
actions. In particular it  is reduced in the space
$\Gamma_{\inv}^{\otimes 2}$. Its left-invariant restriction is
explicitly given by
\begin{lem}\label{lem:B15}
We have
\begin{equation}\label{B65}
\sigma(\eta\otimes\vartheta)=\sum_k\vartheta_k\otimes(\eta\circ c_k)
\end{equation}
for each $\vartheta,\eta\in\Gamma_{\inv}$, where
$\Sum_k\vartheta_k\otimes c_k=\adj(\vartheta)$.
\end{lem}
\begin{pf}
Using the definition \cite{W3} of $\sigma$ and performing direct
transformations we obtain
\begin{multline*}
\sigma(\eta\otimes\vartheta)=\sum_k\sigma(\eta\otimes_{\cal{A}}
\bigl(\vartheta_k
\k(c_k^{(1)})\bigr)c_k^{(2)}=\sum_k\vartheta_k\k(c_k^{(1)})
\otimes_{\cal{A}}\eta c_k^{(2)}\\
=\sum_k\bigl(\vartheta_k\k(c_k^{(1)})c_k^{(2)}\bigr)\otimes_{\cal{A}}
(\eta\circ c_k^{(3)})=\sum_k\vartheta_k\otimes(\eta\circ c_k).\qed
\end{multline*}
\renewcommand{\qed}{}
\end{pf}

Let $\Gamma^\vee$ be the braided exterior algebra \cite{W3}
built over $\Gamma$. In view of the universality of $\Gamma^\wedge$
there exists the unique homomorphism $\Uni\colon\Gamma^\wedge\rightarrow
\Gamma^\vee$ of graded differential algebras reducing to the identity on
$\Gamma$ and $\cal{A}$. In particular
\begin{equation}\label{B66}
S^{\wedge 2}\subseteq\ker(I-\sigma).
\end{equation}
This also follows from \eqref{B40}, \eqref{B65} and Lemma~\ref{lem:B11}.
The map $\Uni$ is surjective, but generally not injective. Moreover, the
algebra $\Gamma^\vee$ is generally not quadratic.
\section{The Minimal Admissible Calculus}
Let $\widehat{\cal{R}}$ be the set of elements $a\in\ker(\e)$ satisfying
\begin{equation}\label{C1}
(X\otimes\id)\ad(a)=0
\end{equation}
for each $X\in\lie(G_{\!cl})$.
\begin{lem}\label{lem:C1}
The space $\widehat{\cal{R}}$ is a right $\cal{A}$-ideal and
\begin{gather}
\ad(\widehat{\cal{R}})\subseteq\widehat{\cal{R}}\otimes\cal{A}\label{Cii}\\
\k(\widehat{\cal{R}})^*=\widehat{\cal{R}}.\label{Ciii}
\end{gather}
\end{lem}
\begin{pf} Let us assume that $a\in\widehat{\cal{R}}$ and
$b\in\ker(\e)$. A direct computation gives
\begin{multline*}
(X\otimes\id)\ad(ab)=X(a^{(2)}b^{(2)})\k(a^{(1)}b^{(1)})a^{(3)}b^{(3)}\\
=X(a^{(2)})\k(b^{(1)})\k(a^{(1)})a^{(3)}b^{(2)}+\e(a)X(b^{(2)})
\k(b^{(1)})b^{(3)}=0.
\end{multline*}
Hence $\widehat{\cal{R}}$ is a right ideal in $\cal{A}$. Properties
\eqref{Cii}--\eqref{Ciii} follow from the definition of
$\widehat{\cal{R}}$, applying elementary properties of maps figuring
in the game.
\end{pf}

Let $\Gamma$ be the left-covariant calculus which canonically, in the
sense of \cite{W3}, corresponds to $\widehat{\cal{R}}$. Then property
\eqref{Cii} implies that $\Gamma$ is bicovariant, while \eqref{Ciii}
shows that $\Gamma$ admits a *-structure. According to
Proposition~\ref{pro:314} the calculus $\Gamma$ is admissible. By
construction,
it is the {\it minimal admissible} left-covariant calculus.

Let $\L^*$ be the dual space of $\lie(G_{\!cl})$. It
turns out that $\Gamma_{\inv}$ can be naturally
embedded in $\L^*\otimes\cal{A}$. As first, let us observe that the
formula
\begin{equation}\label{C2}
\bigl(\nu\pi(a)\bigr)(X)=\nu_X\pi(a)=X(a)
\end{equation}
consistently defines a surjective linear map
$\nu\colon\Gamma_{\inv}\rightarrow\L^*$. Now, according to the
definition of $\widehat{\cal{R}}$, a linear map
$\rho\colon\Gamma_{\inv}\rightarrow\cal{L}^*\otimes\cal{A}$ given by
\begin{equation}\label{C3}
\rho=(\nu\otimes\id)\adj
\end{equation}
is injective.
\begin{lem}\label{lem:C2} The following identities hold
\begin{gather}
(\id\otimes\phi)\rho=(\rho\otimes\id)\adj\label{C4}\\
\rho(\vartheta\circ
a)=\sum_k\varphi_k\otimes\k(a^{(1)})c_ka^{(2)},\label{C5}
\end{gather}
where $\Sum_k\varphi_k\otimes c_k=\rho(\vartheta)$.
\end{lem}
\begin{pf}
Property \eqref{C4} is a direct consequence of the definition of $\rho$,
and the comodule structure property of $\adj$. Equality \eqref{C5}
follows from Lemma~\ref{lem:B8} and the following equation
\begin{equation}\label{C6}
\nu(\vartheta\circ a)=\e(a)\nu(\vartheta),
\end{equation}
which easily follows from \eqref{A5}, \eqref{B40} and \eqref{C2}.
\end{pf}

In the following, $\L^*$ will be endowed with the natural
*-structure, induced from $\lie(G_{\!cl})$. Then maps $\nu$ and $\rho$ are
hermitian.

Let $(\,)_{\!cl}$ be a scalar product in $\cal{L}^*$, with respect to
which the *-operation is antiunitar. Let $h\colon\cal{A}\rightarrow
\Bbb{C}$ be the Haar measure \cite{W2} of $G$. The formula
\begin{equation}\label{C10}
{<}\varphi\otimes a,\psi\otimes b{>}=(\varphi,\psi)_{\!cl}h(a^*b)
\end{equation}
defines a scalar product in $\L^*\otimes\cal{A}$. This enables us
to introduce a scalar product $<\,>$ in $\Gamma_{\inv}$, by requiring
that $\rho$ is isometrical.
\begin{lem}\label{lem:C4}
The introduced scalar product is $\adj$-invariant.\qed
\end{lem}

The above statement follows from the invariance of $h$.
Let $\varkappa\colon\Gamma_{\inv}\rightarrow\Gamma_{\inv}$ be a linear map
defined by
\begin{equation}\label{C12}
\varkappa\pi(a)=\pi\bigl(\k^2(a)\bigr).
\end{equation}
Consistency of this formula is a consequence of the bicovariance of
$\Gamma$.
The following identities hold
\begin{gather*}
\nu\varkappa(\vartheta)=\nu(\vartheta)\qquad
\varkappa(\vartheta)^*=\varkappa^{-1}(\vartheta^*)\qquad
\adj\varkappa=(\varkappa\otimes\k^2)\adj\\
(\vartheta,\varkappa(\eta))=(\varkappa(\vartheta),\eta)\qquad
(\vartheta^*,\eta^*)=(\varkappa^{-1}(\eta),\vartheta).
\end{gather*}

The scalar product on $\Gamma_{\inv}$ can be naturally extended to a
scalar product on $\Gamma_{\inv}^\otimes$, by tensoring and taking the
direct sum. Let us assume that the maps $\varkappa$ and $*$ are
extended from $\Gamma_{\inv}$ to $\Gamma_{\inv}^\otimes$ by requiring
multiplicativity and graded-antimultiplicativity respectively. Such
extended maps, together with the adjoint action $\adj^\otimes$ satisfy
the same relations as initial maps.

Let us assume that the ideal $S_{\inv}^\wedge$ can be orthocomplemented in
$\Gamma_{\inv}^\otimes$, relative to the constructed scalar product. Then
the space $\Gamma_{\inv}^\wedge$ is naturally realizable as the
orthocomplement of $S_{\inv}^\wedge$. In particular, we can introduce an
embedded differential map $\delta\colon\Gamma_{\inv}\rightarrow
\Gamma_{\inv}\otimes\Gamma_{\inv}$. The space
$\Gamma_{\inv}^\wedge=S_{\inv}^\perp$ is invariant under $\adj$, $*$ and
$\varkappa$.

Let $c^\top\colon\Gamma_{\inv}\rightarrow
\Gamma_{\inv}\otimes\Gamma_{\inv}$ be the ``transposed Lie commutator''
map \cite{W3}. This map can be defined by
\begin{equation}\label{C15}
c^\top=(\id\otimes\pi)\adj.
\end{equation}
Maps $\delta$ and $c^\top$ are both covariant with respect to the
adjoint action of $G$. In other words
\begin{lem}\label{lem:C6}
The following identities hold
\begin{equation}
(\delta\otimes\id)\adj=\adj^{\otimes 2}\delta\qquad
(c^\top\otimes\id)\adj=\adj^{\otimes 2}c^\top.\label{C16}
\end{equation}
\end{lem}
\begin{pf}
Applying \eqref{C15} and \eqref{B39} we obtain
\begin{multline*}
\adj^{\otimes 2} c^\top(\vartheta)=\adj^{\otimes 2}\Bigl(
\sum_k\vartheta_k\otimes\pi(c_k)\Bigr)=\sum_k\vartheta_k
\otimes\pi(c_k^{(3)})\otimes c_k^{(1)}\k(c_k^{(2)})c_k^{(4)}\\
=\sum_k\vartheta_k\otimes\pi(c_k^{(1)})\otimes c_k^{(2)}=
(c^\top\otimes\id)\adj(\vartheta),
\end{multline*}
where $\Sum_k\vartheta_k\otimes c_k=\adj(\vartheta)$. The second equality
follows from the covariance of the differential
$d\colon\Gamma^\wedge_{\inv}\rightarrow\Gamma_{\inv}^\wedge$.
\end{pf}
\begin{lem}\label{lem:C7}
\bla{i} For each $\vartheta\in\Gamma_{\inv}$ there exists $a\in\ker(\e)$ such
that
\begin{equation}\label{C17}
\begin{aligned}
\vartheta&=\pi(a)\\
\delta(\vartheta)&=-\pi(a^{(1)})\otimes\pi(a^{(2)}).
\end{aligned}
\end{equation}

\bla{ii} The following identity holds
\begin{equation}\label{C18}
c^\top=\sigma\delta-\delta.
\end{equation}
\end{lem}
\begin{pf} Let us choose $c\in\ker(\e)$ such that $\pi(c)=\vartheta$.
According to Lemma~\ref{lem:B10} we have
$d\vartheta=-\pi(c^{(1)})\pi(c^{(2)})$. According to Lemma~\ref{lem:B11}
there exists $b\in\widehat{\cal{R}}$ such that
$$\delta(\vartheta)=-\pi(c^{(1)})\otimes\pi(c^{(2)})
-\pi(b^{(1)})\otimes\pi(b^{(2)}).$$
Now $a=b+c$ satisfies \eqref{C17}.

To prove \eqref{C18} let us choose, for a given
$\vartheta\in\Gamma_{\inv}$, an element $a\in\ker(\e)$ as above. Applying
\eqref{B39}, \eqref{B40} and \eqref{C15} we obtain
\begin{equation*}
\begin{split}
-\sigma\delta(\vartheta)&=
\sigma\bigl(\pi(a^{(1)})\otimes\pi(a^{(2)})\bigr)\\
&=\pi(a^{(3)})\otimes\pi(a^{(1)})\circ\bigl(\k(a^{(2)})a^{(4)}\bigr)\\
&=\pi(a^{(3)})\otimes\pi\bigl[\bigl(a^{(1)}-\e(a^{(1)})1\bigr)\k(a^{(2)})
a^{(4)}\bigr]\\
&=\pi(a^{(1)})\otimes\pi(a^{(2)})
-\pi(a^{(2)})\otimes\pi\bigl(\k(a^{(1)})a^{(3)}\bigr)
=-c^\top(\vartheta)-\delta(\vartheta).\qed
\end{split}
\end{equation*}
\renewcommand{\qed}{}
\end{pf}
\begin{lem}\label{lem:C8} We have
\begin{equation}\label{C19}
(\nu_X\otimes\id)\delta(\vartheta)-(\id\otimes\nu_X)\delta(\vartheta)=
(\id\otimes X)\adj(\vartheta)
\end{equation}
for each $\vartheta\in\Gamma_{\inv}$ and $X\in\lie(G_{\!cl})$.\qed
\end{lem}
The following lemma gives a rough information about the ``size'' of the
space $\Gamma_{\inv}$. For each $g\in G_{\!cl}$ let $\adj^g\colon\cal{L}^*
\rightarrow\cal{L}^*$ be the induced adjoint action, given by
$$ \adj^g\nu=(\nu\otimes g)\adj.$$
\begin{lem}\label{lem:C10} \bla{i} We have
\begin{equation}\label{C21}
(\adj^g\otimes \zeta_g^\sstar)\rho=\rho
\end{equation}
for each $g\in G_{\!cl}$.

\smallskip
\bla{ii} Let $a\in\ker(\e)$ be an arbitrary $\ad$-invariant element. Then
\begin{equation}\label{C22}
a\bigl(\ker(\e)\bigr)\subseteq\widehat{\cal{R}}.
\end{equation}
\end{lem}
\begin{pf}
The statement \bla{i} directly follows form the definition of $\rho$.
Let us prove \bla{ii}. For arbitrary $b\in\ker(\e)$ and
$\ad$-invariant $a\in\ker(\e)$ we have
$$(X\otimes\id)\ad(ab)=X(ab^{(2)})\k(b^{(1)})b^{(3)}=
X(a)\e(b)1+\e(a)(X\otimes\id)\ad(b)=0. $$
This shows that $ab\in\widehat{\cal{R}}$, and hence
\eqref{C22} holds.
\end{pf}


\begin{thebibliography}{10}
\bibitem{C} Connes A {\it Non-commutative differential geometry},
Extrait des Publications Math\'ematiques--IHES {\bf 62} (1986)
\bibitem{K} Kastler D {\it Introduction to Alain Connes' non-commutative
differential geometry} Karpacz Lectures (1986)
\bibitem{KN} Kobayashi S and Nomizu K {\it Foundations of Differential
Geometry}, Interscience Publishers New York, London (1963)
\bibitem{DV} Daniel M and Viallet C M {\it The geometrical setting of
gauge theories of the Yang-Mills type}, Reviews of Modern Physics {\bf
52} 175--197 (1980)
\bibitem{P} Podle\'s P {\it Quantum Spheres}, Lett Math Phys {\bf 14}
193--202 (1987)
\bibitem{D} \Dj ur\dj evi\'c M {\it Quantum Principal Bundles and
Corresponding Gauge Theories}, Preprint, Belgrade University, Serbia
(1993)
\bibitem{W1} Woronowicz S L {\it Twisted $SU(2)$ group. An example of
a noncommutative differential calculus}, RIMS, Kyoto University 23,
117--181 (1987)
\bibitem{W2}
Woronowicz S L {\it Compact Matrix Pseudogroups}, CMP {\bf 111},
613--665 (1987)
\bibitem{W3}
Woronowicz S L {\it Differential Calculus on Compact Matrix
Pseudogroups (Quantum Groups)}, CMP {\bf 122}, 125--170 (1989)
\bibitem{W4}
Woronowicz S L {\it Tannaka-Krein Duality for Compact Matrix
Pseudogroups. Twisted $SU(n)$ groups}, Invent Math {\bf 93} 35--76
(1988)
\end{thebibliography}
\end{document}